\documentclass[11pt]{article}
\usepackage{bbm}
\usepackage{amsmath,amsthm,amssymb}
\usepackage{booktabs}
\usepackage{graphicx}
\usepackage{pgfplots}
\usepackage{float}
\usepackage{tikz}
\usepackage{tabularx}
\usetikzlibrary{calc}
\usepackage[utf8]{inputenc}
\usepackage{pgfplots}
\usepackage[colorlinks=true,urlcolor=blue,linkcolor=blue,citecolor=blue]{hyperref}
\usepackage[inline]{enumitem} 
\usepackage{amsmath}
\usepackage{amsfonts}
\usepackage{amssymb}
\usepackage[linesnumbered,ruled,lined]{algorithm2e}
\SetKwIF{IfNoEnd}{ElseIfNoEnd}{ElseNoEnd}
  {if}{ then}{ else if}{ else}{} 

\usepackage{setspace}
\usepackage{pict2e,color}
\usepackage{multirow}
\usepackage{wrapfig}
\usepackage{booktabs}
\usepackage{amsfonts}
\usepackage{verbatim} 
\usepackage{sectsty}
\usepackage{url}
\usepackage{tikz}
\usepackage{empheq}
\usepackage[normalem]{ulem}
\usepackage[titletoc,title]{appendix}
\usepackage{pgfplots}
\usepackage{algpseudocode}
\usepackage[flushleft]{threeparttable}
\usepackage{soul} 
\usepackage[numbers]{natbib}
\usepackage{microtype}
\usepackage{subcaption}
\usepackage{pbox}
\usepackage{etoolbox}
\preto\subequations{\ifhmode\unskip\fi}
\usepackage{booktabs}
\usepackage{multirow}
\usepackage{derivative}
\theoremstyle{definition}

\newtheorem{theorem}{Theorem}

\newtheorem{proposition}{Proposition}
\newtheorem{corollary}{Corollary}

\newtheorem{assumption}{Assumption}

\def\argmin{\mathop{\rm arg\,min}}%
\allowdisplaybreaks[4]

\title{Distributionally Robust Optimization for Chemotherapy Scheduling under Asymmetric and Multi-Modal Uncertainty}
\author{Qing Zhu\thanks{Department of Integrated Systems Engineering, The Ohio State University, Columbus, OH, USA, Email: {\tt zhu.2166@osu.edu};}~~~Xian Yu\thanks{Corresponding author; Department of Integrated Systems Engineering, The Ohio State University, Columbus, OH, USA, Email: {\tt yu.3610@osu.edu};}
~~~Yu-Li Huang\thanks{Robert D. and Patricia E. Kern Center for the Science of Health Care Delivery, Mayo Clinic, Rochester, Minnesota, USA, Email: {\tt huang.yuli@mayo.edu};}}
\date{}

\begin{document}
\graphicspath{{figures/}}

\maketitle              
\begin{abstract}
\textit{\textbf{Problem definition:}} We consider a real-world chemotherapy scheduling template design problem, where we cluster patient types into groups and find a representative time-slot duration for each group to accommodate all patient types assigned to that group, aiming to minimize the total expected idle time and overtime. From Mayo Clinic's real data, most patients' treatment durations are asymmetric (e.g., shorter/longer durations tend to have a longer right/left tail). 
\textit{\textbf{Methodology/results:}} Motivated by this observation, we consider a distributionally robust optimization (DRO) model under an asymmetric and multi-modal ambiguity set, where the distribution of the random treatment duration is modeled as a mixture of distributions from different patient types. The ambiguity set captures uncertainty in both the mode probabilities, modeled via a variation-distance-based set, and the distributions within each mode, characterized by moment information such as the empirical mean, variance, and semivariance.
We reformulate the DRO model as a semi-infinite program, which cannot be solved by off-the-shelf solvers. To overcome this, we derive a closed-form expression for the worst-case expected cost and establish lower and upper bounds that are positively related to the variability of patient types assigned to each group, based on which we develop exact algorithms and highly efficient clustering-based heuristics.  \textit{\textbf{Managerial implications:}} The lower and upper bounds on the worst-case cost imply that the optimal cost tends to decrease if we group patient types with similar treatment times. Through numerical experiments based on both synthetic datasets and Mayo Clinic's real data, we illustrate the effectiveness and efficiency of the proposed exact algorithms and heuristics and showcase the benefits of incorporating asymmetric information into the DRO formulation.

~
\end{abstract}
\section{Introduction}
\begin{wrapfigure}{r}{0.5\textwidth}
\centering
\begin{subfigure}[t]{0.24\textwidth}
    \centering
    \includegraphics[width=0.95\textwidth]{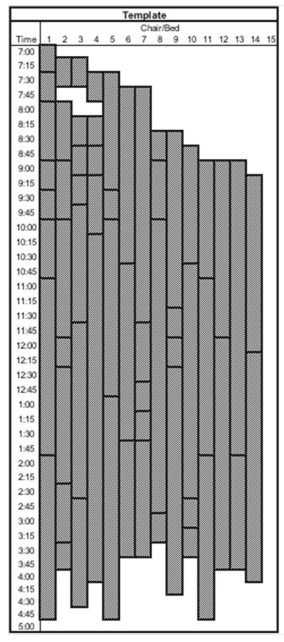} 
    \caption{Template}
    \label{fig:mayo template}
  \end{subfigure}
  \begin{subfigure}[t]{0.24\textwidth}
    \centering
    \includegraphics[width=\textwidth]{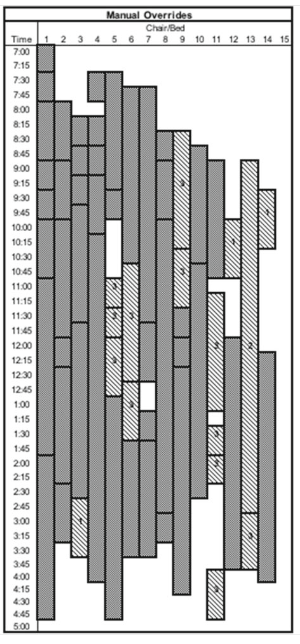} 
    \caption{Override}
    \label{fig:mayo override}
  \end{subfigure}
  \caption{Mayo Clinic example template}
  \label{fig:mayo example template}
\end{wrapfigure}

Chemotherapy is one of the most resource‐intensive services in ambulatory cancer units \cite{huang2018alternative}. Each infusion in chemotherapy treatment requires nurses to perform a coordinated set of tasks—drug verification, chair and line setup, patient education, infusion monitoring, and safe disconnection, resulting in highly uncertain treatment durations. The current scheduling procedure at the Mayo Clinic consists of two phases: (i) chemotherapy template design, where a fixed daily template (e.g., a $10$-hour schedule on each chair) is optimized up front and will be implemented in the next several months (see Figure\ \ref{fig:mayo template}), and (ii) chemotherapy patient scheduling, where patients are placed into the most appropriate available time slots as they call in to schedule their appointments sequentially. When no appropriate slot is available, the hospital operator combines two neighboring slots into a longer one or divides one longer slot into several shorter ones to fit the current patient, leading to the so-called ``overrides'' (see Figure \ref{fig:mayo override}). In this paper, we mainly focus on Phase (i), where we optimize the duration of each contributing block in the daily template. Since these time-slot durations are decided before patients arrive and will be fixed through the next implementation horizon, they should be designed to accommodate variability in patient types and stochastic treatment lengths. 

\begin{figure}[ht!]
\centering

\begin{subfigure}[t]{0.22\textwidth}
  \centering
  \includegraphics[width=\linewidth]{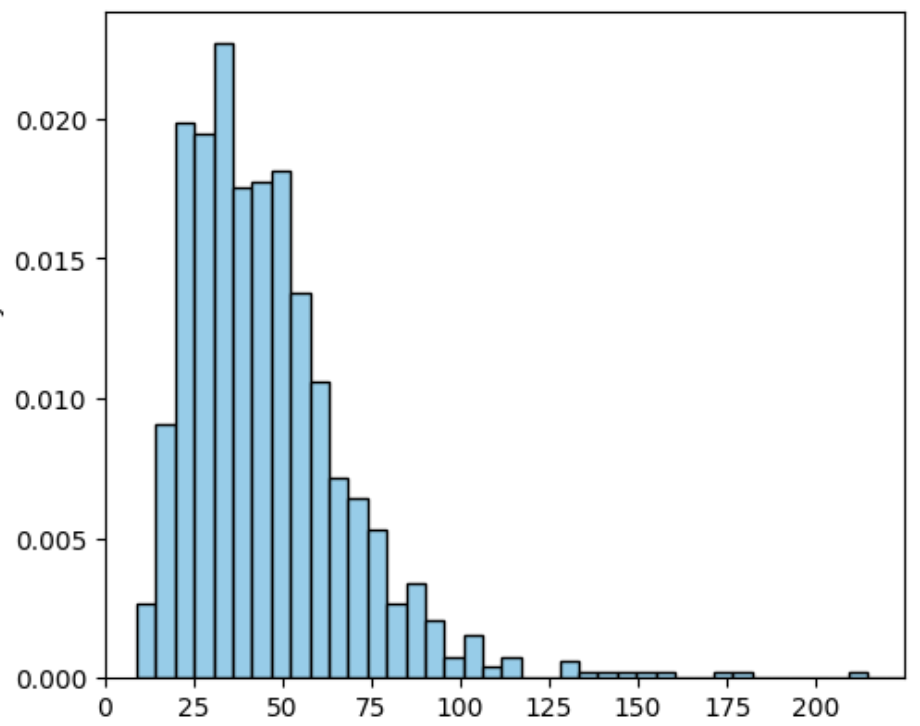}
  \subcaption{30-minute}
  \label{fig:mayo-30}
\end{subfigure}
\begin{subfigure}[t]{0.22\textwidth}
  \centering
  \includegraphics[width=\linewidth]{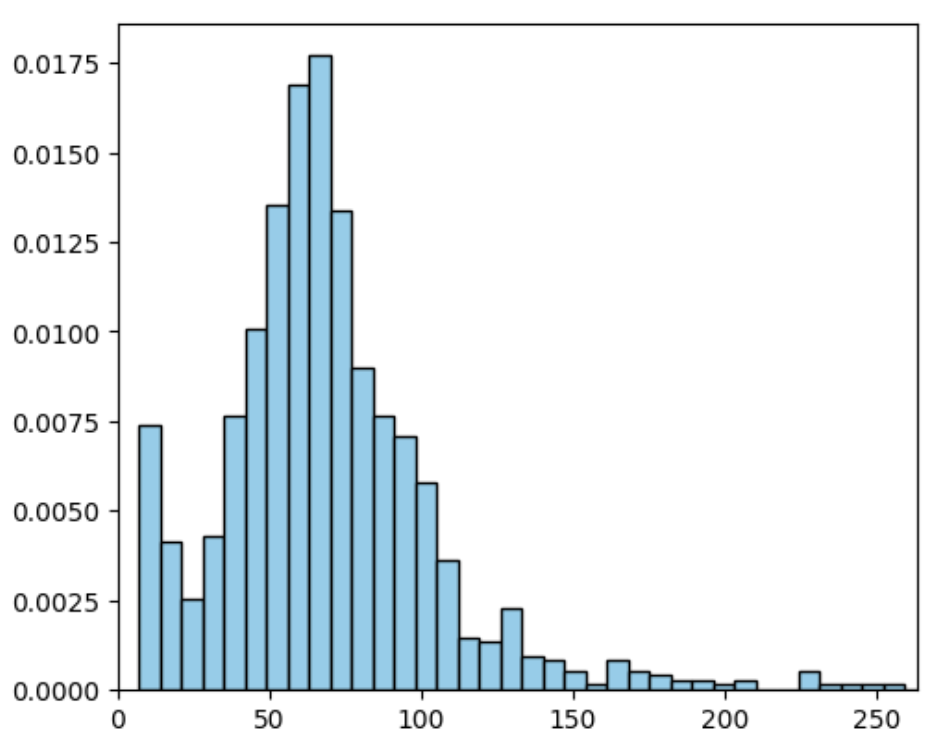}
  \subcaption{60-minute}
  \label{fig:mayo-60}
\end{subfigure}
\begin{subfigure}[t]{0.22\textwidth}
  \centering
  \includegraphics[width=\linewidth]{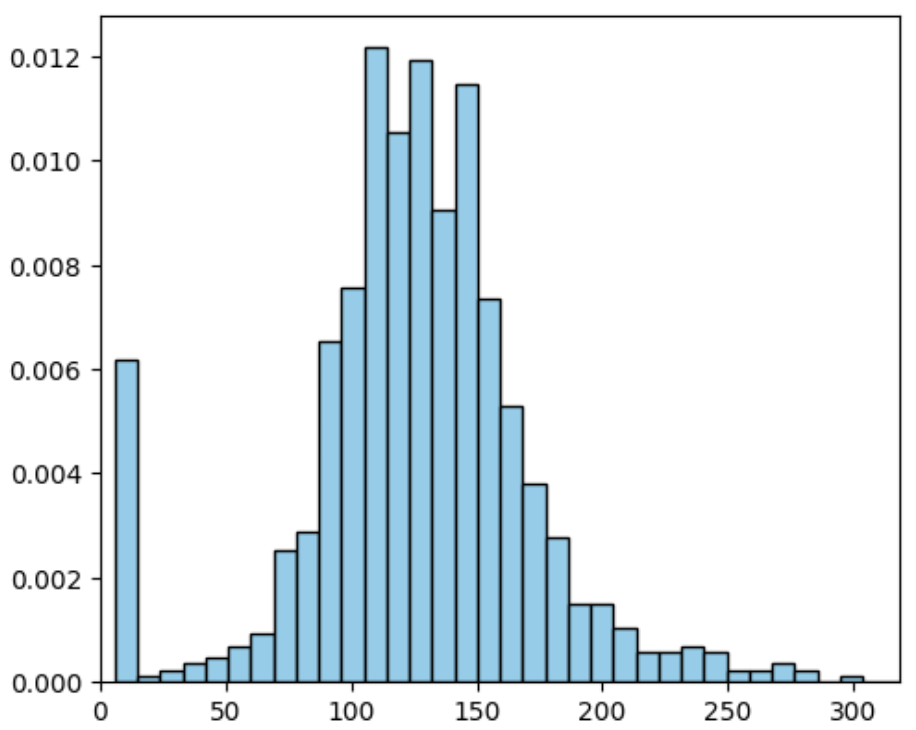}
  \subcaption{120-minute}
  \label{fig:mayo-120}
\end{subfigure}


\begin{subfigure}[t]{0.22\textwidth}
  \centering
  \includegraphics[width=\linewidth]{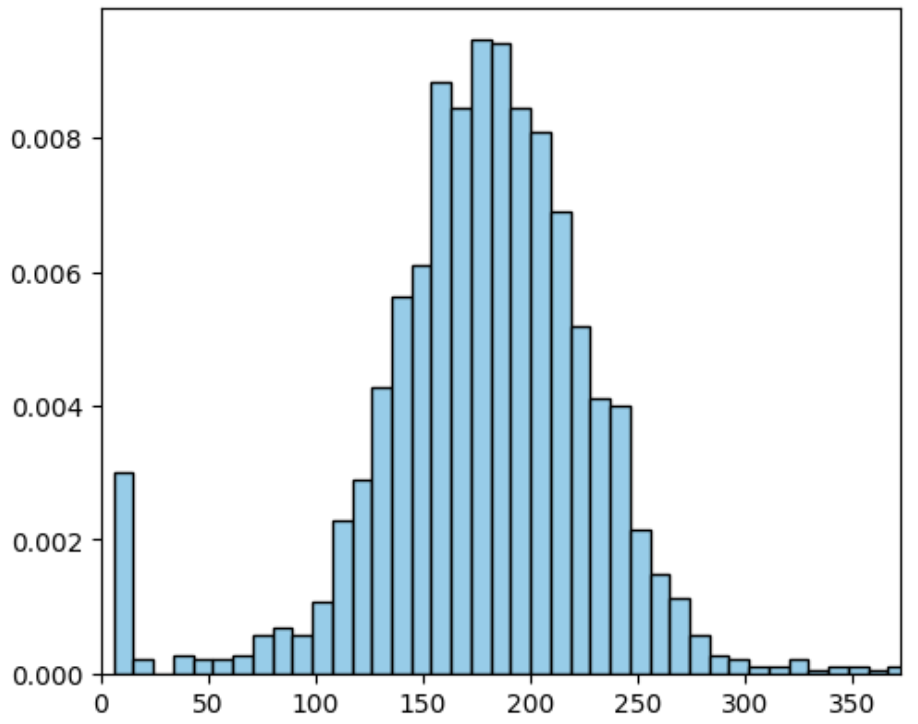}
  \subcaption{180-minute}
  \label{fig:mayo-180}
\end{subfigure}
\begin{subfigure}[t]{0.22\textwidth}
  \centering
  \includegraphics[width=\linewidth]{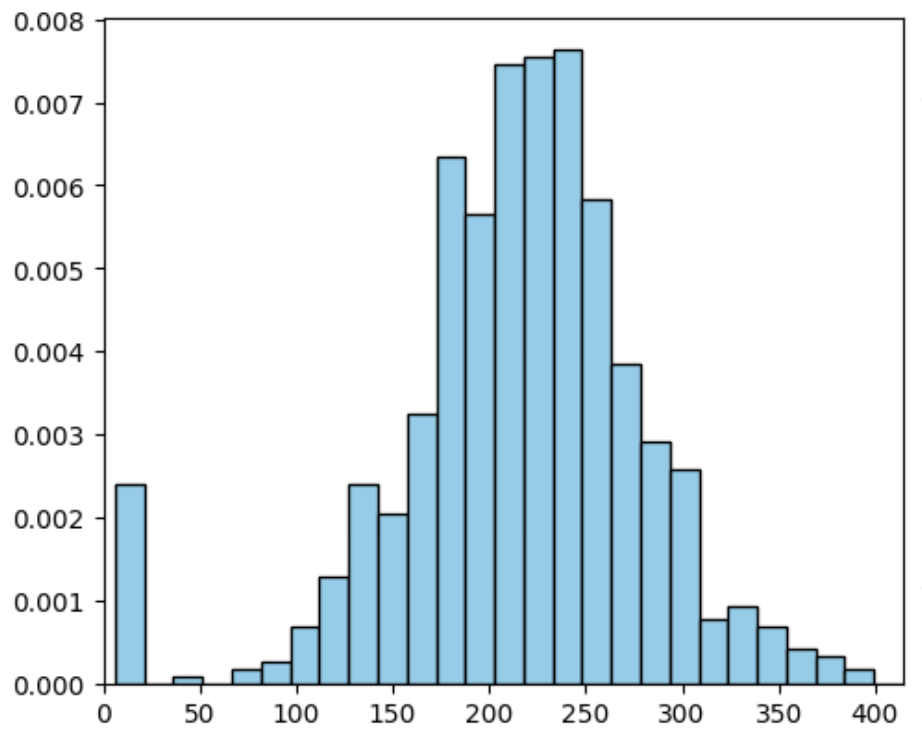}
  \subcaption{240-minute}
  \label{fig:mayo-240}
\end{subfigure}
\begin{subfigure}[t]{0.22\textwidth}
  \centering
  \includegraphics[width=\linewidth]{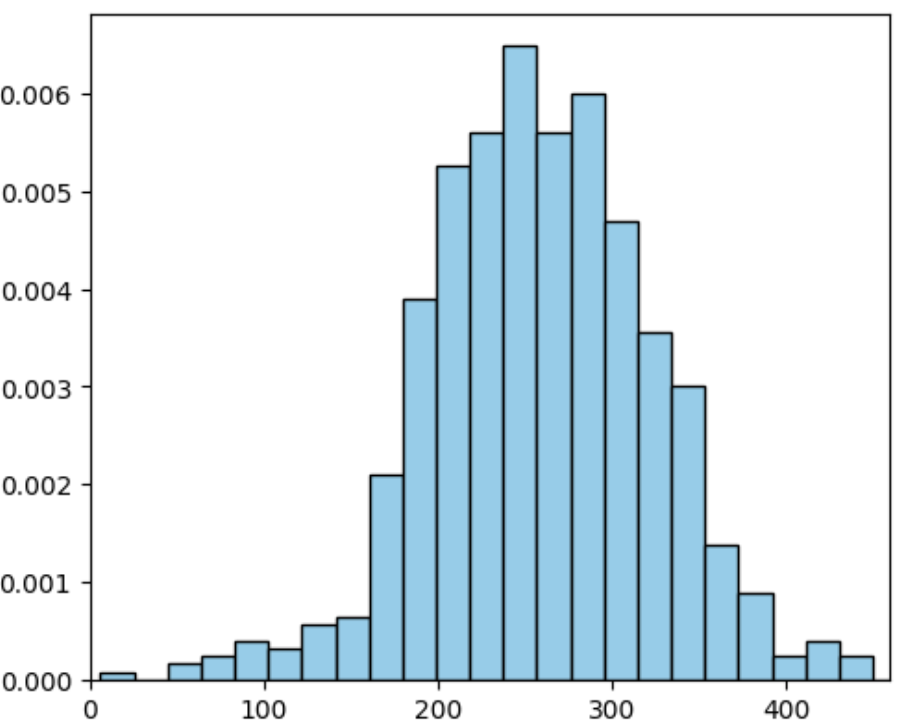}
  \subcaption{300-minute}
  \label{fig:mayo-300}
\end{subfigure}
\begin{subfigure}[t]{0.22\textwidth}
  \centering
  \includegraphics[width=\linewidth]{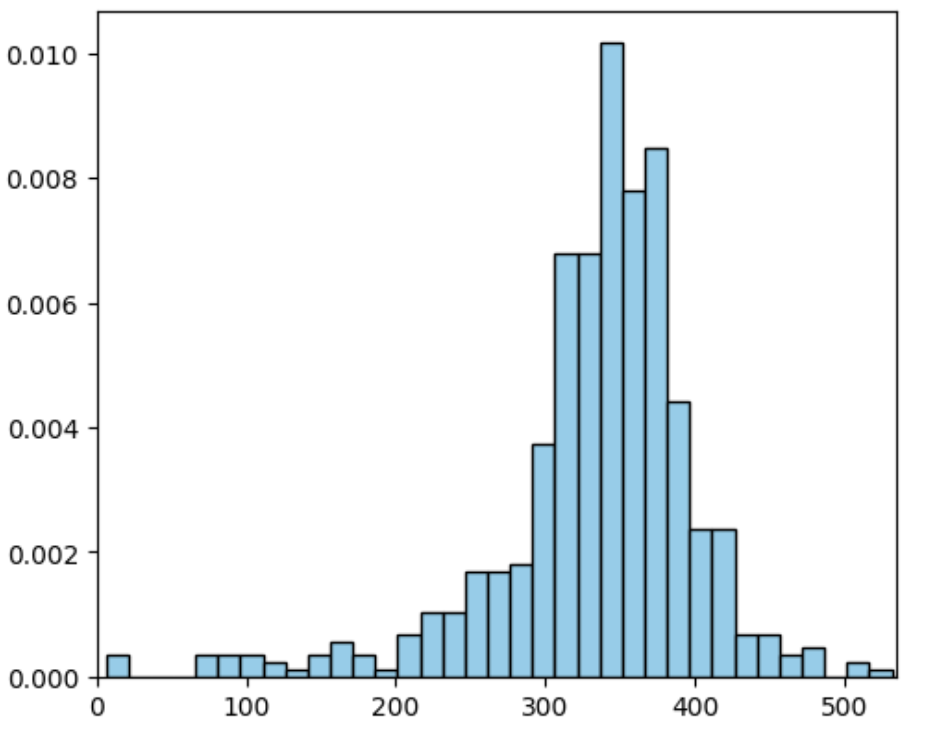}
  \subcaption{360-minute}
  \label{fig:mayo-360}
\end{subfigure}%

\caption{Histograms of treatment durations across seven patient types (Mayo Clinic's real data).}
\label{fig:mayo distribution}
\end{figure}

As demonstrated by Mayo Clinic's real data (See Figure \ref{fig:mayo distribution} and Table \ref{tab:duration-stats}), treatment durations exhibit wide variability. The hospital operator currently classifies all the patients into seven distinct types (e.g., $30$-min type, $60$-min type, etc.), based on their chemotherapy regimens, cancer subtypes, or characteristic patterns of comorbidities. Among these seven patient types, most treatment duration distributions are asymmetric.  For example, Figure \ref{fig:mayo-30} shows that $30$–minute infusions can extend substantially when adverse reactions occur, leading to a ``right-skewed'' or ``positive-skewed'' distribution. On the other hand, Figure \ref{fig:mayo-360} shows that $360$–minute sessions sometimes finish earlier than planned, leading to a ``left-skewed'' or ``negative-skewed'' distribution. This asymmetric property of treatment duration distributions is further confirmed in Table \ref{tab:duration-stats}, where we report the mean, standard deviation, and semivariance for each patient type. Note that a positive semivariance corresponds to a right-skewed distribution while a negative semivariance corresponds to a left-skewed distribution.
This feature creates unique modeling challenges: shorter-than-expected sessions often leave nurses and medical resources idle, while longer-than-expected sessions can generate significant overtime. Both cases directly impact patient care quality and operational efficiency, which emphasizes the importance of explicitly capturing the skewness of treatment durations. 


\captionsetup[sub]{justification=centering,singlelinecheck=false}

\begin{table}[ht]
  \centering
  \caption{Statistics of treatment durations across seven patient types (Mayo Clinic's real data)}
 \begin{tabular}{lrrrrr}
\toprule
\textbf{Patient Type} &\textbf{Count}& \textbf{Proportion} &\textbf{Mean} & \textbf{Std} & \textbf{Semivariance} \\
\midrule
$30$-min  & $923$& $13.83\%$ & $48.70$  & $31.15$ & $0.59$  \\
$60$-min  & $1013$& $15.18\%$ & $73.24$  & $39.65$& $0.47$ \\
$120$-min & $913$ & $13.69\%$ & $133.88$ & $42.98$ & $0.33$ \\
$180$-min & $1851$ & $27.74\%$ & $182.87$ & $47.20$ & $0.07$  \\
$240$-min & $745$ & $11.17\%$ & $223.34$ & $57.07$ & $0.07$  \\
$300$-min & $641$ & $9.61\%$ & $258.78$ & $65.27$ & $-0.005$ \\
$360$-min & $586$ & $8.78\%$ & $336.29$ & $66.66$ & $-0.25$ \\
\bottomrule
\end{tabular}
\label{tab:duration-stats}
\end{table}


However, the current practice at the Mayo Clinic adopts a deterministic approach that ignores the randomness and skewness of the treatment durations. They use a specific time-slot duration to accommodate one patient type (e.g., $30$-min time slots designed solely for $30$-min patient type). The issues of such a deterministic approach are two-fold. First, once a schedule template is fixed (see Figure\ \ref{fig:mayo template}), we have a certain number of time slots for each patient type (e.g., 23 $30$-min time slots, $8$ $60$-min time slots, $10$ $120$-min time slots, etc.). Given the stochastic nature of the patient counts in each type, if the realized count of each patient type exceeds the allocated number of time slots, this may result in excessively large overrides. Second, the deterministic approach fails to capture the asymmetry of the random treatment durations, as reflected in Figure \ref{fig:mayo distribution}, which may result in substantial idle time and overtime in practice. 
Overall, ignoring the randomness and skewness of the treatment durations when designing templates could potentially lead to significant inefficiencies in resource utilization and compromised patient care quality.

In this paper, we adopt a stochastic approach that treats the treatment duration for each patient type as a random variable and explicitly captures its asymmetric information by using second-order statistics (e.g., semivariance). We propose a general framework to consolidate all patient types into a smaller number of groups, designing a representative time-slot duration for each group, such that each duration can accommodate all patient types assigned to that group.
We formulate this chemotherapy patient grouping and time-slot duration design problem as a \underline{M}ulti\underline{m}odal \underline{A}symmetric \underline{D}istributionally \underline{R}obust \underline{O}ptimization (MMA-DRO) model that accounts for (i) multimodal uncertainty, by representing each patient group's treatment time distribution as a mixture of the distributions of the patient types assigned to that group, and (ii) asymmetric information, by incorporating semivariance information into the ambiguity set for each group.
Specifically, the model determines patient-type-to-grouping assignments and the time-slot duration within each group to minimize the group activation cost and the worst-case expected cost of idle time and overtime. Given the assignment decisions, each group finds the optimal time-slot duration against the worst-case scenario within a multi-modal and asymmetric ambiguity set. This data-driven formulation allows us to create more flexible scheduling templates by consolidating multiple patient types and assigning each duration to accommodate the skewness of treatment time distributions from several patient types, potentially reducing idle time, overtime, and overrides in real-world practice.
\begin{figure}[ht]
    \centering
    \includegraphics[width=0.9\linewidth]{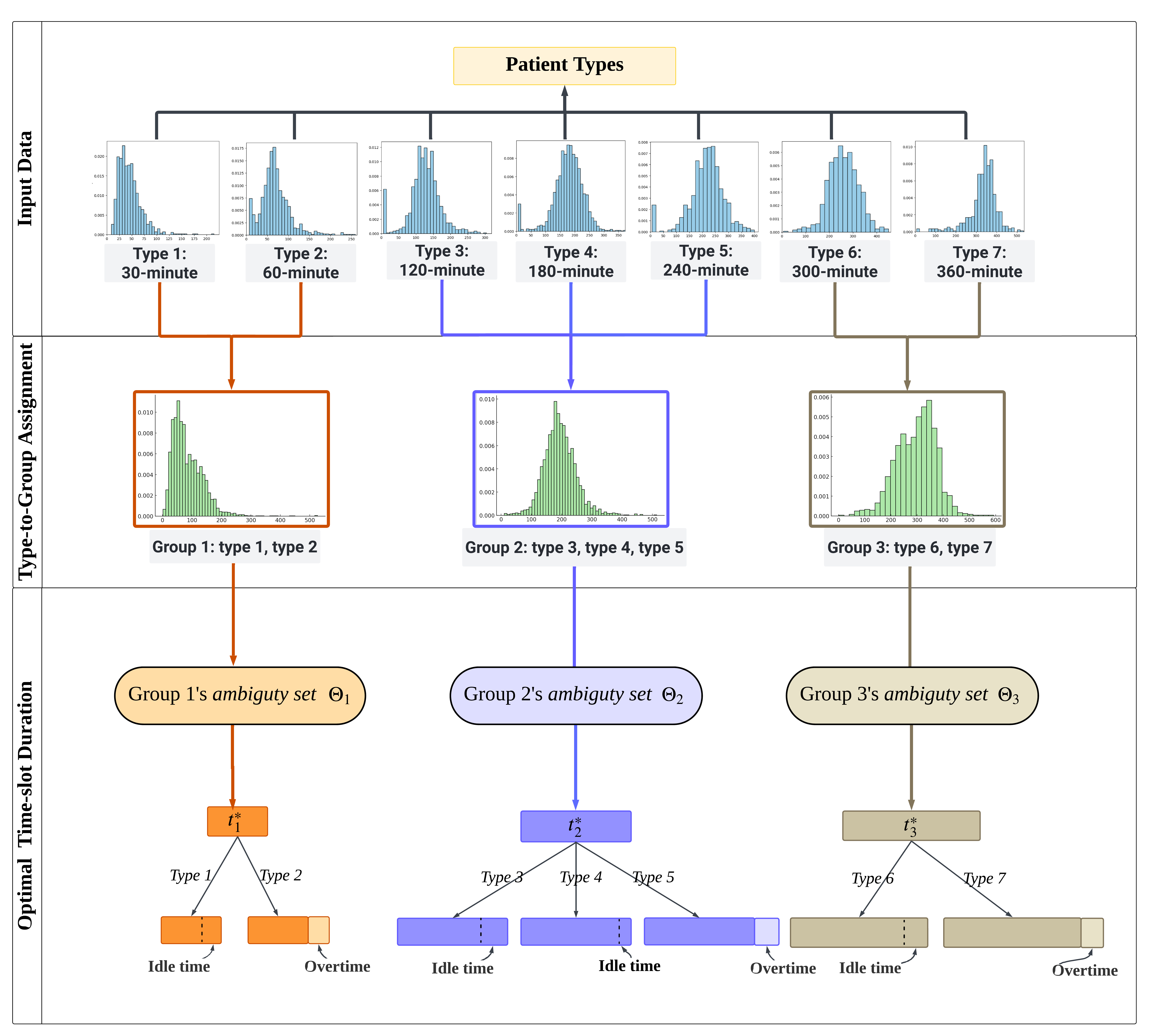}
    \caption{A flowchart of the MMA-DRO model, where $\Theta_1,\ \Theta_2,\ \Theta_3$ denote the multimodal asymmetric ambiguity sets of the uncertain chemotherapy durations and $t_1^*,\ t_2^*,\ t_3^*$ denote the optimal time-slot durations for groups 1, 2, 3, respectively. }
    \label{fig:intro}
\end{figure}


The resulting MMA-DRO model is highly non-convex, as our first-stage assignment decisions can affect the ambiguity sets in a complicated form, leading to decision-dependent ambiguity sets. Using strong duality, we recast the MMA-DRO model as a monolithic minimization problem; however, since we consider a continuous support set for the uncertainty, this results in a semi-infinite program. To overcome the intractability of having infinitely many constraints, we first derive a closed-form expression for the worst-case cost of idle time and overtime given any assignment decision, using which we design exact algorithms to solve this semi-infinite program by enumerating all possible assignments. We also establish both lower and upper bounds on each group's optimal cost, which are positively related to the variability of patient types assigned to that group. This managerial insight further motivates us to leverage clustering approaches (e.g., K-Means and K-Medoids algorithms) to find a near-optimal assignment decision first and then embed it into the closed-form expression to optimize the time-slot durations in a quick fashion.

The remainder of the paper is organized as follows. Section \ref{sec:related_work} reviews the relevant DRO and healthcare scheduling literature, highlighting the existing gaps in modeling multimodality and asymmetry under the DRO framework. Section \ref{sec:problem section} presents our problem formulation and derives a monolithic semi-infinite program using strong duality. Section \ref{sec:exact} introduces exact algorithms to solve the resulting semi-infinite program by using a closed-form expression of the worst-case expected cost of idle time and over time and establishes the lower and upper bounds on each group's optimal cost. Motivated by this insight, Section \ref{sec: heuristic} develops highly efficient heuristics to obtain near-optimal solutions by leveraging clustering approaches (e.g., K-Means and K-Medoids algorithms). In Section\ \ref{sec: synthetic section}, we report experimental results based on synthetic datasets to conduct sensitivity analyses and out-of-sample tests, and in Section\ \ref{sec: mayo section}, we present a case study based on Mayo Clinic's real chemotherapy data. Finally, Section \ref{sec:conclusions} discusses managerial insights and directions for future research.

\subsection{Related Work}\label{sec:related_work}
\paragraph{Chemotherapy scheduling under uncertainty.} We first review the most relevant literature that applies optimization under uncertainty techniques to solve the chemotherapy scheduling problem.  Stochastic programming (SP) emerges as one of the most widely employed frameworks in optimization under uncertainty, which aims to optimize the system performance on average or based on a risk measure. For instance, \cite{castaing2016stochastic,demir2021stochastic,alvarado2018chemotherapy,karakaya2023stochastic,gul2021chemotherapy} proposed two-stage stochastic mixed-integer programs to optimize patient appointment schedules under uncertain treatment durations. Among these papers, decisions include the sequence of patients on a daily schedule, their appointment starting times, patient-to-nurse/chair assignments, and clinic resource schedules. The objective is often set to minimize a trade-off between the expected cost of patient waiting time, chair idle time, and nurse overtime. To solve these problems, they developed both exact and heuristic methods, based on progressive hedging \cite{demir2021stochastic}, scenario bundling-based decomposition \cite{karakaya2023stochastic,gul2021chemotherapy}, etc. Robust optimization (RO), on the other hand, aims to find the optimal decisions against the worst-case scenario within a pre-defined uncertainty set, without access to the probability distributions. Numerous studies have shown its suitability to handle uncertain treatment durations and operational disruptions in chemotherapy scheduling \cite{issabakhsh2021scheduling,behnamian2023multi,bauerhenne2024robust}. \cite{bauerhenne2024robust} and \cite{behnamian2023multi} considered box uncertainty sets to model random appointment durations and no-shows, and service quality, respectively, where \cite{issabakhsh2021scheduling} proposed a Robust Slack Allocation model by inserting a slack variable to each uncertain constraint.
To solve the RO models, \cite{issabakhsh2021scheduling} developed a heuristic based on adaptive large neighborhood search, \cite{behnamian2023multi} proposed a multi-objective particle swarm optimization heuristic, and \cite{bauerhenne2024robust} explored the structural property and reformulated the RO model as a compact mixed-integer linear program. More recently, DRO has become a powerful framework that bridges SP and RO by optimizing decisions against the worst-case expected performance over an ambiguity set of probability distributions. Recent advances have extended DRO to address chemotherapy scheduling challenges. \cite{hristov2025distributionally} developed a DRO framework for radiation therapy treatment planning that accounts for inter-fraction geometric uncertainty of the clinical target volume. \cite{bansal2021distributionally} adopted a DRO model to solve a scheduling problem that requires coordination between clinical and surgical appointments, where they considered two sources of uncertainty: (i) the patient may or may not need surgery after the clinical appointment and (ii) the surgery duration for each patient and procedure is unknown. \cite{keyvanshokooh2022coordinated} proposed an integrated multi-stage stochastic and DRO framework to jointly schedule clinical and surgical appointments under uncertainties in appointment demand, surgery necessity, and surgery durations to minimize overtime. \cite{shehadeh2020distributionally} proposed a DRO framework to determine optimal appointment sequences and minimize the worst-case sum of patient waiting, provider idling, and overtime. Different than the aforementioned papers, which optimize patient appointment start times at any minute, our paper considers a fixed template with pre-determined time slots, which can fit the patients into the most appropriate time slots as they call in to schedule. The use of pre-determined time slots facilitates both the construction of a scheduling template and the real-time assignment of patients to time slots. To the best of our knowledge, only \cite{huang2019chemotherapy} studied a chemotherapy scheduling template design problem; however, they considered a deterministic setting and the time-slot durations are fixed at 30 minutes, 60 minutes, 120 minutes, etc. Instead, the main focus of this paper is to optimize the durations of the contributing blocks in such a template, while considering the asymmetric property of the treatment time distributions.

\paragraph{Modeling asymmetric uncertainty.} 
To model asymmetric uncertainty, \cite{chen2007robust} introduced a robust optimization approach for constructing asymmetric uncertainty sets using forward and backward deviation measures of random variables. \cite{natarajan2008incorporating} then extended it to incorporate asymmetric distributional information in a modified Value-at-Risk (VaR) measure, named Asymmetry-Robust VaR (ARVaR), and proposed a computationally tractable approximation method for minimizing the ARVaR of a portfolio based on robust optimization techniques.
In DRO, traditional approaches to constructing ambiguity sets generally fall into two categories: (i) moment-based ambiguity sets \cite{delage2010distributionally,mehrotra2014cutting, wagner2008stochastic, zhang2018ambiguous,yu2022multistage}, which are defined using information about the moments of the distribution (e.g., mean, variance, or higher-order moments), and (ii) distance-based ambiguity sets \cite{esfahani2018data, blanchet2019quantifying, gao2023distributionally,jiang2018risk,jiang2016data,bayraksan2015data}, which are formed by measuring the distance between the empirical distribution and alternative distributions under a chosen metric (e.g., Wasserstein distance or $\phi$-divergence).
These DRO models, however, often fail to account for the inherent skewness or asymmetric tail behavior of the uncertainty, which limits their effectiveness in modeling the treatment time distributions. 
Several recent works have addressed this gap by directly integrating asymmetry information into the ambiguity sets. \cite{natarajan2018asymmetry} constructed the ambiguity set with mean, variance, and semivariance information to account for the asymmetric demand distribution when solving a newsvendor problem. \cite{chen2024robust} introduced a general moment-dispersion ambiguity set by incorporating a dispersion characteristic function, which can capture the uncertainty's complex attributes, such as sub-Gaussian and asymmetric dispersion. 
\cite{chen2023distributionally} developed a distributionally robust pricing model to optimize the worst-case profit, where they used semivariance to capture the asymmetric information of the consumers' valuation distribution. 

\paragraph{Modeling multimodal uncertainty.} We now turn our attention to recent DRO studies that have emphasized modeling multimodal uncertainty, which is characterized by distributions composed of distinct probability modes and their associated distributional information. \cite{hanasusanto2015distributionally} addressed a risk-averse, multi-dimensional newsvendor problem through a DRO formulation that captures multimodal demand distributions represented by spatially separated probability clusters. Building on this direction, \cite{yu2024distributionallyrobustoptimizationmultimodal} introduced a general two-stage DRO model under multimodal decision-dependent ambiguity sets, where the first-stage decisions could affect both the mode probabilities and the corresponding distributions in each mode. To model bimodal demand distributions in facility location problems, \cite{shehadeh2021distributionally} constructed a scenario-based DRO model that captures the ambiguity of the bimodal demand to optimize the number and location of the facilities.
Moreover, \cite{shehadeh2020distributionally} considered an outpatient colonoscopy scheduling problem and modeled the colonoscopy durations as bimodal distributions in a DRO framework to determine optimal appointment sequences, minimizing the worst-case expected cost of patient waiting, provider idling, and overtime.  In a complementary line of research within machine learning, \cite{sagawa2019distributionally} proposed a group DRO approach to minimize the worst-case expected loss of predictive models over mixtures of predefined distributions, which effectively captured potential multimodalities present in the test data distribution.

To the best of our knowledge, no existing studies have addressed the chemotherapy scheduling template design problem while considering multimodal and asymmetric uncertainty. This paper aims to bridge this gap by introducing an MMA-DRO framework that groups multiple patient types and optimizes the time-slot duration within each group to accommodate the skewness of treatment time distributions assigned to this group. 

\subsection{Contributions}
We summarize the main contributions of this paper as follows:
\begin{enumerate}
\item We formulate the chemotherapy patient grouping and time-slot duration design problem as a distributionally robust optimization (DRO) model, under asymmetric and multi-modal ambiguity sets. Specifically, we optimize the assignment of different patient types to groups and decide the duration of time slots for each group to accommodate various patient types assigned to this group.

\item Under a continuous support set of the uncertainty, we reformulate the DRO model as a monolithic semi-infinite program, which cannot be solved via standard optimization solvers. To overcome this, we derive a closed-form expression for the inner maximization problem, based on which we develop two exact algorithms to solve the DRO model by enumerating all possible patient type-to-group assignment decisions.

\item To obtain managerial insight into optimal grouping assignments, we provide lower and upper bounds for each group's optimal cost, which both depend on the difference between the maximum and minimum average treatment time among all the patient types assigned to this group. This implies that the optimal cost tends to decrease if we group patient types with similar treatment times.

\item Motivated by the above insights, we propose to use K-means clustering algorithms to find a near-optimal assignment based on the similarity of different statistics across all patient types. Once the assignment is given, we leverage the closed-form expression of the inner maximization problem again to find the optimal appointment duration for each group.

\item Finally, we validate our approach through numerical experiments on both real-world Mayo Clinic data and synthetic datasets.. The results demonstrate that explicitly modeling multimodality and asymmetry within the ambiguity set significantly improves the robustness and efficiency of chemotherapy treatment–template design.

\end{enumerate}

\section{Problem Formulations}\label{sec:problem section}
Let $\tilde{\xi} \in \Xi$ represent the random chemotherapy treatment duration that follows an underlying (unknown) distribution $\mathcal{P}$ with a compact support set $\Xi$ (in this paper, we consider a bounded and continuous support set $\Xi$). Let $\mathcal{L}=\{1,2,\ldots,L\}$ denote a set of patient types depending on specific chemotherapy regimens or cancer types (e.g., $30$-min type, $60$-min type, etc.), and $\mathcal{G}=\{1,2,\ldots,G\}$ denote a set of potential patient groups. Note that $G$ can be set to $L$ to represent the maximum number of groups one can create, or it can be set to a lower number depending on the operational constraints and the decision maker's preference. Define $x_{lg}$ as a binary decision variable indicating if the patient type $l \in \mathcal{L}$ is assigned to the group $g \in \mathcal{G}$, and $y_g$ as a binary decision variable indicating whether the group $g \in \mathcal{G}$ is activated or not. Let $t_g \in \mathcal{T}:=[0,T]$ be the time-slot duration to be determined for each group $g \in \mathcal{G}$, which will serve for all the patient types assigned to this group. The upper bound $T$ for the time-slot duration can be set to the maximum treatment time or set according to some operational constraints (e.g., 5 hours due to limited nurse resources). Throughout the paper, we use $a^+$ to denote $\max\{0,a\}$ and use bold symbols to denote vectors.

Let $q$, $b$, and $c_g$ represent the unit cost of overtime, idle time, and activating group $g\in\mathcal{G}$, respectively. Specifically, $q$ is the unit penalty on the overtime when actual chemotherapy durations exceed their scheduled time-slot durations (i.e., $\tilde{\xi} \ge t_g$), and $b$ is the unit penalty on the idle time when treatments conclude earlier than scheduled (i.e., $t_g \ge \tilde{\xi}$). The unit group-activating cost $c_g$ can be set the same across all groups $g$, so that the term $\sum_g c_g y_g$ penalizes the total number of groups. When $c_g$ is set to be a larger value, it encourages a smaller number of activated groups so that each group will accommodate more patient types. This creates a natural trade-off between the number of activated groups and the overall system performance.

Our goal is to simultaneously determine the optimal assignment of patient types to groups and the corresponding time-slot duration within each group in order to minimize the overall cost incurred by overtime, idle time, and group activation. 
To formulate this problem, we adopt the following MMA-DRO framework:
\begin{subequations}\label{eq:MMA-DRO}
\begin{align}
    \min_{x_{lg}, y_g, t_g} \quad & \left\{\sum_{g \in \mathcal{G}} c_g y_g +\frac{1}{\sum_{g \in \mathcal{G}} y_g}  \sum_{g \in \mathcal{G}}\max_{P \in \Theta_g(\boldsymbol{x},\boldsymbol{y})} \mathbb{E}_{\tilde{\xi}\sim P} \left[q(\tilde{\xi}-t_g)^+ + b(t_g-\tilde{\xi})^+\right]\right\} \label{eq:MMA-DRO objective} \\
    \text{s.t.} \quad & \sum_{g \in \mathcal{G}} x_{lg}=1,\ \forall l \in \mathcal{L}\label{eq:sum to one}\\
    \quad & x_{lg} \leq y_g,\ \forall l \in \mathcal{L},\ g \in \mathcal{G}\label{eq:x less than y}\\
    \quad & \sum_{l=1}^L x_{lg} \geq y_g,\ \forall g \in \mathcal{G}\label{eq:feasibility} \\
    \quad & 0\le t_g\le T,\ \forall g \in \mathcal{G}\label{eq:t range},
\end{align}
\end{subequations}
where we aim to minimize the fixed cost of activating groups together with the worst-case expected cost of idle time and overtime, averaged over the activated groups in Eq.\ \eqref{eq:MMA-DRO objective}. Constraints\ \eqref{eq:sum to one} ensure that each patient type is assigned to exactly one group. Constraints\ \eqref{eq:x less than y} ensure that patient types can only be assigned to groups that are activated. Constraints\ \eqref{eq:feasibility} enforce that every activated group must have at least one patient type assigned to it, which eliminates empty groups. Constraints\ \eqref{eq:t range} put a lower and upper bound on the time-slot duration $t_g$ for each $g\in \mathcal{G}$.

Given the group activation and assignment decisions $(\boldsymbol{x},\boldsymbol{y})$, within each group $g$, we consider a multi-modal ambiguity set $\Theta_g(\boldsymbol{x},\boldsymbol{y})$ to capture the uncertainty in treatment-time distributions, where each patient type in the group is treated as a possible mode $l$, weighted by its proportion within the group, referred to as the mode probability $p_l$. Following\ \cite{yu2024distributionallyrobustoptimizationmultimodal}, we adopt a two-layer ambiguity set to describe $\Theta_g(\boldsymbol{x},\boldsymbol{y})$, where the first layer captures the uncertainty in mode probabilities and the second layer represents the uncertainty around the distribution in each mode. Specifically, for each $g\in \mathcal{G}$, any candidate distribution $ P \in \Theta_g(\boldsymbol{x},\boldsymbol{y})$ can be represented as a mixture of distributions from each mode weighted by the mode probability as follows:
\begin{align}\label{eq:ambiguity set}
\Theta_g(\boldsymbol{x},\boldsymbol{y})= \left\{\,P=\sum_{l=1}^L p_{l}\, f_{l} :
\underbrace{\boldsymbol{p}\in \Delta_g\bigl(\hat{\boldsymbol{p}}, \boldsymbol{x},\boldsymbol{y})}_{\text{mode probability vector}},\ \underbrace{f_{l} \in U_{l},\ \forall l\in \mathcal{L}}_{\text{distribution in each mode}}
\right\}.
\end{align}
Here, the probability vector $\boldsymbol{p} = (p_1,\ldots,p_L)^{\mathsf T}$ corresponds to the proportion of each patient type within the group. We use the proportion of each patient type $l$ in the Mayo Clinic's dataset to construct a nominal probability $\hat{p}_l$ for each $l\in\mathcal{L}$ (i.e., the ``Proportion'' column in Table \ref{tab:duration-stats}). Then, according to the group activation and assignment decisions $(\boldsymbol{x},\boldsymbol{y})$, we adopt a decision-dependent variation distance-based ambiguity set $\Delta_g\bigl(\hat{\boldsymbol{p}}, \boldsymbol{x},\boldsymbol{y})$ as follows:  
\begin{align}\label{eq:variation distance}
\Delta_g\bigl(\hat{\boldsymbol{p}}, \boldsymbol{x},\boldsymbol{y}) = \left\{\, \boldsymbol{p} \in \mathbb{R}^L_+ : p_l \leq x_{lg},\ \sum_{l=1}^L \left|p_l(\sum_{l=1}^L x_{lg}\hat{p}_l)-x_{lg}\hat{p}_l\right| \leq \rho \sum_{l=1}^L x_{lg}\hat{p}_l,\ 
\sum_{l=1}^L p_{l} =y_g\right\}.
\end{align}

For a given group $g\in \mathcal{G}$, this formulation restricts mode probabilities $\boldsymbol{p}$ to deviate from the nominal ones within a radius $\rho$. In particular, if group $g$ is not activated ($y_g =0$), then no patient type is assigned to this group ($x_{lg}=0,\ \forall l \in \mathcal{L}$), which enforces all mode probabilities to be $0$, i.e., $p_l=0,\ \forall l\in\mathcal{L}$. As a result, this group will not contribute any idle time or overtime. If group $g$ is activated ($y_g=1$), the nominal probability for type $l$ within group $g$ is given by the conditional distribution $\frac{x_{lg}\hat{p}_l}{\sum_{l=1}^Lx_{lg}\hat{p}_l}$, which rescales the proportions $\hat{p}_l$ so that they sum to one across all patient 
types assigned to group $g$. Dividing both sides of the second condition in \eqref{eq:variation distance} by $\sum_{l=1}^L x_{lg}\hat{p}_l$, we require the total variation distance between the mode probability $p_l$ and the nominal one $\frac{x_{lg}\hat{p}_l}{\sum_{l=1}^Lx_{lg}\hat{p}_l}$ for all $l=1,\ldots,L$ to be within $\rho$. Moreover, when patient type $l$ is not assigned to the group $g$ ($x_{lg}=0$), the corresponding probability $p_l$ must be zero. 

Given a mode $l$, the second layer characterizes the uncertainty within each mode $l$ by describing the random duration through the unknown density $f_l$ that belongs to a moment-based ambiguity set $U_l$. To capture the asymmetric information of the duration distribution within each mode, we follow \cite{natarajan2018asymmetry} by including a second-order statistic (i.e., semivariance). For any random variable $\tilde{\xi}$, the normalized semivariance is defined as $
s \;=\; 
\frac{
\mathbb{E}\Bigl[ \bigl( (\tilde{\xi} - m)^{+} \bigr)^2 \Bigr]
\;-\;
\mathbb{E}\Bigl[ \bigl( (m - \tilde{\xi})^{+} \bigr)^2 \Bigr]
}{
\sigma^2
}$, where $m$ denotes the mean of $\tilde{\xi}$ and $\sigma$ denotes the standard deviation of $\tilde{\xi}$. The normalized semivariance $s$ captures asymmetry by measuring the imbalance between the variability above and below the mean $m$. A positive $s$ indicates a right-skewed distribution, while a negative $s$ corresponds to a left-skewed distribution.
Consequently, we define the asymmetric moment-based ambiguity set $U_l$ as
\begin{align}\label{eq:ul ambiguity set}
U_l:= \left\{f_l \in \mathcal{P}(\Xi): \mathbb{E}_{f_l}\bigl[\tilde{\xi}\bigr] = m_l,\ \mathbb{E}_{f_l}\bigl[(\tilde{\xi}-m_l)^2\bigr]=\sigma_l^2,\ \mathbb{E}_{f_l}\bigl[(\tilde{\xi}-m_l)^{+2}\bigr]-\mathbb{E}_{f_l}\bigl[(m_l-\tilde{\xi})^{+2}\bigr]=s_l \sigma_l^2\right\},
\end{align}
where the parameters $m_l,\ \sigma_l$ and $s_l$ represent the mean, variance, and normalized semivariance of $\tilde{\xi}$ in mode $l$, respectively, and $\mathcal{P}(\Xi)$ denotes all the probability distributions supported on $\Xi$. This layer of ambiguity captures the essential second-order information and asymmetric tail behavior of the distribution of chemotherapy duration. 
 

Under the ambiguity set \eqref{eq:ambiguity set}, we further reformulate our Problem \eqref{eq:MMA-DRO} as
\begin{align}\label{eq:two stage problem objective}
    \min_{x_{lg},y_g ,t_g} \quad &  \left\{\sum_{g \in \mathcal{G}} c_gy_g + \frac{1}{\sum_{g \in \mathcal{G}} y_g}\sum_{g \in \mathcal{G}} \max_{\boldsymbol{p} \in \Delta_g(\hat{\boldsymbol{p}}, \boldsymbol{x},\boldsymbol{y})} \sum_{l=1}^L p_l \max_{f_l \in U_l} \mathbb{E}_{\tilde{\xi} \sim f_l}\left[q(\tilde{\xi}-t_g)^+ + b(t_g-\tilde{\xi})^+\right]\right\} \\
    \textrm{s.t.} \quad & \textrm{Constraints}\  \eqref{eq:sum to one},\ \eqref{eq:x less than y},\ \eqref{eq:feasibility},\ \eqref{eq:t range}.\nonumber
\end{align}

To better represent the problem formulation, we first denote the innermost maximization problem (i.e., worst-case cost of idle time and overtime) in Problem\ \eqref{eq:two stage problem objective} for a given group $g$ and mode $l$ as 
\begin{align}\label{eq:pl known inner}
    \Pi_l(t_g):=\max_{f_l \in U_l}\mathbb{E}_{\tilde{\xi} \sim f_l}\left[q(\tilde{\xi}-t_g)^+ + b(t_g-\tilde{\xi})^+\right],
\end{align}
and make the following assumption to derive a dual representation of $\Pi_l(t_g)$.

\begin{assumption}\label{ass: feasible assumption}
    For each type $l \in \mathcal{L}$, the treatment duration has mean $m_l >0$, standard deviation $\sigma_l >0$ and normalized semivariance $s_l$ satisfying $s_l \in \Bigl[\frac{\sigma_l^2 - m_l^2}{\sigma_l^2 + m_l^2},1\Bigr)$.
\end{assumption}
According to Proposition 2.1 in \cite{natarajan2018asymmetry}, Assumption \ref{ass: feasible assumption} provides a necessary and sufficient condition on the mean, standard deviation, and normalized semivariance of a nonnegative random variable, which ensures that the ambiguity set $U_l$ is non-empty.

In the following theorem, we present a dual representation of $\Pi_l(t_g)$. Throughout the paper, all the proofs are presented in Appendix\ \ref{sec:omitted proofs}.

\begin{theorem}[Dual representation of $\Pi_l(t_g)$]\label{thm: dual of Pi}
Under Assumption \ref{ass: feasible assumption}, $\Pi_l(t_g)$ admits the following dual formulation
\begin{align}\label{eq:continuous ul ambiguity set}
\Pi_l(t_g)= \min_{\gamma_{lg},\nu_{lg},\kappa_{lg},\zeta_{lg}} \quad &
  \gamma_{lg} + \nu_{lg} m_l + \kappa_{lg} \sigma_l^2 + \zeta_{lg} s_l \sigma_l^2 \nonumber \\
\text{s.t.}\quad &
  \gamma_{lg} + \nu_{lg}\tilde{\xi} + \kappa_{lg}(\tilde{\xi}-m_l)^2 
  + \zeta_{lg}\!\left((\tilde{\xi}-m_l)^{+2} - (m_l-\tilde{\xi})^{+2}\right) \nonumber \\
& \ge q(\tilde{\xi}-t_g)^+ + b(t_g-\tilde{\xi})^+,\quad \forall\, \tilde{\xi}\in\Xi, \nonumber \\
& \gamma_{lg},\nu_{lg},\kappa_{lg},\zeta_{lg}\in\mathbb{R}.
\end{align}
\end{theorem}

Next, we investigate the overall maximization problem by including the mode probability ambiguity set $\Delta_g(\boldsymbol{\hat{p}},\boldsymbol{x}, \boldsymbol{y})$. For each group $g \in \mathcal{G}$, we denote $$\Omega_g(t_g) : = \max_{\boldsymbol{p} \in \Delta_g(\boldsymbol{\hat{p}},\boldsymbol{x}, \boldsymbol{y})} \sum_{l=1}^L p_l\Pi_l(t_g),$$ and derive its dual representation in the next theorem.

\begin{theorem}[Dual representation of $\Omega_g(t_g)$]\label{thm:dual of Omega}
Under Assumption \ref{ass: feasible assumption}, $\Omega_g(t_g)$
admits the following dual formulation
\begin{subequations}\label{eq: dual of overall inner max}
\begin{align}
    \Omega_g(t_g) = \min_{\mu_g,\lambda_{lg},\tau_g, \alpha_l,\beta_l} \quad & \mu_g y_g + \sum_{l=1}^L \lambda_{lg} x_{lg} + \sum_{l=1}^L (\alpha_{lg} - \beta_{lg}) \hat{p}_lx_{lg} +  \rho \tau_g \sum_{l=1}^L x_{lg}\hat{p}_l\\
    \textrm{s.t.} \quad & \mu_g + \lambda_{lg} + (\alpha_{lg}-\beta_{lg})(\sum_{l=1}^L x_{lg}\hat{p}_l) \geq \Pi_l(t_g),\ \forall l \in \mathcal{L}\label{eq:rhs pi}\\
    \quad & -\alpha_{lg} - \beta_{lg} + \tau_g \geq 0,\ \forall l \in \mathcal{L} \\
    \quad & \mu_g \in \mathbb{R},\ \lambda_{lg},\  \alpha_{lg},\ \beta_{lg},\ \tau_g \geq 0.
\end{align}
\end{subequations}
\end{theorem}

Combining Theorems \ref{thm: dual of Pi} and \ref{thm:dual of Omega}, we derive a monolithic minimization problem of our MMA-DRO model\ \eqref{eq:MMA-DRO} as follows
\begin{subequations}\label{eq:continuous overall min}
\begin{align}
\min \quad & \Biggl\{\sum_{g \in \mathcal{G}} c_gy_g + \frac{1}{\sum_{g \in \mathcal{G}} y_g}\sum_{g \in \mathcal{G}} \Bigl(\mu_g y_g + \sum_{l=1}^L \lambda_{lg} x_{lg} + \sum_{l=1}^L (\alpha_{lg} - \beta_{lg}) \hat{p}_l x_{lg}+\tau_g (\rho \sum_{l=1}^L x_{lg}\hat{p}_l) \Bigr) \Biggr\} \nonumber\\
    \textrm{s.t.}  \quad & \mu_g + \lambda_{lg} + (\alpha_{lg}-\beta_{lg})(\sum_{l=1}^L x_{lg}\hat{p}_l) \geq \gamma_{lg} + \nu_{lg} m_l + \kappa_{lg} \sigma_l^2 + \zeta_{lg} s_l \sigma_l^2,\ \forall l \in \mathcal{L},\ g\in \mathcal{G}\label{eq:discretized approximation first constraint}\\
    \quad & \gamma_{lg} + \nu_{lg}\tilde{\xi} + \kappa_{lg}(\tilde{\xi}-m_l)^2 + \zeta_{lg} \left((\tilde{\xi}-m_l)^{+2} - (m_l-\tilde{\xi})^{+2}\right) \nonumber \\
    \quad & \hspace{5.2cm} \geq q(\tilde{\xi}-t_g)^+ + b(t_g-\tilde{\xi})^+,\ \forall l \in \mathcal{L},\ \forall g \in \mathcal{G},\ \forall \tilde{\xi} \in \Xi \label{eq:discrete change}\\
    \quad & -\alpha_{lg} - \beta_{lg} + \tau_g \geq 0,\ \forall l \in \mathcal{L},\ g\in \mathcal{G} \\
    \quad & \text{Constraints}\ \eqref{eq:sum to one},\ \eqref{eq:x less than y},\ \eqref{eq:feasibility},\ \eqref{eq:t range} \nonumber \\
    \quad & \gamma_{lg},\ \nu_{lg},\ \kappa_{lg},\ \zeta_{lg},\ \mu_g \in \mathbb{R},\ \lambda_{lg},\  \alpha_{lg},\ \beta_{lg},\ \tau_g,\ t_g \geq 0 \\
    \quad & x_{lg},\ y_g \in \{0,1\},\ \forall l \in \mathcal{L},\ g\in \mathcal{G}
\end{align}
\end{subequations}

This formulation, however, results in an intractable semi-infinite optimization problem due to the continuous support set $\Xi$ in Constraints\ \eqref{eq:discrete change}, thereby preventing the direct use of standard optimization solvers. To overcome this intractability, we first derive a closed-form expression for $\Pi_l(t_g)$, based on which we develop two exact algorithms for solving this semi-infinite program \eqref{eq:continuous overall min} in Section\ \ref{sec:exact}. Although these exact algorithms guarantee to find optimal solutions, they become computationally expensive when the problem size increases. To further improve the computational time, in Section\ \ref{sec: heuristic}, we develop several heuristics based on clustering algorithms (e.g., K-Means and K-Medoids), which can find near-optimal solutions in a quick fashion. We also consider a generalized ambiguity set $U_l$ that allows the mean, variance, and semivariance to deviate from the empirical ones within a given radius and derive its corresponding dual form under a discrete support set in Appendix\ \ref{app:mmadro_mis}. We will compare it with other benchmarks in Section\ \ref{sec: synthetic section}.


\section{Exact Algorithms for Solving MMA-DRO Model \eqref{eq:continuous overall min}}\label{sec:exact}

In the semi-infinite program \eqref{eq:continuous overall min}, the challenging part lies in Constraints \eqref{eq:discretized approximation first constraint}--\eqref{eq:discrete change} (or, equivalently, Constraints \eqref{eq:rhs pi}), which involve the worst-case expected cost $\Pi_l(t_g)$ for mode $l$ in group $g$.  In this section, we first derive a closed-form expression for $\Pi_l(t_g)$, which takes the form of a piecewise convex function in $t_g$. Then given a specific assignment decision $(\boldsymbol{x},\boldsymbol{y})$, we leverage this closed-form expression to develop an exact algorithm to solve the semi-infinite program\ \eqref{eq:continuous overall min} by exploiting its first‑order conditions when $\rho=0$ (i.e., the variation distance set\ \eqref{eq:variation distance} reduces to a singleton) in Section\ \ref{sec:pl known}. When $\rho\not=0$, we replace the right-hand side of\ \eqref{eq:discretized approximation first constraint} and Constraint\ \eqref{eq:discrete change} using the closed-form expression of $\Pi_l(t_g)$ and develop an exact algorithm to solve\ \eqref{eq:continuous overall min} in Section\ \ref{sec:pl unknown}.

Specifically, given an assignment decision $(\boldsymbol{x}, \boldsymbol{y})$, Problem\ \eqref{eq:MMA-DRO} can be decomposed by group $g$ and each group $g$ needs to solve the following min-max problem
\begin{align}\label{eq:heuristic problem}
    \min_{t_g \in \mathcal{T}}\max_{\boldsymbol{p} \in \Delta_g(\hat{\boldsymbol{p}}, \boldsymbol{x},\boldsymbol{y})} \sum_{l=1}^L p_l \Pi_l(t_g).
\end{align}   

Next, we derive a closed-form expression for $\Pi_l(t_g)$ in Theorem\ \ref{thm:thm1}.


\begin{theorem}\label{thm:thm1}
Under Assumption \ref{ass: feasible assumption}, $\Pi_l(t_g)$ admits the following closed‑form expression for each mode $l$ and group $g$:
\begin{align}
     & \Pi_l(t_g)=f_{l1}(t_g):=\left(\frac{(b+q)(1-s_l)\sigma_l^2}{2m_l^2}-q\right)t_g + qm_l,\ \text{if}\ t_g \in \mathcal{T}_{l1}:=\left[0,\frac{m_l}{2}\right] \nonumber\\
    & \Pi_l(t_g)=f_{l2}(t_g):=\frac{(b+q)(1-s_l)\sigma_l^2}{8(m_l-t_g)}-qt_g+qm_l,\ \text{if}\ t_g \in \mathcal{T}_{l2}:=\left[\frac{m_l}{2},m_l-\frac{\sigma_l}{2}\sqrt{\frac{1-s_l}{1+s_l}}\right] \nonumber \\
     & \Pi_l(t_g)=f_{l3}(t_g):=\frac{(b+q)\sigma_l}{2}\sqrt{1-s_l^2} + (m_l-t_g)\frac{(q-b)-(q+b)s_l}{2},\nonumber\\
     &\hspace{8cm}\text{if}\ t_g \in \mathcal{T}_{l3}:=\left[m_l-\frac{\sigma_l}{2}\sqrt{\frac{1-s_l}{1+s_l}},m_l+\frac{\sigma_l}{2}\sqrt{\frac{1+s_l}{1-s_l}}\right] \nonumber\\
     & \Pi_l(t_g)=f_{l4}(t_g):=\frac{(b+q)(1+s_l)\sigma_l^2}{8(t_g-m_l)}+bt_g-bm_l,\ \text{if}\ t_g \in \mathcal{T}_{l4}:=\left[m_l+\frac{\sigma_l}{2}\sqrt{\frac{1+s_l}{1-s_l}},m_l + \frac{m_l(1+s_l)}{2(1-s_l)}\right] \nonumber \\
    & \Pi_l(t_g)=f_{l5}(t_g):=bt_g+qm_l-\frac{b+q}{2}\left[m_l+\beta_l t_g-(t_g-m_l)\sqrt{\beta_l \left(\beta_l+\frac{m_lw_{1l}-2w_{2l}(t_g-m_l)}{m_l(t_g-m_l)^2}\right)}\right],\nonumber \\ 
     &\hspace{9.5cm} \text{if}\ t_g \in \mathcal{T}_{l5}:=\left[m_l + \frac{m_l(1+s_l)}{2(1-s_l)},T \right].
\end{align}
where $w_{1l} = \frac{(1+s_l)\sigma_l^2}{2},\ w_{2l} = \frac{(1-s_l)\sigma_l^2}{2},\ \beta_l = 1-\frac{(1-s_l)\sigma_l^2}{2m_l^2}$.
\end{theorem}

The proof of Theorem\ \ref{thm:thm1} is shown in Appendix\ \ref{proof: thm1}, which consists of constructing a dual feasible solution and a primal feasible distribution with the same objective values. By strong duality, this establishes optimality. The construction reduces to analyzing the relationship between the piecewise function 
$\Pi_l(t_g)$ and the quadratic function of the dual formulation of the ambiguity set $U_l$.

Note that in Theorem\ \ref{thm:thm1}, for each mode $l$, we partition the feasible set $\mathcal{T}$ into five pieces $\mathcal{T}_{li},\ i=1,\ldots,5$ and derive a closed-form expression for $\Pi_l(t_g)$ on each piece, i.e., $f_{li}(t_g),\ i=1,\ldots,5$. From the primal formulation in \eqref{eq:pl known inner}, we observe that $\Pi_l(t_g)$ is convex in $t_g$ because maximum and expectation operators preserve convexity. Therefore, these five pieces $f_{li}(t_g),\ i=1,\ldots,5$ are linear or convex functions in $t_g$. 

\subsection{$\rho =0$}\label{sec:pl known}
When we are confident in the nominal mode probability estimated from the real-world dataset, we can set $\rho=0$. In this case, the ambiguity set\ \eqref{eq:variation distance} collapses to a singleton, i.e., $p_l = \frac{x_{lg}\hat{p}_l}{\sum_{l=1}^L x_{lg}\hat{p}_l},\ \forall l\in\mathcal{L}$ and Problem\ \eqref{eq:heuristic problem} becomes
$\min_{t_g\in \mathcal{T}} \Pi(t_g):=\sum_{l=1}^L p_l\Pi_{l}(t_g)$.

Since each $\Pi_{l}(t_g)$ is a piece-wise convex function with intervals $\mathcal{T}_{li}$ for $i=1,\ldots,5$, taking a convex combination of $\Pi_{l}(t_g)$ for all $l=1,\ldots,L$ leads to a convex function $\Pi(t_g)$. All such intervals $
\bigcup_{\substack{l:\, x_{lg}=1 \\ i=1,\ldots,5}} \{\mathcal{T}_{li}\}$ will partition the feasible set $\mathcal{T}$ into disjoint intervals $\{\mathcal{T}_j\}_{j=1}^J$, where $\mathcal{J}$ denotes the total number of pieces in $\Pi(t_g)$. For each interval $\mathcal{T}_j$, we construct an explicit representation of $\Pi(t_g)$, using the closed-form expression in $\Pi_l(t_g)$. Consequently, the optimal solution $t_g^*$ can be characterized by the subgradient optimality condition, i.e. $0 \in \partial_{t_g} \Pi(t_g^*)+ \mathcal{N}_{\mathcal{T}}(t_g^*)$, where $\partial_x f(x)$ denotes the subgradient of $f(x)$ at $x$ and $\mathcal{N}_{\mathcal{T}}(t_g^*)$ denotes the normal cone of $\mathcal{T}$.

\begin{theorem}\label{thm:convexity}[Convexity and global optimality condition.]
$\Pi(t_g)$ is convex in $t_g$, and $t_g^*$ is optimal if and only if $0 \in \partial_{t_g} \Pi(t_g^*)+\mathcal{N}_{\mathcal{T}}(t_g^*)$. Specifically, the optimal solution $t_g^*$ is either obtained at a stationary point of one of the convex function pieces or at one of the endpoints of all intervals $\{\mathcal{T}_j\}_{j=1}^J$.
\end{theorem}

According to Theorem \ref{thm:convexity}, we transition from solving the intractable semi-infinite problem~\eqref{eq:continuous overall min} to a simple root‑finding problem, which is more computationally efficient. Rather than solving a generic
root-finding problem on all intervals and checking all breakpoints, we first analyze the functional form on each piece and identify a set of conditions on the parameters that will lead to a closed-form optimal $t_g^*$.

\begin{proposition}\label{prop:parameter analysis pieces}
    When $\sum_{l:x_{lg}=1} p_l (1-s_l)\left(\frac{\sigma_l}{m_l}\right)^2 \;\ge\; \frac{2q}{\,b+q\,}$, the optimal time-slot duration for group $g$ is $t_g^*=0$. When $\frac{2b}{b+q} \le \sum_{l:x_{lg}=1} p_l \left(\beta_l-\sqrt{S_l(T)}
+\frac{\beta_l}{\sqrt{S_l(T)}}\!\left(
\frac{(1+s_l)\sigma_l^2}{2\,\Delta_{l}^{2}}
-\frac{(1-s_l)\sigma_l^2}{2\,m_l\,\Delta_{l}}
\right)\right)$, where $S_l(T)=\beta_l\!\left(\beta_l+\frac{(1+s_l)\sigma_l^2}{2\,\Delta_T^{2}}
-\frac{(1-s_l)\sigma_l^2}{m_l\,\Delta_T}\right)$, the optimal time-slot duration for group $g$ is $t_g^*=T$.
\end{proposition}

\begin{corollary}\label{coro:parameter}
    When the unit penalty of overtime vanishes ($q=0$) or the unit penalty of idle time goes to infinity ($b\to+\infty$), $\sum_{l:x_{lg}=1} p_l (1-s_l)\left(\frac{\sigma_l}{m_l}\right)^2 \;\ge\; \frac{2q}{\,b+q\,}\to0$. As a result, setting $t_g^*=0$ avoids large idle time penalties. On the other hand, when the unit penalty of idle time vanishes ($b=0$) or the unit penalty of overtime goes to infinity ($q\to+\infty$), $\frac{2b}{\,b+q\,}\to0$ and the condition in Proposition \ref{prop:parameter analysis pieces} is satisfied. Therefore, the minimum is obtained at $t_g^* =  T$ to avoid large overtime penalties.
\end{corollary}

\begin{figure}[ht!]
    \centering
\includegraphics[width=0.97\linewidth]{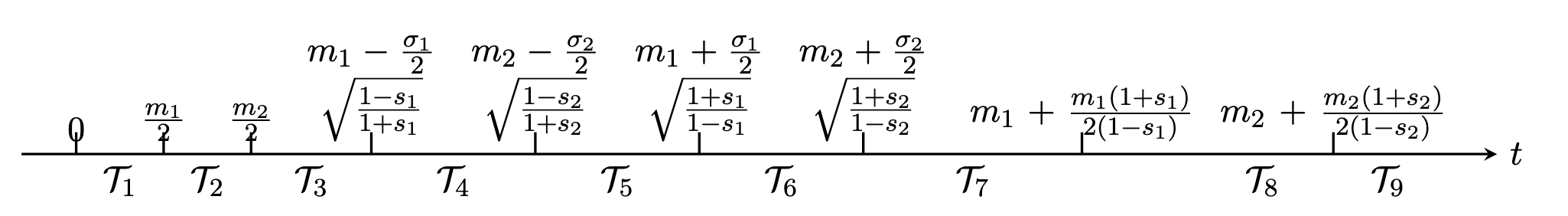}
    \caption{A toy example illustration for $L=2$.}
    \label{fig:schematic}
\end{figure}

For non-trivial parameter setups, we provide a road map for finding the closed-form solution $t_g^*$ for group $g$ given the grouping assignment decision $(\boldsymbol{x},\boldsymbol{y})$ in Algorithm \ref{algo:theorem2}.
\begin{algorithm}[ht]
\SetKwInOut{Input}{Input}
\SetKwInOut{Output}{Output}

\Input{Grouping assignments $(\boldsymbol{x},\boldsymbol{y})$, treatment–duration data for group $g$: $\{(m_l, \sigma_l, s_l),\hat{p}_l\}_{l: x_{lg}=1}$, and cost parameters $b,\ q$.}

\BlankLine
Construct nominal mode probability $p_l=\frac{x_{lg}\hat{p}_l}{\sum_{l=1}^L x_{lg}\hat{p}_l}$ for all $l\in\mathcal{L}$.

Compute intervals $\mathcal{T}_{li}$ for all $l$ such that $x_{lg}=1$ and $i=1,\ldots, 5$ following Theorem\ \ref{thm:thm1}. 

All of these intervals $\bigcup_{\substack{l:\, x_{lg}=1 \\ i=1,\ldots,5}} \{\mathcal{T}_{li}\}$ partition the feasible
set $\mathcal{T}$ into disjoint intervals
$\bigl\{\mathcal{T}_j\bigr\}_{j=1}^{J}$. 

\For {$j=1,\dots,J$}
{
    
      Write out $\Pi(t_g)=\sum_{l: x_{lg}=1} p_l \Pi_l(t_g)$ for $t_g \in \mathcal{T}_j$ explicitly based on the closed-form expression of $\Pi_l(t_g)$ in Theorem\ \ref{thm:thm1}.
      
      Take the derivative of $\Pi(t_g)$ with respect to $t_g$ on interval $\mathcal{T}_j$, and solve $\odv{\Pi(t_g)}{t_g}=0$ to get the root $r_{jg}$.
    

\IfNoEnd{$r_{jg} \in \mathcal{T}_j$}{
  Set $t_g^* = r_{jg}$\;
  \textbf{Break}
}
\Else{
  Denote the endpoints of interval $\mathcal{T}_j$ as $\{\underline t_j,\overline t_j\}$\;
  Compute $J(\underline t_j)=\sum_{l:x_{lg}=1} p_l \Pi_l(\underline t_j)$ and 
  $J(\overline t_j)=\sum_{l:x_{lg}=1} p_l \Pi_l(\overline t_j)$\;
  Set $r_{jg} = \arg\min_{t \in \{\underline t_j,\overline t_j\}} J(t)$\;
}
}

If $t_g^*$ has not been found,
then calculate the optimal time-slot duration for group $g$ as $t_g^*=\arg\min_{t_g\in\{r_{jg}\}}\sum_{l:x_{lg}=1} p_l \Pi_l(t_{g})$.

\Output{Optimal duration $t^*_g$ and optimal cost $\sum_{l:x_{lg}=1} p_l \Pi_l(t_{g}^*)$ for group $g$.}
\caption{Finding the optimal duration and cost for group $g$ given $(\boldsymbol{x},\boldsymbol{y})$ and $\rho=0$.}
\label{algo:theorem2}
\end{algorithm}
For illustration, consider a toy example with $l=1,2$ are assigned to group $g$ and the corresponding parameters satisfying the following relationship
\[
\begin{aligned}
\tau_{1}&:= 0 < \tau_{2} := \frac{m_{1}}{2} < \tau_{3} := \frac{m_{2}}{2} < \tau_{4} := m_{1} - \frac{\sigma_{1}}{2} \sqrt{\frac{1-s_1}{1+s_1}}< \tau_{5} := m_{2} - \frac{\sigma_{2}}{2}\sqrt{\frac{1-s_2}{1+s_2}} \\
&< \tau_{6} := m_{1} + \frac{\sigma_{1}}{2}\sqrt{\frac{1+s_1}{1-s_1}} < \tau_{7} := m_{2} + \frac{\sigma_{2}}{2} \sqrt{\frac{1+s_2}{1-s_2}}< \tau_{8} := m_{1} + \frac{m_{1}(1+s_{1})}{2(1-s_{1})} < \tau_{9} := m_{2} + \frac{m_{2}(1+s_{2})}{2(1-s_{2})}
\end{aligned}
\]
where all the breakpoints partition the feasible region $\mathcal{T}$ into nine disjoint intervals $\mathcal{T}_j=[\tau_j,\tau_{j+1})$, $j=1,\ldots,9$ with 
$\tau_{10}:=T$ (see Figure\ \ref{fig:schematic}).
For group $g$, we write $\Pi(t_g) = p_1 \Pi_1(t_g) + p_2 \Pi_2(t_g)$ and substitute the corresponding explicit formulas for $\Pi_1(t_g)$ and $\Pi_2(t_g)$ on each interval $\mathcal{T}_j$ based on Theorem \ref{thm:thm1}. For example, on $\mathcal{T}_1$, both modes achieve the first linear piece $f_{11}$ and $f_{21}$ and we have $\Pi(t_g) =p_1 f_{11}(t_g) + p_2 f_{21}(t_g)$; on $\mathcal{T}_2$, mode $l=1$ achieves the second piece $f_{12}$ while mode $l=2$ still achieves the first linear piece $f_{21}$ and we have $\Pi(t_g) =p_1 f_{12}(t_g) + p_2 f_{21}(t_g)$. We then solve $\frac{d}{dt_g}\Pi(t_g)=0$ within each $\mathcal{T}_j$ and check if the obtained stationary point is within $\mathcal{T}_j$. If it is true, then we return the stationary point as $t_g^*$ and evaluate $\Pi(t_g^*)$ to obtain the optimal cost for group $g$. Otherwise, the minimizer is attained at one of the
breakpoints. We then evaluate $\Pi_l(t_g)$ at each breakpoint $\tau_i$ for $i=1,\ldots,10$ and set the one with the smallest cost as $t_g^*$.

Algorithm\ \ref{algo:theorem2} returns the optimal time-slot duration $t_g^*$ for each group $g$ given a grouping assignment $(\boldsymbol{x},\boldsymbol{y})$. Next, we consider the overall problem\ \eqref{eq:MMA-DRO}, where the type-to-group assignment is also part of the decisions. We propose Algorithm\ \ref{algo:enumerate all groupings} to find the optimal overall cost by enumerating all possible grouping assignments (i.e., set partitions of the $L$ patient types) and applying Algorithm\ \ref{algo:theorem2} for each specific grouping assignment. 
\begin{algorithm}[ht]
\SetKwInOut{Input}{Input}
\SetKwInOut{Output}{Output}

\Input{The moment $m_l,\ \sigma_l,\ s_l$ and $\hat{p}_l$ for each patient type $l$; cost parameters $b,\ q$.}

\BlankLine

  Put all possible grouping assignments of the $L$ patient types into $\mathcal{P}$.
  
  \For{each grouping assignment $P\in\mathcal{P}$ with exactly $|P|$ groups}{
  \For{$g = 1,\ldots,|P|$}{
    Run Algorithm\ \ref{algo:theorem2} for group $g$ to obtain the optimal duration $t_g^*$ and the optimal cost $\sum_l p_l\Pi_l(t_g^*)$.}

    Calculate the overall cost with respect to assignment $P$: $Cost(P) = c_g |P| + \frac{1}{|P|} \sum_{g=1}^{|P|} \sum_{l=1}^L p_l \Pi_l(t_g^*)$
  }

Select the best grouping assignment $P^* = \argmin_{P\in\mathcal{P}} Cost(P)$.

\Output{optimal overall objective function value $\min_{P\in\mathcal{P}} Cost(P)$, the corresponding assignment decision in $P^*$, and the optimal time-slot duration for each group.}
\caption{Solving MMA-DRO\ \eqref{eq:continuous overall min} by enumerating all grouping assignments}
\label{algo:enumerate all groupings}
\end{algorithm}
Specifically, we enumerate all set partitions of the $L$ patient types. For instance, when $L=3$, this will generate $5$ possible grouping assignments, $\{\{1\}, \{2\}, \{3\}\}, \{\{1,2\},\{3\}\},\{\{1,3\},\{2\}\},\{\{2,3\},\{1\}\}$, and $\{1,2,3\}$. Next, for every possible partition, we apply Algorithm\ \ref{algo:theorem2} on each group $g$ to compute the optimal time-slot duration $t_g^*$ and the optimal cost for group $g$. We then compute the overall cost for this partition P as $Cost(P)$ by including the assignment cost. Finally, we select partition $P$ that achieves the minimum overall cost. This procedure is exact since it searches all possible partitions of $L$ patient types and is tractable for small $L$, where the total number of partitions is the Bell number $B_L$ (\cite{mezo2011r}). The details are outlined in Algorithm\ \ref{algo:enumerate all groupings}.

\subsection{$\rho>0$}\label{sec:pl unknown}
When $\rho=0$, the ambiguity set\ \eqref{eq:variation distance} reduces to a singleton, and we leverage the closed-form expression of $\Pi_l(t_g)$ to find the optimal solutions to MMA-DRO model\ \eqref{eq:MMA-DRO} in Section\ \ref{sec:pl known}. However, when we are not confident in the estimated mode probability, we need to set $\rho>0$. Specifically, Eq.\ \eqref{eq:variation distance} specifies the ambiguity set for the mode probabilities $p_l$, where 
the radius $\rho$ controls the allowable deviation between the true mode probabilities 
$p_l$ and the nominal ones $\hat{p}_l$. In this case, the inner problem also involves the maximization over the ambiguity set $\Delta_g(\boldsymbol{\hat{p}},\boldsymbol{x},\boldsymbol{y})$, and we cannot apply Algorithm\ \ref{algo:theorem2} or\ \ref{algo:enumerate all groupings} to solve it. 

Instead, Theorem\ \ref{thm:dual of Omega} provides a dual representation of this maximization problem over $\Delta_g(\boldsymbol{\hat{p}},\boldsymbol{x},\boldsymbol{y})$, which involves $\Pi_l(t_g)$ on the right-hand side of Eq.\ \eqref{eq:rhs pi}. As a result, for a given grouping assignment 
$(\boldsymbol{x},\boldsymbol{y})$, we can substitute the closed-form expression of $\Pi_l(t_g)$ into the right-hand side of Eq.\ \eqref{eq:rhs pi} and solve a finite-dimensional convex optimization problem using off-the-shelf solvers, instead of solving the semi-infinite program\ \eqref{eq:continuous overall min}. 

Algorithm\ \ref{algo:eq14 algo1} provides a road map for finding the optimal solutions and the optimal objective value for each group $g$ when $\rho>0$. 
\begin{algorithm}[ht]
\SetKwInOut{Input}{Input}
\SetKwInOut{Output}{Output}

\Input{Grouping assignment $(\boldsymbol{x},\boldsymbol{y})$, treatment–duration data for group $g$: $\{(m_l, \sigma_l, s_l),p_l\}_{l: x_{lg}=1}$, and cost parameters $b,\ q$.}

\BlankLine
Construct nominal mode probability $p_l=\frac{x_{lg}\hat{p}_l}{\sum_{l=1}^L x_{lg}\hat{p}_l}$ for all $l\in\mathcal{L}$.

Compute intervals $\mathcal{T}_{li}$ for all $l$ such that $x_{lg}=1$ and $i=1,\ldots, 5$ following Theorem\ \ref{thm:thm1}. 

All of these intervals $\bigcup_{\substack{l:\, x_{lg}=1 \\ i=1,\ldots,5}} \{\mathcal{T}_{li}\}$ partition the feasible
set $\mathcal{T}$ into disjoint intervals
$\bigl\{\mathcal{T}_j\bigr\}_{j=1}^{J}$. 

\For {interval
$\mathcal{T}_j\;(j=1,\dots,J)$,}
{
    
    Solve $\min_{t_g\in\mathcal{T}_j} \Omega_g(t_g)$, where we replace the right-hand side of Eq.\ \eqref{eq:rhs pi} with the corresponding piece of the function $\Pi_l(t_g)$ with $t_g \in \mathcal{T}_j$ for each $l=1,\ldots, L$ according to Theorem\ \ref{thm:dual of Omega}.

    Obtain the optimal time-slot duration $t_g^*$ and optimal cost for group $g$,\ $Cost_g(t_g^*)$.
    }

Calculate optimal cost overall cost as $\sum_{g} c_g y_g + \frac{1}{\sum_g y_g} \sum_g Cost_g(t_g^*)$.

\Output{Optimal treatment time $t^*_g$ for each group $g$ and optimal overall cost}
\caption{Finding the optimal duration and cost for group $g$ given $(\boldsymbol{x},\boldsymbol{y})$ and $\rho>0$.}
\label{algo:eq14 algo1}
\end{algorithm}
Given a grouping assignment ($\boldsymbol{x},\boldsymbol{y}$), all the intervals $\{T_{lj}\}_{l: x_{lg}=1, j=1,...,5}$, partition the feasible set $\mathcal{T}$ into $J$ intervals $\{\mathcal{T}_j\}_{j=1}^J$ and we create $J$ subproblems of\ \eqref{eq: dual of overall inner max} by restricting the feasible set for $t_g$ to each $T_j$. Then in each subproblem, we solve $\min_{t_g\in\mathcal{T}_j} \Omega_g(t_g)$, where we replace the right-hand side of\ \eqref{eq:rhs pi} using the piecewise-convex form from Theorem\ \ref{thm:thm1}. We then select the overall optimal solution by taking the minimum cost across all $J$ subproblems. This 
allows us to solve a convex optimization problem each time with finitely many constraints, which is much easier compared to the semi-infinite program\ \eqref{eq:continuous overall min}. 

Similar to Algorithm\ \ref{algo:enumerate all groupings}, the following algorithm provides an enumeration approach to obtain the grouping assignment and uses Algorithm\ \ref{algo:eq14 algo1} to obtain the optimal treatment duration $t_g^*$ and the corresponding optimal cost. 

\begin{algorithm}[ht!]
\SetKwInOut{Input}{Input}
\SetKwInOut{Output}{Output}

\Input{The moment $m_l,\ \sigma_l,\ s_l$ for each patient type $l$; cost parameters $b,\ q$.}

\BlankLine
  Put all possible grouping assignments of the $L$ patient types into $\mathcal{P}$.\;
  \For{each grouping assignment $P\in\mathcal{P}$ with exactly $|P|$ groups}{\For{$g = 1,\ldots,|P|$}{
    Run Algorithm\ \ref{algo:eq14 algo1} for group $g$ to obtain the optimal duration $t_g^*$ and the optimal cost $\Omega_g(t_g^*)$.}

        Calculate the overall cost with respect to assignment P: $Cost(P) = c_g |P| + \frac{1}{|P|} \sum_{g =1}^{|P|} \Omega_g(t_g^*)$
  }

Select the best grouping assignment $P^* = \argmin_{P\in\mathcal{P}} Cost(P)$.

\Output{The optimal overall objective function value $\min_{P\in\mathcal{P}} Cost(P)$, the corresponding assignment decision in $P^*$, and the optimal time-slot duration for each group.}
\caption{Solving MMA-DRO\ \eqref{eq:MMA-DRO} by enumerating all grouping assignments}
\label{algo:eq14 algo}
\end{algorithm}

\subsection{Lower and Upper Bounds on Each Group's Optimal Cost}\label{sec:insight}
To obtain insights into the optimal grouping assignment decisions, we provide lower and upper bounds on each group's optimal cost, which depend on the difference between the maximum and minimum average treatment time among all the patient types assigned to this group. This further sheds light on how to find the optimal grouping assignments: the optimal grouping should group the adjacent types with similar treatment durations (with closer values of $m_l$), and avoid the clustering of patient types that are distant (with widely different $m_l$).  The proof of the following Propositions\ \ref{prop:lb} and\ \ref{prop:upper bound} are provided in the Appendices\ \ref{app: lower bound pf} and\ \ref{app:ub proof}, respectively.
We first derive a lower bound of the optimal cost of Problem\ \eqref{eq:heuristic problem} in the next proposition.

\begin{proposition}[Lower bound of Problem\ \eqref{eq:heuristic problem}]\label{prop:lb}
For any given group $g\in\mathcal G$, define
\[
m_{\min}(g):=\min_{l:\,x_{lg}=1}m_l,\qquad
m_{\max}(g):=\max_{l:\,x_{lg}=1}m_l,
\]
and
\[
\overline p_{\min}(g):=\max_{\boldsymbol{p}\in\Delta_g(\hat{\boldsymbol p},\boldsymbol x,\boldsymbol y)}
\sum_{\substack{l:\,x_{lg}=1\\ m_l=m_{\min}(g)}} p_l,\qquad
\overline p_{\max}(g):=\max_{\boldsymbol{p}\in\Delta_g(\hat{\boldsymbol p},\boldsymbol x,\boldsymbol y)}
\sum_{\substack{l:\,x_{lg}=1\\ m_l=m_{\max}(g)}} p_l .
\]
Then
\[
\min_{t_g \in \mathcal{T}}\ \max_{\boldsymbol{p} \in \Delta_g(\hat{\boldsymbol{p}}, \boldsymbol{x},\boldsymbol{y})}\ \sum_{l:\,x_{lg}=1} p_l\,\Pi_l(t_g)
\;\ge\;
\frac{b\,\overline p_{\min}(g) q\,\overline p_{\max}(g)}{b\,\overline p_{\min}(g)+q\,\overline p_{\max}(g)} (m_{\max}(g)- m_{\min}(g)).
\]
\end{proposition}

In the next proposition, we provide an upper bound of\ \eqref{eq:heuristic problem}, which is also related to the span $m_{\max}(g) - m_{\min}(g)$ for group $g$.
\begin{proposition}[Upper bound of Problem\ \eqref{eq:heuristic problem}]\label{prop:upper bound}
For any group $g \in \mathcal{G}$, define $\sigma_{\max}(g) = \max_{l:x_{lg}=1}\sigma_l$.
Then the group-$g$ objective in\ \eqref{eq:heuristic problem} satisfies
\begin{align}
    \min_{t_g \in \mathcal{T}} \max_{\boldsymbol{p} \in \Delta_g(\boldsymbol{\hat{p}},\boldsymbol{x},\boldsymbol{y})} \sum_{l:x_{lg}=1} p_l \Pi_l(t_g) \le y_g\max\{b,q\} \Bigl(\sigma_{\max}(g) + \frac{1}{2} \bigl(m_{\max}(g) - m_{\min}(g)\bigr)\Bigr).
\end{align}
\end{proposition}

From Propositions \ref{prop:lb} and \ref{prop:upper bound}, both lower and upper bounds of the optimal objective function value of Problem \eqref{eq:heuristic problem} are positively related to the span $m_{\max}(g) - m_{\min}(g)$. This implies that the optimal cost tends to decrease when patient types with similar average treatment times are assigned together. Leveraging this structural insight, we will propose heuristic approaches that cluster patient types according to the similarity of their treatment time statistics in the next section.

\section{Heuristic Approaches}\label{sec: heuristic}
In Section \ref{sec:exact}, we introduced Algorithms\ \ref{algo:enumerate all groupings} and \ref{algo:eq14 algo}, which enumerate all grouping assignments (set partitions) and achieve the global optimum for $\rho=0$ and $\rho>0$, respectively. However, such enumeration-based approaches can quickly become computationally inefficient as the number of patient types $L$ grows (the number of possible partitions scales with the Bell number $B_L$). Motivated by Propositions \ref{prop:lb} and \ref{prop:upper bound}, in this section, we propose heuristic approaches that leverage clustering algorithms (i.e., K-Means and K-Medoids) to group patient types with similar treatment time statistics, which can produce near-optimal solutions in a quick fashion. 

\paragraph{K-Means-based Heuristic. } We first use the K-Means method to group patient types based on the similarity of their input feature. Let $P_l$ represent input feature of patient $l\in \mathcal{L}$ (e.g., $m_l, \sigma_l, s_l$, or any combination of these moment information). Given the number of activated groups $K$, our goal is to cluster the $L$ patient types to $K$ groups $(K \le L)$ such that each within-cluster sum of squares objective is minimized. The following problem can be seen as a simplification of Model\ \eqref{eq:MMA-DRO}, where we approximate the downstream cost of an assignment by the within-cluster sum of squares. According to Propositions\ \ref{prop:lb} and\ \ref{prop:upper bound}, if the input feature among all the patient types assigned to a group is similar to each other, then it will potentially lead to a lower overall cost.
\begin{align}
\min_{\boldsymbol{x},\boldsymbol{\mu}} \quad & \sum_{l \in \mathcal{L}} \sum_{g=1}^K x_{lg} \|P_l - \mu_k\|^2 \label{eq:wcss}\\
    \textrm{s.t.} \quad & \sum_{g=1}^K x_{lg} = 1,\ \forall l \in \mathcal{L} \nonumber \\
    \quad & x_{lg} \in \{0,1\},\ \forall l \in \mathcal{L},\ \forall g=1,\ldots,K \nonumber 
\end{align}
In Model\ \eqref{eq:wcss}, $x_{lg}$ is the assignment decision as we defined before, and $\mu_g$ is the centroid of cluster $g$. It is non-convex in general because the objective function involves bilinear terms. K-Means algorithm solves\ \eqref{eq:wcss} by iterating between the following two steps: \textbf{(i)} assignment step: given $\boldsymbol{\mu}$, assign each patient type $l$ to its closest cluster $g=\argmin_g ||P_l - \mu_g||^2$, and
\textbf{(ii)} update step: given the assignment $x_{lg}^*$, update the centroids by minimizing the objective in\ \eqref{eq:wcss} given by $\mu_g = \frac{\sum_{l=1}^L x_{lg}P_l}{\sum_{l=1}^L x_{lg}},\ \forall g=1,\ldots. K$.
We iterate between these two steps until assignments stabilize, which yields a local optimum of Problem \eqref{eq:wcss}\ \cite{bottou1994convergence}. 
The final clusters yield a feasible assignment decision $(\boldsymbol{x},\boldsymbol{y})$ for our model\ \eqref{eq:heuristic problem}.

Once we obtain a grouping assignment decision $(\boldsymbol{x},\boldsymbol{y})$, we apply Algorithm\ \ref{algo:theorem2} and Algorithm\ \ref{algo:eq14 algo1} to find the optimal duration $t_g^*$ for each group $g$ when $\rho=0$ and when $\rho>0$, respectively. To tune the number of activated groups $K$, we apply cross-validation on the range $K=1,\ldots, L$ to find the best $K$ that minimizes the average cost on the validation dataset.

\paragraph{K-Medoids-based Heuristic.}Different than the K-Means algorithm, where the centroid can fall outside of the input points $\{P_l\}_{l\in\mathcal{L}}$, the K-Medoids algorithm picks a point from the input dataset as a medoid. Specifically, we cluster patient types using K-Medoids on the same feature vectors $P_l$ for $l \in \mathcal{L}$ with a chosen dissimilarity denoted as $\mathrm{dist}(\cdot,\cdot)$. Given the number of activated groups $K$, we partition $L$ patient types into $K$ clusters by selecting $K$ medoids from the actual data points $\{P_l\}$. Using binary variables $z_i\in\{0,1\}$ to indicate whether type $i$ is selected as a medoid and $r_{li}\in\{0,1\}$ to indicate whether we assign type $l$ to medoid $i$, the K-Medoids problem is formulated as
\begin{align}\label{eq:kmedoids obj}
\min_{\boldsymbol{r},\boldsymbol{z}}\quad 
 & \sum_{l\in\mathcal L}\sum_{i\in\mathcal L} r_{li}\,\mathrm{dist}\!\big(P_l,P_i\big) \nonumber \\
\text{s.t.}\quad
 & \sum_{i\in\mathcal L} r_{li} = 1,\ \forall\, l\in\mathcal L,\nonumber\\
 & r_{li} \le z_i,\ \forall\, l,\ i\in\mathcal L,\nonumber\\
 & \sum_{i\in\mathcal L} z_i = K,\nonumber\\
 & r_{li},\ z_i\in\{0,1\},\ \forall\, l,\ i\in\mathcal L.
\end{align}

In the assignment step, given the current medoids $M=\{\,P_i : z_i=1\,\}$ (i.e., $z_i$ fixed),
for each type $l$ solve\ \eqref{eq:kmedoids obj} to obtain the optimal assignment $r_{li}^*$, which assigns each $l$ to its nearest medoid.
In the update step, given the assignments $r_{li}^*$,
for each cluster $g$ choose the new medoid index and solve for $\boldsymbol{z}$ medoids.
Repeat the assignment and update until the objective no longer changes, which yields a local optimum of Problem\ \eqref{eq:kmedoids obj}\ \cite{schubert2019faster}.

K-Means clusters the observations into $K$ groups by iteratively assigning points to the nearest centroid and iteratively updating centroids to minimize the within-cluster sum of squares. By contrast, K-Medoids partitions the observations into $K$ clusters by selecting the representative medoids from the input data points and assigning each observation to its nearest medoid to minimize total within-cluster dissimilarity. Compared with K-Means, K-Medoids is less sensitive to outliers and to skewed or heavy-tailed 
treatment durations. We will test both methods and compare them with the exact algorithms in Sections\ \ref{sec: synthetic section} and\ \ref{sec: mayo section}.

\section{Computational Results based on Synthetic Datasets}\label{sec: synthetic section}
In this section, we conduct sensitivity analyses on synthetic datasets and evaluate the in-sample and out-of-sample performance of the exact algorithms (Algorithms\ \ref{algo:enumerate all groupings} and\ \ref{algo:eq14 algo}) and the heuristics (K-Means and K-Medoids) on these instances.

\subsection{Data Generation}\label{sec:data generation}
We consider the number of patient types $L$ ranging from $5$ to $8$, and set $G=L$. The group activation cost is set the same across all the groups, i.e., $c_g=80$. We also set the unit penalty of idle time and overtime as $q=30,\ b=20$ at default. The random treatment duration $\tilde{\xi}$ is assumed to have a continuous support $\Xi=[0,720]$, which covers the $12$-hour operational window.  

For each mode $l$, we assume that the treatment time follows a log-normal distribution whose log-mean and log-standard deviation are sampled according to
$\mu_l \sim \mathcal U(ln(100), ln(600))$ and
$\sigma_l \sim \mathcal U(0.5, 1.5)$. We draw a nominal mode probability vector $\boldsymbol{\hat{p}}=(\hat{p}_1,\ldots, \hat{p}_L)$ and consider a training dataset of size $100$. We then allocate these $100$ training data points according to the mode probability and draw $K_l=100\hat{p}_l$ training samples $\{\xi_{lk}\}_{k=1}^{K_l}$ for each mode $l$. These samples are then clipped into the $[0,720]$ range.
From the resulting samples, we compute the empirical mean $m_l$, standard deviation $\sigma_l$, and semi-variance $s_l$ for each mode $l\in\mathcal{L}$; these three moments form the input parameters of the MMA-DRO ambiguity set $U_l$\ \eqref{eq:ul ambiguity set}. We apply the exact Algorithms\ \ref{algo:enumerate all groupings} and\ \ref{algo:eq14 algo} and heuristics (K-Means and K-Medoids) to solve for the grouping assignments $\boldsymbol{x}^*,\boldsymbol{y}^*$, the treatment duration $\boldsymbol{t}^*$, and their corresponding in-sample (IS) costs.

Next, we generate a testing dataset of $1,000$ scenarios by sampling from the same ground-truth log-normal distribution of $\tilde{\xi}$ according to the same nominal mode probability vector $\hat{\boldsymbol{p}}$. We will test the algorithm performance when we perturb the out-of-sample (OOS) distribution in Section \ref{sec:performance-comparison}. We evaluate the optimal solutions from solving the MMA-DRO model\ \eqref{eq:MMA-DRO} on the testing dataset to obtain the OOS costs, which allow us to assess the robustness of our algorithms. We report results by averaging over $10$ independent runs of this procedure. We also consider the following two benchmarks:

\paragraph{Benchmark I: Sample Average Approximation (SAA).} As a benchmark to the MMA-DRO model, we consider an \textit{ambiguity-free} setting where we set $\rho=0$ to reduce the variation distance set $\Delta_g(\boldsymbol{\hat{p}},\boldsymbol{x},\boldsymbol{y})$ to a singleton $\left(\frac{x_{1g}\hat{p}_1}{\sum_{l=1}^L x_{lg}\hat{p}_l},\ldots,\frac{x_{Lg}\hat{p}_L}{\sum_{l=1}^L x_{lg}\hat{p}_l}\right)^{\mathsf T}$ and set the ambiguity set $U_l$ to be the empirical distribution $\frac{1}{K_l} \sum_{k=1}^{K_l} \delta_{\xi_{lk}}$, which is supported on the $K_l$ training data points. Under these settings, Model\ \eqref{eq:MMA-DRO} becomes the SAA counterpart.

\paragraph{Benchmark II: Discretized Approximation.} As another benchmark, we replace the continuous support set $\xi\in\Xi$ in \eqref{eq:discrete change} with a discrete support $\{0,1,\ldots, 720\}$, which can be solved by the off-the-shelf solvers (e.g., Gurobi) directly. Because this discretization approximation enforces only a subset of the constraints in Model\ \eqref{eq:continuous overall min}, it provides a valid lower bound to Model \eqref{eq:continuous overall min}.

Numerical tests are conducted on a MacBook Pro with 8 GB RAM and an Apple M1 Pro chip. Since the discretized approximation of Model \eqref{eq:continuous overall min} involves non-convex terms, we use Gurobi 10.0.0 coded in Python 3.11.0 for solving all non-convex programming models (with Non-Convex parameter set to $2$). The time limit for solving this non-convex program is set to 10 hours.

\subsection{Sensitivity Analysis Results}
We perform a sensitivity analysis for $L=5$ and vary the cost parameters $b$ and $q$, the input feature for the heuristic algorithms, the ambiguity set radius parameter $\rho$, and the skewness of the log-normal distribution to examine how these parameters affect the overall performance. 
Summary statistics of both the training and testing datasets are presented in Table \ref{tab:stats} below.

\begin{table}[ht!]
\centering
\caption{Summary statistics for $L=5$}
\resizebox{\textwidth}{!}{
\begin{tabular}{l|rrrrr|rrrrr}
\toprule
& \multicolumn{5}{c}{Training Data} & \multicolumn{5}{c}{Testing Data}\\
\cmidrule(lr){2-6}\cmidrule(lr){7-11}
Statistic & mode $1$ & mode $2$ & mode $3$ & mode $4$ & mode $5$
         & mode $1$ & mode $2$ & mode $3$ & mode $4$ & mode $5$ \\
\midrule
Min& $26.64$  & $43.55$  & $32.88$  & $54.98$  & $100.46$ & $21.81$ & $9.92$ & $4.24$ & $20.92$ & $26.76$ \\
Max  & $560.74$ & $720.00$ & $720.00$ & $720.00$ & $720.00$ & $720.00$ & $720.00$ & $720.00$ & $720.00$ &  $720.00$\\
Mean & $231.58$ & $302.44$ & $366.66$ & $378.14$ & $423.30$ &  $276.97$ & $308.72$ & $333.83$ & $353.10$ & $380.01$  \\
Std.& $170.16$ & $243.91$ & $252.87$ & $248.99$ & $200.03$ &  $200.59$& $257.10$ & $257.25$ & $253.57$ & $229.99$ \\
Semi-variance  & $0.25$   & $0.28$   & $0.16$   & $0.08$   & $0.11$   & $0.36$ & $0.25$ &$0.15$  &$0.14$  & $0.13$ \\
\bottomrule
\end{tabular}}
\label{tab:stats}
\end{table}


\paragraph{Comparison of Different Cost Parameters.}
We first vary the unit penalty of idle time $b$ from $0$ to $30$ and vary the unit penalty of overtime $q$ from $20$ to $50$. The corresponding optimal time-slot durations, IS and OOS costs, and the runtime are reported in Table\ \ref{tab:cost parameter comparison}, where we put the patient types assigned to each group in parentheses.



\begin{table}[ht]
\centering
\caption{Comparison of different cost parameters $b,\ q$.}
\resizebox{\textwidth}{!}{ 
\begin{tabular}{@{} l c c c c c c @{}}
\toprule
\multicolumn{1}{c}{$b,\ q$} &
\multicolumn{1}{c}{Optimal Duration (min.)} &
\multicolumn{1}{c}{IS} &
\multicolumn{3}{c}{OOS} &
\multicolumn{1}{c}{Runtime (sec.)} \\
\cmidrule(lr){4-6}
 &  &   & Idle time & Overtime& Total cost & \\
\midrule
$b=20,\ q=30$ & $164\ (1),\ 516\ (2,3,4,5)$  & $5112.99$ & $189.04$ & $72.82$ & $6125.17$ & $0.99$ \\
$b=30,\ q=20$ & $149\ (1),\ 241\ (2,3,4,5)$ & $4404.14$ & $46.25$ & $154.02$  & $4627.90$ & $1.34$ \\
$b=25,\ q=25$ & $158\ (1),\ 258\ (2,3,4,5)$&  $4899.70$ &$53.37$ & $145.84$  & $5140.26$ &  $0.97$ \\
$b=0,\  q=50$ & $720\ (1,2,3,4,5)$&$80$ & $389.95$ & $0$  & $80$ &  $0.37$ \\
$b=50,\  q=0$ & $0\ (1,2,3,4,5)$&$80$ & $0$ & $329.05$  & $80$ &  $0.49$ \\
\bottomrule
\end{tabular}}
\label{tab:cost parameter comparison}
\end{table}

We begin with the default setting, where the overtime cost parameter $q$ is higher than the idle-time cost $b$ to reflect clinical reality in practice. In this case, the optimal assignment is to group the last four modes, while keeping the first mode a separate group. The second row shows that when the penalty is reversed so that $\frac{b}{q}>1$, the optimal grouping is unchanged, but the optimal durations decrease, leading to smaller IS and OOS costs. With $\frac{b}{q}<1$, the model favors longer slots to avoid costly overtime, whereas with $\frac{b}{q}>1$, idle time is more costly and yields shorter optimal durations. When $\frac{b}{q}=1$, the optimal durations fall between these two cases. 

An extreme setting ($b=0,\ q=50$) illustrates this effect better: the model stretches the duration to the maximum value ($720$ minutes), which yields no idle cost because $b=0$ nor overtime cost because $(\tilde{\xi}-t_g)^+ = 0$, and the lowest IS cost is obtained when only one group is open. Another extreme setting ($b=50,\ q=0$) drives the optimal duration to the minimum value ($0$ minutes), which yields no overtime cost because $q=0$ nor idle cost because $(t_g-\tilde{\xi})^+ = 0$. Both agree with our theoretical results in Corollary \ref{coro:parameter}.
These observations confirm that the relative magnitudes of $b$ and $q$ govern the trade-off between idle and overtime cost, and therefore shifting the optimal treatment duration $t_g^*$ accordingly. 

\paragraph{Comparison of Different Input Features in the Heuristic Algorithms.}
We compare different input features for the heuristic algorithms (K-Means and K-Medoids) when $\rho = 0$ in this section and present the corresponding results when $\rho =0.1$ in Appendix\ \ref{app:rho =0.1 feature selection}. For each input feature combination, we perform cross-validation (CV) over $k = 1,\ldots,L$ on the training data with $L=5$ to select the best $k$. With this $k$ fixed, we run the algorithms to obtain the IS solution $\boldsymbol{t}^*$ and its IS cost. Then we evaluate $\boldsymbol{t}^*$ on the testing dataset to compute the OOS cost. We repeat the entire procedure for $10$ independent runs and report the averages in the Tables\ \ref{tab: compare input feature kmeans fsolve} and\ \ref{tab: compare input feature kmedoids fsolve}, for K-Means and K-Medoids algorithms, respectively, where we mark the best IS and OOS costs in bold below.

\begin{table}[ht!]
\centering
\caption{Comparison of different input features for K-Means}

\begin{tabular}{lccccc}
\toprule
{Feature} & {IS Cost} & {OOS Cost} & {CV Time (sec.)} & {Solve Time (sec.)} & {Total Time (sec.)} \\
\midrule
$(m_l,\sigma_l,s_l)$ & $\bf{4439.88} $ & $4444.85$ & $0.90$ & $0.03$ & $0.93$ \\
$(m_l)$  & $  4448.02$ & $ 4450.49$ & $0.96$ & $0.04$ & $1.00$ \\
$(\sigma_l)$  & $4579.20$& $4517.30$ & $1.04$ & $0.04 $& $1.08$ \\
$(s_l)$   & $4537.95$ &$ 4474.84$ &$ 0.89$ & $0.04$ & $0.93$ \\
$(m_l,\sigma_l)$  & $4469.55 $ & $4457.63$ & $0.92$ & $0.03$ & $0.95$ \\
$(m_l,s_l)$ & $4443.89$ & $\bf{ 4434.26}$ & $0.90$ & $0.03$ & $0.93$ \\
$(\sigma_l,s_l)$  & $4458.68$ & $4465.55$ & $0.89$ & $0.03$ & $0.92$ \\
\bottomrule
\end{tabular}
\label{tab: compare input feature kmeans fsolve}
\end{table}

\begin{table}[ht!]
\centering
\caption{Comparison of different input features for K-Medoids}

\begin{tabular}{lccccc}
\toprule
{Feature} & {IS} & {OOS} & {CV Time (sec.)} & {Solve Time (sec.)} & {Total Time (sec.)} \\
\midrule
$(m_l,\sigma_l,s_l)$ & $4628.61  $ & $\bf{4486.70}$ & $0.56$ & $0.02$ & $0.58$ \\
$(m_l)$  & $  4546.78$ & $ 4492.57$ & $0.53$ & $0.02$ & $0.55$ \\
$(\sigma_l)$  & $ \bf{4513.23}$& $4514.86$ & $0.52$ & $0.02 $& $0.54$ \\
$(s_l)$   & $4640.32$ &$ 4511.92 $ &$  0.53$ & $0.02$ & $0.55$ \\
$(m_l,\sigma_l)$  & $4604.95  $ & $ 4495.80$ & $0.50$ & $0.02$ & $0.52$ \\
$(m_l,s_l)$ & $4559.03$ & $ 4496.78$ & $0.49$ & $0.02$ & $0.51$ \\
$(\sigma_l,s_l)$  & $4557.52$ & $4521.64$ & $0.50$ & $0.02$ & $0.52$ \\
\bottomrule
\end{tabular}
\label{tab: compare input feature kmedoids fsolve}
\end{table}

Tables\ \ref{tab: compare input feature kmeans fsolve} and\ \ref{tab: compare input feature kmedoids fsolve} show that the feature sets $(m_l,s_l)$ and $(m_l,\sigma_l,s_l)$ yield the lowest OOS cost for K-Means and K-Medoids-based heuristic algorithms, respectively. Moreover, based on Tables\ \ref{tab: compare input feature kmeans eq14} and \ref{tab: compare input feature kmedoids eq14} in Appendix\ \ref{app:rho =0.1 feature selection}, for $\rho=0.1$, K-Means returns the lowest OOS cost using \((m_l,s_l)\), and K-Medoids prefers the full feature set \((m_l,\sigma_l,s_l)\). These results confirm the managerial insights derived in Section \ref{sec:insight}, which suggest grouping patient types with similar average treatment times. For consistency purposes, we adopt these feature sets for our heuristic algorithms in the remainder of the paper.
Furthermore, both heuristic algorithms have a solve time within $0.05$ seconds and a CV time of about one second, demonstrating the computational efficiency.

\paragraph{Comparison of Different $\rho$.}\label{para:rho}
We vary the ambiguity set radius $\rho$ over the values $0,\ 0.1,\ 0.5,\ 1$. By construction, $\rho<1$ because it bounds the total variation distance between the unknown mode probabilities $p_l$ and their nominal references according to\ \eqref{eq:variation distance}. The table below reports the average IS costs, OOS costs, and runtime for $10$ independent runs.

\begin{table}[ht!]
\centering
\caption{Comparison for different $\rho$ for $L=5$}
\begin{tabular}{lcccc}
\toprule
$\rho$ & Method & IS & OOS & Runtime (sec.) \\
\midrule
$0$   & Algorithm\ \ref{algo:enumerate all groupings} & $\bf{4031.36}$ & $4659.24$ & $1.69$ \\
$0.1$ & Algorithm\ \ref{algo:eq14 algo}  & $4071.34$ & $4644.4$5 & $1.78$ \\
$0.5$ & Algorithm\ \ref{algo:eq14 algo} & $4185.46$ & $\bf{4582.08}$ & $1.56$ \\
$1$   & Algorithm\ \ref{algo:eq14 algo}  & $4262.82$ & $4683.53$ & $1.46$ \\
\bottomrule
\end{tabular}
\label{tab:rho comparison}
\end{table}

Table\ \ref{tab:rho comparison} shows that the IS cost rises monotonically with respect to the ambiguity set radius $\rho$. When $\rho=0$, the ambiguity set \eqref{eq:variation distance} collapses to the nominal mode probability $p_l = \frac{x_{lg}\hat{p}_l}{\sum_{l=1}^L x_{lg}\hat{p}_l}$ for each $l$. As $\rho$ increases, the ambiguity set \eqref{eq:variation distance} enlarges by including more candidate distributions. As a result, the decision maker becomes more conservative, leading to a higher IS cost. As $\rho$ increases, the OOS cost initially decreases, which demonstrates improved robustness to mode probability misspecification. However, when $\rho=1$, the OOS increases because the ambiguity set admits too many undesired probabilities and therefore is too conservative.

We put the detailed results (the optimal grouping assignments, optimal treatment time-slot durations, IS costs, OOS costs, and runtime) in Appendix\ \ref{app:one run rho} to investigate their changes with respect to $\rho$ for a fixed seed.




\paragraph{Effect of Asymmetry.}
Next, we investigate the effect of asymmetry in our model. Recall that this asymmetry is captured through the semivariance constraint in the ambiguity set $U_l$ defined in\ \eqref{eq:ul ambiguity set}. By removing the last semivariance constraint from this set and reformulating the resulting model, we obtain a counterpart, named \underline{M}ulti\underline{m}odal \underline{D}istributionally \underline{R}obust \underline{O}ptimization (MM-DRO), which does not consider the asymmetry in the distribution of $\tilde{\xi}$. In this case, the closed-form expression for $\Pi_l(t_g)$ in Theorem\ \ref{thm:thm1} does not hold. Therefore, we replace the support set $\xi\in\Xi$ in \eqref{eq:discrete change} with a discrete support $\{0,1,\ldots, 720\}$ and solve the discretized approximation of MMA-DRO and MM-DRO using Gurobi, respectively.

We adjust the skewness by varying the $\sigma$ parameter of the log-normal distribution, with a higher $\sigma$ denoting a higher skewness. 
We draw from a log-normal distribution with log-mean and log-standard deviation sampled from uniform ranges $\mu_l \sim \mathcal U(ln(60), ln(100))$ and
$\sigma_l \sim \mathcal U(1.5, 2.5)$ at the default setting and change the $\sigma$ range to examine its impact on the performance. We generate $100$ training samples and $1,000$ testing samples from this log-normal distribution. We report the average IS cost, OOS cost, and run time for $10$ independent runs to compare the performance of both models under varying skewness levels in Table\ \ref{tab:asy-comparison}.



\begin{table}[ht]
\centering
\caption{Comparison of different skewness under different models}
\begin{tabular}{lccccc}
\toprule
$\sigma_l$ & Model & IS & OOS & OOS Gap & Runtime (sec.)  \\
\midrule
$\mathcal U(1.5,2.5)$  & MMA-DRO & $\bf{3187.43}$ & $\bf{3625.39}$ & $-$ & $196.21$ \\
(Baseline)              & MM-DRO & $3588.14$ & $3765.48$ & $140.09$ & $314.81$ \\
\midrule
$\mathcal U(1.5,2.5)\times 0.8$  & MMA-DRO       & $\bf{3267.95}$ & $\bf{3537.10}$ & $-$ & $179.25$ \\
                        & MM-DRO & $3747.73$ & $3681.34$ & $144.24$ & $398.48$ \\
\midrule
$\mathcal U(1.5,2.5)\times 1.5$  & MMA-DRO        & $\bf{3025.15}$ & $\bf{3509.50}$ & $-$ & $213.44$ \\
                        & MM-DRO & $3463.25$ & $3881.54$ & $372.04$ & $235.33$ \\
\midrule
$\mathcal U(1.5,2.5)\times 2$    & MMA-DRO         & $\bf{2872.48}$ & $\bf{3384.79}$ & $-$ & $198.01$ \\
                        & MM-DRO & $3277.41$ & $3855.75$ & $470.96$ & $249.08$ \\
\bottomrule
\end{tabular}
\label{tab:asy-comparison}
\end{table}

As shown in Table\ \ref{tab:asy-comparison}, for a given skewness level, the skew-aware MMA-DRO produces both lower IS and OOS costs than MM-DRO. Moreover, as the log-standard deviation $\sigma$ increases, the gap of the OOS costs between MMA-DRO and MM-DRO becomes larger. These results highlight that the distributional asymmetry affects both the IS and OOS costs and that ignoring it reduces efficiency and accuracy.

\subsection{Performance Comparison between Different Approaches}\label{sec:performance-comparison}
\paragraph{Model Comparison.}
We now compare the performance of the exact Algorithm\ \ref{algo:enumerate all groupings}, the heuristic algorithms (K-Means and K-Medoids) with their best input features, and the discretized approximation solved by Gurobi when the number of modes $L$ ranges from $5$ to $8$ in Table \ref{tab:synthetic result}. When $L$ is small $(L=5,6)$, we report the results under two discretized approximation benchmarks: 
\textbf{(i)} solving the discretized version of Model \eqref{eq:continuous overall min} by treating $\boldsymbol{y}$ as a decision variable and
\textbf{(ii)} enumerating all possible grouping assignments $P\in\mathcal{P}$ and solving \eqref{eq:continuous overall min} with $\boldsymbol{y}$ fixed, both of which are solved by Gurobi directly. For the former, the time limit is set to 10 hours; for the latter, we set the time limit for each subproblem with a fixed $\boldsymbol{y}$ to one hour. For larger instances $(L=7,8)$, we only report the results under the enumeration approach, as treating $\boldsymbol{y}$ as a decision variable is computationally prohibitive within the time limit.

\begin{table}[ht!]
\centering
\caption{Performance comparison across different modes}
\scriptsize
\setlength{\tabcolsep}{4pt}
\renewcommand{\arraystretch}{1.2}
\begin{tabularx}{0.95\textwidth}{c >{\raggedright\arraybackslash}p{0.22\textwidth}
>{\centering\arraybackslash}p{0.10\textwidth}
>{\centering\arraybackslash}p{0.10\textwidth}
>{\centering\arraybackslash}p{0.17\textwidth}
>{\centering\arraybackslash}p{0.105\textwidth}
>{\centering\arraybackslash}p{0.105\textwidth}}
\toprule
\multirow[c]{2}{*}{$L$} &
\multirow[c]{2}{*}{Method} &
\multirow[c]{2}{*}{IS} &
\multirow[c]{2}{*}{OOS} &
\multirow[c]{2}{*}{Runtime (sec.)} &
\multicolumn{2}{c}{Relative error (\%)} \\
\cmidrule(lr){6-7}
&&&&& IS & OOS \\
\midrule
\multirow{5}{*}{5}
  & Gurobi + $y$ dec var. & $\bf{3997.22}$&$4684.64$ &$183.48$ & -$0.86$ & $0.54$\\
  & Gurobi + enumerate $y$ & $\bf{3997.22}$&$4684.64$ &$247.48$ & -$0.86$ & $0.54$\\
  & Algorithm\ \ref{algo:enumerate all groupings} & $4031.36$ & $4659.24$ & $1.69$ & $-$ & $-$\\
  & K-Means&$4443.89$ & $\bf{4434.26}$& CV: $0.89$, Solve: $0.03$ & $10.21$ & -$4.82$\\
  & K-Medoids &$4628.61$ & $ 4486.70$&CV: $0.51$, Solve: $0.02$ & $14.79$ & -$3.70$\\
\midrule
\multirow{5}{*}{6}
  & Gurobi + $y$ dec var.& $\bf{3808.05}$ & $4830.13$& $1966.07$& -$1.79$ & $0.16$\\
  & Gurobi + enumerate $y$ & $\bf{3808.05}$ & $4830.13$& $936.17$& -$1.79$ & $0.16$\\
  & Algorithm\ \ref{algo:enumerate all groupings} & $3877.58$ & $4822.42$ & $5.09$ & $-$ & $-$\\
  & K-Means &$ 4373.98$ &$\bf{4606.98} $ & CV: $1.30$, Solve: $0.04$ &$12.80$ & -$4.47$\\
  & K-Medoids &$4266.21$ &$4633.12$ &CV: $0.96$, Solve: $0.02$ &$10.02$ & -$3.93$\\
\midrule
\multirow{4}{*}{7}
  & Gurobi + enumerate $y$ & $\bf{3499.36}$& $4864.61$&$10996.48\ (\text{gap:}\ 1744.44\%)$ & -$1.25$ & $0.31$\\
  & Algorithm\ \ref{algo:enumerate all groupings} & $3543.62$ & $4849.74$ & $23.78$ & $-$ & $-$\\
  & K-Means &$4273.22$ & $4828.71$&CV: $1.61$, Solve: $0.04$& $20.58$ & -$0.43$\\
  & K-Medoids& $4325.10$& $\bf{4825.83}$& CV: $0.98$, Solve: $0.02$& $22.05$ & -$0.49$\\
\midrule
\multirow{4}{*}{8}
  & Gurobi + enumerate $y$     & $\bf{3638.13}$& $4811.49$&$20634.21\ (\text{gap:}\ 2943.98\%)$ & -$0.74$ & -$0.35$\\
  & Algorithm\ \ref{algo:enumerate all groupings} & $3665.21$ & $4828.39$ & $124.60$ & $-$ & $-$\\
  & K-Means&$4430.98$ & $\bf{4760.46}$&CV: $1.85$, Solve: $0.04$ &$20.89$ & -$1.41$\\
  &K-Medoids &$4555.65$ & $4762.56$& CV: $1.21$, Solve: $0.03$& $24.29$ & -$1.36$\\
\bottomrule
\end{tabularx}
\label{tab:synthetic result}
\end{table}

From Table\ \ref{tab:synthetic result}, for smaller instances ($L=5, 6$), both discretized approximation benchmarks are solved to optimality, which provides a lower bound to Model \eqref{eq:continuous overall min}. For larger instances (\(L=7,8\)), the Gurobi solver hits the time limit and returns only the incumbent with a nonzero optimality gap. On the other hand, the exact Algorithm \ref{algo:enumerate all groupings} can solve Model \eqref{eq:continuous overall min} to optimality more efficiently.
K-Means and K-Medoids, when running with their best input features, produce competitive IS and OOS costs while significantly reducing the runtime. As expected, the heuristic algorithms yield higher IS costs, since clustering does not guarantee optimal grouping assignments, with an IS relative error of $10\% \sim 24\%$. However, the heuristics achieve lower OOS costs with an OOS relative error of -$1\% \sim$ -$5\%$. Overall, these clustering-based heuristic approaches are highly efficient, with strong out-of-sample performance on large-scale instances.

\paragraph{Effect of Misspecified Distribution.}
In the previous sections, we generated the testing data from the same distribution as the training data. We now examine the impact of distributional misspecification of $\tilde{\xi}$. 
Specifically, for each patient type $l$, we perturb the log-mean and log-standard deviation by the perturbation parameter $\epsilon$ when generating the testing data, i.e.,
\[
\mu_l^{\mathrm{mis}} = (1+\epsilon)\,\mu_l
\qquad
\sigma_l^{\mathrm{mis}} = (1+\epsilon)\,\sigma_l
\]
with $\epsilon \in \{0.1,0.2,0.3,0.4,0.5\}$. To account for the misspecified distribution, we modify the ambiguity set $U_l$ to allow the mean, variance, and semivariance to deviate from the empirical ones within a given radius $\delta$ and reformulate the corresponding discretized MMA-$\rm{DRO^{mis}}$ model (full details are provided in Appendix\ \ref{app:mmadro_mis}).

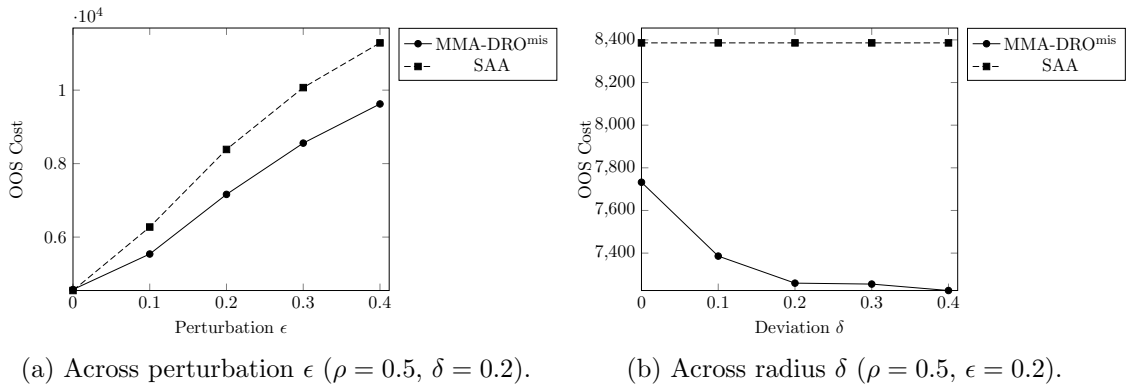
\begin{figure}[ht!]
  \centering
  \begin{subfigure}{0.45\textwidth}
    \centering
    \resizebox{\textwidth}{!}{
    \begin{tikzpicture}
      \begin{axis}[
        axis lines=box,
        xlabel={Perturbation $\epsilon$},
        ylabel={OOS Cost},
        xtick={0,0.1,0.2,0.3,0.4},
        xmin=0, xmax=0.4,
        enlarge x limits={upper, value=0.03},
        enlarge y limits={upper, value=0.06},
        cycle list name=black white,
        legend pos= outer north east
      ]
        \addplot+[mark=*, mark options={solid}]
          coordinates {(0.0,4572.47) (0.1,5539.12) (0.2,7162.40) (0.3,8557.98) (0.4,9623.93)};

        \addplot+[dash pattern=on 4pt off 2pt, mark=square*, mark options={solid}]
          coordinates {(0.0,4545.88) (0.1,6274.02) (0.2,8385.90) (0.3,10068.96) (0.4,11286.67)};
        \legend{MMA-DRO$^{\text{mis}}$, SAA}
      \end{axis}
    \end{tikzpicture}}
    \caption{Across perturbation $\epsilon$ ($\rho=0.5$, $\delta=0.2$).}
    \label{fig:misspecified-plot:eps}
  \end{subfigure}
  \begin{subfigure}{0.45\textwidth}
    \centering
    \resizebox{\textwidth}{!}{
    \begin{tikzpicture}
      \begin{axis}[
        axis lines=box,
        xlabel={Deviation $\delta$},
        ylabel={OOS Cost},
        xtick={0,0.1,0.2,0.3,0.4},
        xmin=0, xmax=0.4,
        enlarge x limits={upper, value=0.03},
        enlarge y limits={upper, value=0.06},
        cycle list name=black white,
        legend pos= outer north east
      ]
        \addplot+[mark=*, mark options={solid}]
          coordinates {(0.0,7732.31) (0.1,7386.06) (0.2,7259.20) (0.3,7254.81) (0.4,7224.63)};

        \addplot+[dash pattern=on 4pt off 2pt, mark=square*, mark options={solid}]
          coordinates {(0.0,8385.90) (0.1,8385.90) (0.2,8385.90) (0.3,8385.90) (0.4,8385.90)};
        \legend{MMA-DRO$^{\text{mis}}$, SAA}
      \end{axis}
    \end{tikzpicture}}
    \caption{Across radius $\delta$ ($\rho=0.5$, $\epsilon=0.2$).}
    \label{fig:misspecified-plot:delta}
  \end{subfigure}
  \caption{OOS comparisons of MMA-DRO$^{\text{mis}}$ and SAA.}
  \label{fig:misspecified-plot}
\end{figure}

We then compare MMA-$\rm{DRO^{mis}}$ with its SAA counterpart. Given a fixed $\rho=0.5$, the OOS cost comparison over $10$ independent runs of MMA-$\rm{DRO^{mis}}$ and SAA is shown in Figure\ \ref{fig:misspecified-plot}.

Figure\ \ref{fig:misspecified-plot:eps} shows that, for a fixed radius ($\delta=0.2$), the OOS cost increases with the perturbation level $\epsilon$ for both models. When $\epsilon=0$, SAA and MMA-$\rm{DRO^{mis}}$ achieve roughly the same OOS cost. 
However, as $\epsilon$ increases, MMA-$\rm{DRO^{mis}}$ achieves lower OOS costs compared to SAA, and the performance gap becomes larger.
On the other hand, Figure\ \ref{fig:misspecified-plot:delta} considers a fixed perturbation ($\epsilon=0.2$) and varies the radius $\delta$ in the $U_l$ ambiguity set from $0$ to $0.4$. As $\delta$ increases, the set $U_l$ admits a wider ambiguity in the distributions. Consequently, MMA-$\rm{DRO^{mis}}$ attains a lower OOS cost than SAA.

\section{Case Study: Mayo Clinic Chemotherapy Dataset}\label{sec: mayo section}
In this section, we evaluate the exact algorithms and proposed heuristics on the Mayo Clinic's real chemotherapy dataset. This dataset contains $1,851$ historical chemotherapy treatment duration records, and each treatment was scheduled into one of the seven predefined time‑slot templates (i.e., $30$-min, $60$-min, etc.). These templates define the seven patient types $l=1,\ldots,L$ in our model with $L=7$. In this section, we present 
(i) a comparison between the exact Algorithms\ \ref{algo:enumerate all groupings}, \ref{algo:eq14 algo}, the heuristic algorithms K-Means and K-Medoids, and the discretized approximation to solve the MMA-DRO model\ \eqref{eq:continuous overall min} in Section\ \ref{sec:rho=0 model comparison mayo}; (ii) a comparison across different radii of the mode probability ambiguity set in Section\ \ref{sec:rho>0 model comparison mayo}; (iii) a comparison between the MMA-DRO \eqref{eq:MMA-DRO} and its MM-DRO counterpart, which ignores distributional asymmetry in Section\ \ref{sec:effect of asy mayo}; and (iv) the real-world override comparison in Section\ \ref{sec:override}.

\subsection{Comparison of Different Models}\label{sec:rho=0 model comparison mayo}

We apply Algorithm\ \ref{algo:enumerate all groupings}, K-Means, K-Medoids, and the discretized approximation solved by Gurobi to compare the optimal treatment durations $t_g^*$ and the corresponding optimal cost. The comparative results are summarized in Table\ \ref{tab:mayo table}.

\begin{table}[ht!]
\centering
\caption{Comparison of Different Algorithms}
\resizebox{\textwidth}{!}{
\begin{tabular}{lcccc}
\toprule
Method & Optimal Durations (min.) & Optimal Cost & Runtime (sec.) & Relative Error ($\%$)\\ 
\midrule
  Gurobi + $y$ dec var. &$39,\ 60,\ 216$ &$\bf{1268.19}$&$36000 $ & $0$\\ 
 Gurobi + enumerate $y$&$39,\ 60,\ 216$ & $1268.19$&$31682$&$0$ \\ 
  Algorithm\ \ref{algo:enumerate all groupings} &$39,\ 60,\ 216$ &$1268.28$&$27.16$ & $-$\\ 
 K-Means&$58,\ 118,\ 217,\ 368$&$1409.32$&CV: $1.38$, Solve: $0.05$ &$11.12$ \\ 
K-Medoids&$39,\ 60,\ 118,\ 210,\ 261,\ 368$&$1503.47$&CV: $0.91$, Solve: $0.04$ &$18.54$ \\ \bottomrule
\end{tabular}}
\label{tab:mayo table}
\end{table}

For the discretized approximation where we treat $\boldsymbol{y}$ as a decision variable, we set a time limit of $10$ hours, within which Gurobi produced the best incumbent solution with a cost of $1268.19$ with an optimality gap of $1577116.52\%$. In contrast, the enumeration method achieved the same overall cost of $1268.19$ in $31682$ seconds ($\sim 8.8$ hours) with an optimality gap of $481246.74\%$. The exact Algorithm\ \ref{algo:enumerate all groupings} finds the same optimal assignment $(\boldsymbol{x}^*,\boldsymbol{y}^*)$ decisions and optimal time-slot durations with a cost of $1268.28$ more efficiently (within 30 seconds), where the duration $t_1^* = 39$ serves the $30$‑minute mode, duration $t_2^* = 60$ serves the $60$‑minute mode, and duration $t_3^* = 216$ jointly accommodates the $\{120,180,240,300,360\}$‑minute modes. 


On the other hand, K-Means and K-Medoids both produce near-optimal solutions in roughly two seconds, with the actual solve time around $0.05$ seconds. Using this approach, K-Means achieves a cost of
$1409.32$ with a relative error $11.12\%$ and K-Medoids yields a cost of
$1503.47$ with a relative error $18.54\%$. Specifically, K-Means generates four groups with a new time-slot template of $(58,118,217,368)$ minutes ($t_1^*=58$ to serve the $\{30,60\}$-minute types, $t_2^*=118$ to serve the $120$-minute type, $t_3^*=217$ to serve the $\{180,240,300\}$-minute types, and $t_4^*=368$ to serve the $360$-minute type, respectively), and K-Medoids yields $6$ groups with a new time-slot template of $(39,60,118,210,261,368)$ minutes.
These results demonstrate that our heuristic algorithms both achieve strong performance while significantly reducing computational time.

Next, we investigate how changing the size of the mode probability ambiguity set $\Delta_g\bigl(\hat{\boldsymbol{p}}, \boldsymbol{x},\boldsymbol{y})$ affects the overall model performance.

\subsection{Changing the radius of the ambiguity set $\rho$}\label{sec:rho>0 model comparison mayo}

We vary the robustness parameter $\rho \in \{0,0.1,0.5,1\}$ to enlarge the ambiguity set and solve the MMA-DRO model using the exact Algorithms\ \ref{algo:enumerate all groupings} and\ \ref{algo:eq14 algo} and the heuristic algorithms. We compare performance between the algorithms and the discretized approximation solved by the Gurobi enumeration method. The comparative results are summarized in Table\ \ref{tab:rho comparison for Mayo}.

\begin{table}[ht!]
\centering
\caption{Comparison for different $\rho$}

\resizebox{\textwidth}{!}{%
\begin{tabular}{c l c c c c}
\toprule
\multirow[c]{2}{*}{$\rho$} &
\multirow[c]{2}{*}{Method} &
\multirow[c]{2}{*}{Optimal duration (min.)} &
\multirow[c]{2}{*}{Optimal Cost} &
\multirow[c]{2}{*}{Runtime (sec.)} &
\multirow[c]{2}{*}{Relative Error (\%)} \\
&&&&& \\
\midrule
\multirow{4}{*}{$0$}
  & Gurobi                               & $39,\ 60,\ 216$                      & $\bf{1268.19}$ & $36000$           & $0$ \\
  & Algorithm~\ref{algo:enumerate all groupings} & $39,\ 60,\ 216$                      & $1268.28$ & $27.16$           & $-$ \\
  & K-Means                               & $58,\ 118,\ 217,\ 368$               & $1409.32$ & CV: $1.38$, Solve: $0.05$ & $11.12$ \\
  & K-Medoids                             & $39,\ 60,\ 118,\ 210,\ 261,\ 368$    & $1503.47$ & CV: $0.91$, Solve: $0.04$ & $18.54$ \\
\midrule
\multirow{4}{*}{$0.1$}
  & Gurobi                               & $39,\ 60,\ 221$                      & $\bf{1313.19}$ & $31056$           & $0$ \\
  & Algorithm~\ref{algo:eq14 algo}        & $39,\ 60,\ 221$                      & $1313.28$ & $43.65$           & $-$ \\
  & K-Means                               & $58,\ 118,\ 219,\ 368$               & $1419.94$ & CV: $2.79$, Solve: $0.50$ & $8.12$ \\
  & K-Medoids                             & $39,\ 60,\ 118,\ 210,\ 262,\ 368$    & $1505.64$ & CV: $2.35$, Solve: $0.47$ & $14.62$ \\
\midrule
\multirow{4}{*}{$0.5$}
  & Gurobi                               & $39,\ 60,\ 218,\ 368$                & $\bf{1386.86}$ & $33923$           & $0$ \\
  & Algorithm~\ref{algo:eq14 algo}        & $39,\ 60,\ 218,\ 368$                & $1386.93$ & $41.54$           & $-$ \\
  & K-Means                               & $59,\ 118,\ 225,\ 368$               & $1456.65$ & CV: $2.53$, Solve: $0.48$ & $5.03$ \\
  & K-Medoids                             & $39,\ 60,\ 118,\ 210,\ 278,\ 368$    & $1512.49$ & CV: $2.72$, Solve: $0.49$ & $9.05$ \\
\midrule
\multirow{4}{*}{$1$}
  & Gurobi                               & $39,\ 60,\ 238,\ 368$                & $\bf{1429.12}$ & $39132$           & $0$ \\
  & Algorithm~\ref{algo:eq14 algo}        & $39,\ 60,\ 238,\ 368$                & $1429.20$ & $41.98$           & $-$ \\
  & K-Means                               & $60,\ 118,\ 238,\ 368$               & $1486.36$ & CV: $2.53$, Solve: $0.49$ & $3.99$ \\
  & K-Medoids                             & $39,\ 60,\ 118,\ 210,\ 280,\ 368$    & $1512.52$ & CV: $2.72$, Solve: $0.49$ & $5.83$ \\
\bottomrule
\end{tabular}}
\label{tab:rho comparison for Mayo}
\end{table}

Table\ \ref{tab:rho comparison for Mayo} shows that the optimal cost increases monotonically with $\rho$, since a larger ambiguity set admits more adversarial mode probabilities, making the inner maximization problem to have higher costs and leading to more conservative decisions. For example, from $\rho=0.5$ to $\rho=1$, the grouping of the $7$ patient types remains unchanged, but the inner maximization shifts probability mass toward longer-duration modes, raising the optimal durations $t_3^*$ (e.g., the third group’s duration increases from $218$ to $238$). Also, from $\rho=0.1$ to $\rho=0.5$, Algorithm\ \ref{algo:eq14 algo} opens a fourth group to trade additional assignment cost for reduced overtime by isolating the $360$-minute mode. This result also illustrates how increasing the size of the ambiguity set $\rho$ values drives more conservative solutions.

We also observe that, the optimal costs returned by the exact Algorithms\ \ref{algo:enumerate all groupings} and\ \ref{algo:eq14 algo} are upper bounds on those from the discretized approximation solved by Gurobi, since the latter enforces only a subset of constraints. Moreover, both K-Means and K-Medoids algorithms perform well, achieving relative errors approximately between $4\%$ and $15\%$ while substantially reducing the computational time.

\subsection{Effect of Asymmetry}\label{sec:effect of asy mayo}
As described previously, once we remove the last equality constraint from the $U_l$ ambiguity set in\ \eqref{eq:ul ambiguity set}, an exact algorithm is no longer available. Therefore, we compare solutions from solving the discretized approximation model under two specifications: explicitly modeling the distributional asymmetry, and omitting it. We compare the optimal decisions $t_g^*$, the corresponding optimal cost, and the run time in Table\ \ref{tab:mayo asymmetry table}.
\begin{table}[ht!]
\centering
\caption{Effect of Asymmetry Comparison}
\resizebox{\textwidth}{!}{%
\begin{tabular}{lcccc}
\toprule
\textbf{Model} & \textbf{Optimal Durations (min.)} & \textbf{Optimal Cost} & \textbf{Runtime (sec.)} & \textbf{Relative Error (\%)} \\
\midrule
MMA-DRO
  & $39,\ 60,\ 216$ 
  & $\bf{1268.19}$ 
  & $31682\ (\text{gap:}\ 1577116.52\%)$ 
  & $-$ \\
MM-DRO 
  & $55,\ 81,\ 217$ 
  & $1448.63$ 
  & $36000\ (\text{gap:}\ 1136385.25\%)$ 
  & $14.23$ \\
\bottomrule
\end{tabular}%
}
\label{tab:mayo asymmetry table}
\end{table}

Table\ \ref{tab:mayo asymmetry table} shows that omitting the asymmetry by removing the semivariance constraint from the ambiguity set $U_l$ causes the optimal cost to increase by $14.23\%$ compared to the original MMA-DRO model. This gap demonstrates the value of including the asymmetry information that better hedges against the skewness in treatment duration distributions. 

\subsection{Real-world Override Comparison between the Optimized Template and the Current Template}\label{sec:override}
In chemotherapy scheduling, an override occurs when the allocated slots in a daily template are insufficient to accommodate the realized number of patients in a given category. This typically arises when patients call to schedule appointments, and the hospital operator assigns them to the most appropriate available time slot in the daily template. If no suitable slot exists within the same category as the patient type, the operator may split a longer slot or merge two shorter ones, resulting in an override.

Following the Mayo Clinic's current practice, we evaluate the number of overrides by calculating the difference between the allocated number of time slots and the realized number of patient counts in each patient type based on a three-month validation dataset. Specifically, Mayo Clinic’s current daily template contains $79$ time slots for 19 chairs in total, having a capacity of \(8{,}970\) minutes. Recall that the optimal time-slot durations produced by Algorithm \ref{algo:enumerate all groupings} in Table \ref{tab:mayo table} are \(39\), \(60\), and \(216\) minutes to accommodate the \(30\)-min, \(60\)-min, and \(120,180,240,300,360\)-min patient types, respectively. Using the same capacity, we allocate the \(8{,}970\) minutes to the three patient groups according to the percentage of these groups shown in the one-month training dataset, resulting in \((25,14,33)\) slots for \((39,60,216)\) minutes, respectively.
Then, on the three-month validation set, we compare the allocated time slots with the realized counts for each patient group on each day. When the template allocates fewer slots than the realized demand, we record an override. We apply the same procedure to Mayo Clinic’s current template to compute its override counts.

\begin{table}[ht!]
  \centering
  \caption{Comparison of overrides under the current and new templates.}
  \begin{subtable}[t]{0.48\linewidth}
    \centering
     \caption{Current template for 19 chairs}
    \begin{tabular}{lccc}
      \toprule
      & Duration (min.) & Slots & Overrides \\
      \midrule
     & $30$  & $27$ & $2$ \\
     & $60$  & $17$ & $0$ \\
      &$120$ & $11$ & $4$ \\
      &$180$ & $10$ & $27$ \\
      &$240$ & $7$  & $3$ \\
      &$300$ & $3$  & $1$ \\
      &$360$ & $4$  & $0$ \\
      \midrule
Total & 8,970 & 79 & 37 \\
      \bottomrule
    \end{tabular}
    \label{tab:mayo-template}
  \end{subtable}
  \hfill
  \begin{subtable}[t]{0.48\linewidth}
    \centering
    \caption{New template generated by MMA-DRO}
    \begin{tabular}{lccc}
      \toprule
      &Duration (min.) & Slots & Overrides \\
      \midrule
      &$39$  & $25$ & $15$ \\
      &$60$  & $14$ & $2$ \\
      &$216$ & $33$ & $0$ \\
      \midrule
      Total & 8,970 & 72 & 17\\
      \bottomrule
    \end{tabular}
    
    \label{tab:new-template}
  \end{subtable}
  \label{tab:template-comparison}
\end{table}

Table\ \ref{tab:template-comparison} compares override counts between Mayo Clinic’s current template and the optimized template, where the last row calculates the total duration, number of time slots, and number of overrides, respectively. The current template has a total of \(37\) overrides, whereas the new template only has \(17\) overrides, representing a \(54\%\) reduction in overrides. This decrease indicates that reallocating slots based on solving the MMA-DRO model\ \eqref{eq:MMA-DRO} achieves a better alignment between planned capacity and realized durations.

\section{Conclusion}\label{sec:conclusions}
In this paper, we developed a multimodal asymmetric DRO (MMA-DRO) model for chemotherapy template design that jointly determines the patient-type-to-grouping assignments and the corresponding optimal time-slot durations for each group. The ambiguity set captured both multi-modality across patient types and asymmetry in treatment-duration distributions. 
To handle the intractability of the resulting semi-infinite program, we proposed
\textbf{(i)} exact algorithms that enumerate all possible grouping assignments and solve for an optimal duration for each group based on the closed-form expression of the inner maximization problem;
\textbf{(ii)} heuristics that first cluster patient types via K-Means and K-Medoids to obtain the grouping assignments and then solve for the time-slot duration for each group; and
\textbf{(iii)} a discretized approximation solved by off-the-shelf solvers, in which the continuous support is replaced by a discrete support.

We evaluated the performance of the above algorithms on synthetic datasets and the Mayo Clinic's real chemotherapy dataset. Compared with the discretized approximation solved by Gurobi, our exact algorithms achieved optimality much more efficiently, and the heuristic approaches (K-Means and K-Medoids) yielded competitive IS and OOS costs while substantially reducing runtime on large instances. 
Sensitivity analyses further suggested that as the radius of the ambiguity set grows, the IS cost produced by the exact algorithms increases monotonically, whereas the OOS cost tends to decrease, indicating an improved robustness to distribution misspecification. Compared with the asymmetry-unaware counterpart MM-DRO, MMA–DRO consistently outperforms MM–DRO, which highlights the value of explicitly modeling asymmetry. 
We observed similar insights on the Mayo Clinic's real dataset. Furthermore, the optimized template generated by our MMA-DRO model significantly reduced overrides compared to the Mayo Clinic’s current template. 
In future work, we will extend our optimization framework to optimize the placement of time slots in a daily template that explicitly minimizes the downstream overrides. This paper mainly considered moment-based ambiguity sets to capture distributional uncertainty.  We will also consider distance-based ambiguity sets to incorporate the asymmetric information in the underlying treatment distributions.

\bibliographystyle{apalike}
\bibliography{bibl,refs}

\appendix
\section{Omitted Proofs}\label{sec:omitted proofs}
\subsection{Proof for Theorem\ \ref{thm: dual of Pi}}\label{pf:theorem 1}
\begin{proof}

We first expand the ambiguity set by writing out all moment constraints in $U_l$ as follows
\begin{align}\label{eq:ul ambiguity moment constrs}
\Pi_l(t_g) = \quad & \max_{f_{l}} \int_{\tilde{\xi} \in \Xi} f_{l}(\tilde{\xi}) \left(q(\tilde{\xi}-t_g)^+ + b(t_g-\tilde{\xi})^+\right) \nonumber\\
    \textrm{s.t.} \quad &  \int_{\tilde{\xi} \in \Xi} f_l(\tilde{\xi}) \, d\tilde{\xi}=1,\nonumber\\
    \quad & \int_{\tilde{\xi} \in \Xi} \tilde{\xi} f_l(\tilde{\xi}) \, d\tilde{\xi} = m_l,\nonumber \\
    \quad & \int_{\tilde{\xi} \in \Xi} (\tilde{\xi}-m_l)^2 f_l(\tilde{\xi}) \, d\tilde{\xi}= \sigma_l^2,\nonumber\\
    \quad & \int_{\tilde{\xi} \in \Xi} f_{l}(\tilde{\xi}) (\tilde{\xi}-m_l)^{+2} \, d\tilde{\xi} - \int_{\tilde{\xi} \in \Xi} f_{l}(\tilde{\xi}) (m_l-\tilde{\xi})^{+2} \, d\tilde{\xi} = s_l \sigma_l^2, \nonumber \\
    \quad & f_l(\tilde{\xi}) \geq 0,\ \forall \tilde{\xi} \in\Xi.
\end{align}
Since any feasible $t_g$ is bounded and $\Xi$ is a compact set, the objective function of\ \eqref{eq:ul ambiguity moment constrs} is bounded. Furthermore, according to Proposition 2.1 in \cite{natarajan2018asymmetry}, Assumption \ref{ass: feasible assumption} guarantees that the ambiguity set $U_l$ is nonempty. Following Theorem 1 in\ \cite{isii1962sharpness}, strong duality holds. Assigning dual variables $\gamma_{lg},\ \nu_{lg},\ \kappa_{lg},\ \zeta_{lg}$ to constraints in \eqref{eq:ul ambiguity moment constrs}, the dual formulation of\ \eqref{eq:ul ambiguity moment constrs} can then be written as
\begin{align}
    \Pi_l(t_g) = \min_{\gamma_{lg},\nu_{lg},\kappa_{lg},\zeta_{lg}} \quad & \gamma_{lg} + \nu_{lg} m_l + \kappa_{lg} \sigma_l^2 + \zeta_{lg} s_l \sigma_l^2 \nonumber \\
    \textrm{s.t.} \quad & \gamma_{lg} + \nu_{lg} \tilde{\xi} + \kappa_{lg}(\tilde{\xi}-m_l)^2 + \zeta_{lg} \left((\tilde{\xi}-m_l)^{+2} - (m_l-\tilde{\xi})^{+2}\right) \nonumber \\
    \quad & \geq q(\tilde{\xi}-t_g)^+ + b(t_g-\tilde{\xi})^+,\ \forall \tilde{\xi} \in \Xi \nonumber \\
    \quad & \gamma_{lg},\ \nu_{lg},\ \kappa_{lg},\ \zeta_{lg} \in \mathbb{R} \nonumber.
\end{align}
This completes the proof.
\end{proof}

\subsection{Proof for Theorem\ \ref{thm:dual of Omega}}\label{pf:theorem 2}
\begin{proof}
Expand $\Omega_g(t_g)$ under the ambiguity set $\Delta_g(\boldsymbol{\hat{p}},\boldsymbol{x}, \boldsymbol{y})$ defined in\ \eqref{eq:variation distance} as follows
\begin{subequations}\label{eq: pl ambiguity}
\begin{align}
    \max_{p_l} \quad & \sum_{l=1}^L p_l \Pi_l(t_g) \nonumber\\
    \textrm{s.t} \quad & \sum_{l=1}^L p_l = y_g,\label{eq: constr1}\\
    \quad & 0 \le p_l \leq x_{lg},\ \forall l \in \mathcal{L},\label{eq: constr2}\\
    \quad & \sum_{l=1}^L \left|p_l\Bigl(\sum_{l=1}^L x_{lg}\hat{p}_l\Bigr)-x_{lg}\hat{p}_l\right| \leq \rho \sum_{l=1}^L x_{lg}\hat{p}_l. \label{eq: constr3}
\end{align}
\end{subequations}

We have shown that $\Pi_l(t_g)$ is bounded in Theorem\ \ref{thm: dual of Pi} under Assumption\ \ref{ass: feasible assumption}. Moreover, for any feasible $(\boldsymbol{x}, \boldsymbol{y}, \boldsymbol{t})$, we can always find a feasible $\boldsymbol{p}$ that satisfies Problem\ \eqref{eq: pl ambiguity}. For instance, if group $g$ is not activated ($y_g = 0$), then no patient type is assigned to this group ($x_{lg}=0,\ \forall l\in \mathcal{L}$). As a result, $\boldsymbol{p}=(0,\ldots,0)^{\mathsf T}$ is a feasible solution to Problem\ \eqref{eq: pl ambiguity}. If group $g$ is activated ($y_g = 1$), then at least one patient type is assigned to group $g$. As a result, simply setting $p_l=\frac{x_{lg}\hat{p}_l}{\sum_{l=1}^L x_{lg}\hat{p}_l}$ for all $l\in\mathcal{L}$ yields a feasible solution.

To formulate the dual, we first linearize the absolute value term in Constraint\ \eqref{eq: constr3} by introducing an auxiliary nonnegative variable $d_{lg}$ as 
\begin{subequations}
\begin{align}
    \quad & d_{lg} \ge p_l \Bigl(\sum_{l=1}^L x_{lg}\hat{p}_l \Bigr) - x_{lg}\hat{p}_l,\label{eq:constr4} \\
    \quad & d_{lg} \ge x_{lg}\hat{p}_l - p_l \Bigl(\sum_{l=1}^L x_{lg}\hat{p}_l \Bigr),\label{eq:constr5}  \\
    \quad & \sum_{l=1}^L d_{lg} \le  \rho \sum_{l=1}^L x_{lg}\hat{p}_l, \label{eq:constr6}\\
    \quad & d_{lg} \ge 0. \label{eq:constr7}
\end{align}
\end{subequations}

Replacing Constraint\ \eqref{eq: constr3} with Constraints\ \eqref{eq:constr4}--\eqref{eq:constr7}, we obtain a linear programming (LP) formulation of the primal problem\ \eqref{eq: pl ambiguity}. Assigning dual variables $\mu_g,\ \lambda_{lg},\ \alpha_{lg},\ \beta_{lg}$, and $\tau_g$ to Constraint\ \eqref{eq: constr1},\ \eqref{eq: constr2},\ \eqref{eq:constr4},\ \eqref{eq:constr5}, and\ \eqref{eq:constr6}, respectively, we derive the dual formulation of Problem\ \eqref{eq: pl ambiguity} for group $g$ as follows
\begin{align}
\min_{\mu_g,\lambda_{lg},\tau_g, \alpha_l,\beta_l} \quad & \mu_g y_g + \sum_{l=1}^L \lambda_{lg} x_{lg} + \sum_{l=1}^L (\alpha_{lg} - \beta_{lg}) \hat{p}_lx_{lg} +  \rho \tau_g \sum_{l=1}^L x_{lg}\hat{p}_l\nonumber\\
    \textrm{s.t.} \quad & \mu_g + \lambda_{lg} + (\alpha_{lg}-\beta_{lg})(\sum_{l=1}^L x_{lg}\hat{p}_l) \geq \Pi_l(t_g),\ \forall l \in \mathcal{L}\nonumber\\
    \quad & -\alpha_{lg} - \beta_{lg} + \tau_g \geq 0,\ \forall l \in \mathcal{L} \nonumber \\
    \quad & \mu_g \in \mathbb{R},\ \lambda_{lg},\  \alpha_{lg},\ \beta_{lg},\ \tau_g \geq 0. \nonumber
\end{align}
Since this LP is feasible and bounded, we obtain the strong duality.
\end{proof}

\subsection{Proof for Theorem\ \ref{thm:thm1}}\label{proof: thm1}
\begin{proof}
To reformulate problem\ \eqref{eq:pl known inner}, we first rewrite the $U_l$ ambiguity as follows\ \eqref{eq:continuous ul ambiguity set}  as
\begin{align}\label{eq:new primal}
 \quad & \max_{f_{l}} \int_{\tilde{\xi} \in \Xi} f_{l}(\tilde{\xi}) \left(q(\tilde{\xi}-t_g)^+ + b(t_g-\tilde{\xi})^+\right) \nonumber\\
    \textrm{s.t.} \quad &  \int_{\tilde{\xi} \in \Xi} f_l(\tilde{\xi}) \, d\tilde{\xi}=1,\nonumber\\
    \quad & \int_{\tilde{\xi} \in \Xi} (\tilde{\xi} - m_l)^+ f_l(\tilde{\xi}) d\tilde{\xi} - \int_{\tilde{\xi} \in \Xi} (m_l - \tilde{\xi})^+ f_l(\tilde{\xi}) d\tilde{\xi}= 0,\nonumber \\
    \quad & \int_{\tilde{\xi} \in \Xi} f_{l}(\tilde{\xi}) (\tilde{\xi}-m_l)^{+2} \, d\tilde{\xi} = 
    \frac{(1+s_l)\sigma_l^2}{2}, \nonumber \\
    \quad & \int_{\tilde{\xi} \in \Xi} f_{l}(\tilde{\xi}) (m_l-\tilde{\xi})^{+2} \, d\tilde{\xi} = \frac{(1-s_l)\sigma_l^2}{2}, \nonumber \\
    \quad & f_l(\tilde{\xi}) \geq 0,\ \forall \tilde{\xi} \geq 0.
\end{align}
Assign dual variables $\alpha,\ \beta,\ y_1$ and $y_2$ to derive the dual formulation of\ \eqref{eq:new primal} to be
\begin{align}\label{eq:dual1}
     \min_{\alpha,\beta,y_1,y_2} \quad &\beta+\frac{(1+s_l)\sigma_l^2}{2}y_1 + \frac{(1-s_l)\sigma_l^2}{2}y_2 \nonumber\\
    \textrm{s.t} \quad & \beta+\alpha(x-m_l) + y_1(x-m_l)^2 \geq q (x-t_g)^+ + b (t_g-x)^+,\ \forall x \geq m_l\nonumber\\
    \quad & \beta+\alpha(x-m_l) + y_2(x-m_l)^2 \geq q (x-t_g)^+ + b (t_g-x)^+,\ \forall 0\leq x < m_l
\end{align}
The goal is to construct primal and dual feasible solutions with the same objective values and obtain optimal solutions. 
Let $f_1(x),\ g_1(x),\ f_2(x)$, and $g_2(x)$ denote the LHS and RHS of the first and second constraint, respectively. Let $w_1,\ w_2$ denote $\frac{(1+ s_l)\sigma_l^2}{2},\ \frac{(1- s_l)\sigma_l^2}{2}$, respectively. Checking feasibility is equivalent to verifying $f_1(x) \leq g_1(x)$ for all $x \geq m_l$, and $f_2(x) \leq g_2(x)$ for all $ 0 \leq x \leq m_l$, where the equalities are satisfied when the primal distribution has positive masses. The optimal cost is obtained by evaluating either the optimal primal or the optimal dual solution on its respective objective function.

\paragraph{Case $1$.} Suppose $m_l > t_g$. In this case, $g_1(x) = q(x-t_g)$ and \[
g_2(x)=
\begin{cases}
q(x-t_g), & t_g \le x < m_l,\\
b(t_g-x), & 0 \le x < t_g\\
\end{cases}
\]

\subparagraph{$(1a)$} Let $f_1(x)$ coincide with $g_1(x)$ for all $x \geq m_l$. Then on this domain we have
\[
\beta + \alpha(x - m_\ell) + y_1(x - m_\ell)^2 = q(x - t_g).
\]
Since $f_1$ is linear in this region, we set $y_1 = 0$. Then the matching coefficients are
$\alpha = q,\ \beta = q(m_l - t_g)$. Let $f_2(x)$ intersect with $g_2(x)$ on $x_1 = 0$ (See Figure\ \ref{fig:case1a}). Then we can get $y_2 = \frac{(b+q)t_g}{m_l^2}$. 

Since $f_1$ and $g_1$ coincide on infinitely many, but an unknown number of points, the primal distribution must consist of at least two support points. Suppose it is a two-point distribution. Denote $x_1$ to be the $x$-coordinate of the unique intersection of $f_1$ and $g_1$. Then the two-point primal distribution is
\[
\tilde{\xi} =
\begin{cases}
 0,\ \text{w.p.}\ p_1\\
 x_2,\ \text{w.p.}\ p_2,
\end{cases}
\]
which does not satisfy the primal feasibility condition in Eq.\ \eqref{eq:new primal}.

Consequently, denote $x_2>m_l,\ x_3>m_l$ to be the $x$-coordinates of the two intersection points and consider the following three-point distribution
\[
\tilde{\xi} =
\begin{cases}
 0,\ \text{w.p.}\ p_1\\
 x_2,\ \text{w.p.}\ p_2\\
 x_3,\ \text{w.p.}\ p_3,
\end{cases}
\]
and substitute it into Eq.\ \eqref{eq:new primal}. Since the system of equations contains five unknowns $
(x_2,\; x_3,\; p_1,\; p_2,\; p_3)$
and only four equalities, we introduce a parameter $\pi$ and express each variable in terms of $\pi$. As a result, this three-point primal distribution is 
\begin{align}
\tilde{\xi} = 
          \begin{cases}
  0 & \text{w.p.} \frac{w_2}{m^2}\\
  \frac{m_l}{m^2-w_2}(m^2-\sqrt{w_1(m^2-w_2)-w_2^2}\sqrt{\frac{1-\pi}{\pi}}) & \text{w.p.} \pi(1-\frac{w_2}{m^2}) \\
  \frac{m_l}{m^2-w_2}(m^2+\sqrt{w_1(m^2-w_2)-w_2^2}\sqrt{\frac{\pi}{1-\pi}}) & \text{w.p.} (1-\pi)(1-\frac{w_2}{m^2})
\end{cases}    
\end{align}
with $\pi \in \left[1-\frac{w_2^2}{w_1(m_l^2-w_2)},1\right)$. $t_g$ lies in the interval $[0,\frac{m_l}{2})$, ensuring that $f_2'(0) > -b$.

\subparagraph{$(1b)$}
Assume $f_1(x)$ coincide with $g_1(x)$ for all $x \geq m_l$. We have $y_1 = 0,\ \alpha = q$, and $\beta = q(m_l - t_g)$ from the previous. 
Now let $f_2(x)$ be tangent to $g_2(x)$ on $x \in (0,m_l)$ (See Figure\ \ref{fig:case1b}). Let $x_1$ denote the $x$-coordinate of the tangent point $x^*$, we can write out the system of equations for $f_2(x)$ and $g_2(x)$ as
\begin{align}
\begin{cases}
qm_l - qt_g - q(m-x_1) + y_2(m-x_1)^2 = b(t-x_1)\\
2y_2(x_1-m) + q = -b\\
0\leq x_1 \leq t_g
\end{cases}
\end{align}
Solving this equation gives us $x_1 = 2t_g - m_l,\ y_2 =\frac{b+q}{4(m_l-t_g)}$.

Suppose it is a two-point distribution. Let $x_2 \geq m_l$ denote the $x$-coordinate of the unique intersection point of $f_1$ and $g_1$. Then we can write
\[
\tilde{\xi} =
\begin{cases}
 2t_g - m_l, \text{w.p.}\ p_1\\
 x_2,  \text{w.p.}\ p_2,
\end{cases}
\]
which does not satisfy the primal feasibility condition in Eq.\ \eqref{eq:new primal}.

Consequently, denote $x_2>m_l,\ x_3>m_l$ to be the $x$-coordinates of the two intersection points and consider the following three-point distribution
\[
\tilde{\xi} =
\begin{cases}
  2t_g - m_l,\ \text{w.p.}\ p_1 \\
  x_2,\ \text{w.p.}\ p_2 \\
  x_3,\ \text{w.p.}\ p_3
\end{cases}
\]
and substitute it into Eq.\ \eqref{eq:new primal}. Using a similar approach and introducing the variable $\pi$, this three-point primal distribution is 
\[\tilde{\xi} = 
          \begin{cases}
  2t_g-m_l & \text{w.p.} \frac{w_2}{(2m_l-2t_g)^2}\\
  m_l+\frac{2(m_l-t_g)}{4(m_l-t_g)^2-w_2}(w_2-\sqrt{4w_1(m_l-t_g)^2-w_2(w_1+w_2)}\sqrt{\frac{1-\pi}{\pi}}) & \text{w.p.} \pi(1-\frac{w_2}{(2m_l-2t_g)^2}) \\
  m_l+\frac{2(m_l-t_g)}{4(m_l-t_g)^2-w_2}(w_2+\sqrt{4w_1(m_l-t_g)^2-w_2(w_1+w_2)}\sqrt{\frac{\pi}{1-\pi}}) & \text{w.p.} (1-\pi)(1-\frac{w_2}{(2m_l-2t_g)^2})
\end{cases}
\]
with $\pi \in [1-\frac{w_2^2}{w_1[4(m-t_g)^2-w_2]},1)$. $t_g$ lies in the interval $[\frac{m_l}{2}, m_l-\frac{1}{2}\sqrt{\frac{w_2(w_1+w_2)}{w_1}})$, ensuring that $f_2(0) \geq g_2(0)$.

\subparagraph{$(1c)$} Let $f_1(x)$ be tangent to $g_1(x)$ at $x_1$ and $f_2(x)$ be tangent to $g_2(x)$ at $x_2$, where
$
x_1 \geq m_l,\ 0 \leq x_2 \leq t_g$ and $y_1,y_2>0$ (See Figure\ \ref{fig:case1c}).
Then the primal admits a two–point distribution
\[
\tilde{\xi} =
\begin{cases}
x_1,\ \text{w.p}\ p_1\\
x_2,\ \text{w.p}\ p_2
\end{cases}
\]

Substituting into Eq.\ \eqref{eq:new primal} yields the system
\begin{align}\label{eq: 1c equation}
\begin{cases}
\beta + \alpha(x_1 - m_l) + y_1(x_1 - m_l)^2 = q(x_1 -t_g)\\
2y_1 ( x_1 - m_l) + \alpha = q \\
\beta + \alpha(x_2 - m_l) + y_2(x_2 - m_l)^2 =b(t_g - x_2)\\
2y_2 ( x_2 - m_l) + \alpha = -b
\end{cases}
\end{align}
which contains four variables and two equations. This system involves four variables and only two equations. Introducing a parameter $\pi$ and setting $p_2 = \pi,\ p_1 = 1 - \pi$ then substituting into Eq.\ \eqref{eq: 1c equation} gives $\pi = \frac{w_1}{w_1 + w_2}$. Hence, the optimal primal distribution is
\begin{align}
\tilde{\xi} = 
          \begin{cases}
  m_l+\sqrt{\frac{w_1(w_1+w_2)}{w_2}} & \text{w.p.} \frac{w_2}{w_1+w_2}\\
m_l-\sqrt{\frac{w_2(w_1+w_2)}{w_1}} & \text{w.p.} \frac{w_1}{w_1+w_2}\\
\end{cases}
\end{align}
Substituting this back into Eq.\ \eqref{eq: 1c equation}, we obtain the dual optimal solution
\begin{align}
\begin{cases}
    \quad & \beta = \frac{(m_l-t_g)(qw_2-bw_1)}{w_1+w_2}+\frac{b+q}{2}\sqrt{\frac{w_1w_2}{w_1+w_2}},\nonumber\\
\quad &\alpha=q-\frac{b+q}{w_1+w_2}[2(t_g-m_l)\sqrt{\frac{w_1w_2}{w_1+w_2}}-w_2], \nonumber\\
\quad & y_1=\frac{w_2(t_g-m_l)(b+q)}{(w_1+w_2)^2}+\frac{b+q}{2}\frac{w_2}{w_1+w_2}\frac{w_2}{w_1(w_1+w_2)},\nonumber\\
\quad & y_2=\frac{w_1(m_l-t_g)(b+q)}{(w_1+w_2)^2}+\frac{b+q}{2}\frac{w_1}{w_1+w_2}\frac{w_1}{w_2(w_1+w_2)},\nonumber
\end{cases}
\end{align}
and $t_g$ lies in the interval $[m_l - \frac{1}{2}\sqrt{\frac{w_2(w_1+w_2)}{w_1}},, m_l)$, ensuring that $y_1,\ y_2 > 0$.

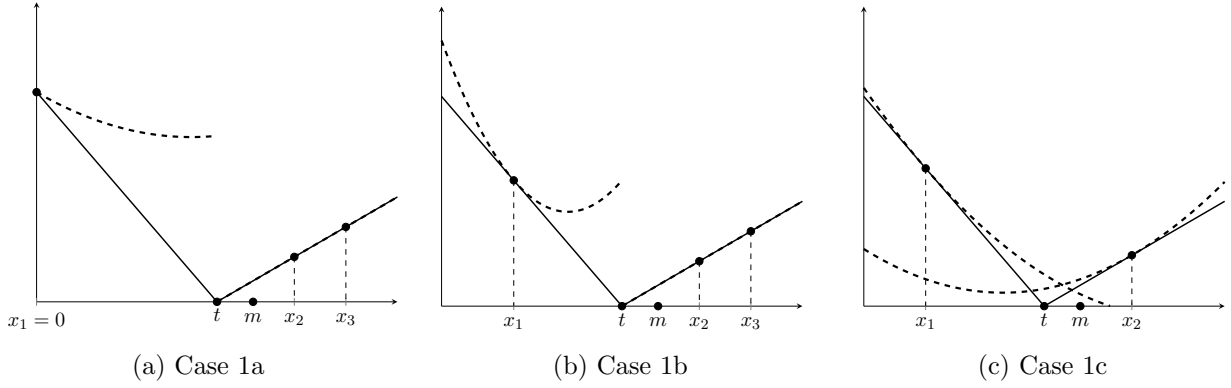
\begin{figure}[ht]
\centering

\begin{subfigure}{0.32\textwidth}
\centering
\begin{tikzpicture}[scale=0.7]
  \def\t{3.5}  
  \def\q{1.2}   
  \def\b{0.6}  

  \pgfmathsetmacro{\Cint}{\q*\t} 
  \def\Aint{0.10}                
  \def\Bint{-0.60}               

  \def\m{4.2}

  \begin{axis}[
      axis lines=middle,
      xmin=0, xmax=7,
      ymin=0, ymax=6,
      xtick=\empty, ytick=\empty,
      samples=300, clip=false,
      extra x ticks={0,\m,5.0,6.0},
      extra x tick labels={$x_1=0$,$m$,$x_2$,$x_3$}
    ]

    \addplot[thick, domain=0:\t] {\q*(\t - x)};
    \addplot[thick, domain=\t:7] {\b*(x - \t)};
    \addplot[very thick, dashed, domain=\t:7] {\b*(x - \t)};

    \addplot[only marks, mark=*] coordinates {(\t,0)};
    \node[below] at (axis cs:\t,0) {$t$};

    \addplot[very thick, dashed, domain=0:\t-0.05]{\Aint*x^2 + \Bint*x + \Cint};

    \addplot[only marks, mark=*, mark options={fill=black}] coordinates {(0,\Cint)};

    \addplot[only marks, mark=*] coordinates {(\m,0)};

    \def\xii{5.0}\def\xiii{6.0}
    \pgfmathsetmacro{\yii}{\b*(\xii-\t)}
    \pgfmathsetmacro{\yiiii}{\b*(\xiii-\t)}

    \addplot[only marks, mark=*] coordinates {(\xii,\yii)};
    \addplot[only marks, mark=*] coordinates {(\xiii,\yiiii)};

    \draw[dashed] (axis cs:\xii,\yii) -- (axis cs:\xii,0);
    \draw[dashed] (axis cs:\xiii,\yiiii) -- (axis cs:\xiii,0);

  \end{axis}
\end{tikzpicture}
\caption{Case 1a}
\label{fig:case1a}
\end{subfigure}
\hfill
\begin{subfigure}{0.32\textwidth}
\centering
\begin{tikzpicture}[scale=0.7]
  \def\t{3.5}  
  \def\q{1.2}   
  \def\b{0.6}  
  \def\a{1.4}  
  \pgfmathsetmacro{\ya}{\q*(\t-\a)}                
  \pgfmathsetmacro{\A}{\q/(\t-\a)}                 
  \pgfmathsetmacro{\xv}{\a + \q/(2*\A)}            
  \pgfmathsetmacro{\yv}{\A*(\xv-\a)^2 - \q*(\xv-\a) + \ya}

  \def\m{4.2}

  \begin{axis}[
      axis lines=middle,
      xmin=0, xmax=7,
      ymin=0, ymax=6,
      xtick=\empty, ytick=\empty,
      samples=300, clip=false,
      extra x ticks={1.4,\m,5.0,6.0},
      extra x tick labels={$x_1$,$m$,$x_2$,$x_3$}
    ]

    \addplot[thick, domain=0:\t] {\q*(\t - x)};
    \addplot[thick, domain=\t:7] {\b*(x - \t)};
    \addplot[very thick, dashed, domain=\t:7] {\b*(x - \t)};

    \addplot[only marks, mark=*] coordinates {(\t,0)};
    \node[below] at (axis cs:\t,0) {$t$};

    \addplot[very thick, dashed, domain=0:\t-0.05]{\A*(x-\a)^2 - \q*(x-\a) + \ya};

    \addplot[only marks, mark=*, mark options={fill=black}] coordinates {(\a,\ya)};
    \draw[dashed] (axis cs:\a,\ya) -- (axis cs:\a,0);

    \addplot[only marks, mark=*] coordinates {(\m,0)};

    \def\xii{5.0}\def\xiii{6.0}
    \pgfmathsetmacro{\yii}{\b*(\xii-\t)}
    \pgfmathsetmacro{\yiiii}{\b*(\xiii-\t)}

    \addplot[only marks, mark=*] coordinates {(\xii,\yii)};
    \addplot[only marks, mark=*] coordinates {(\xiii,\yiiii)};

    \draw[dashed] (axis cs:\xii,\yii) -- (axis cs:\xii,0);
    \draw[dashed] (axis cs:\xiii,\yiiii) -- (axis cs:\xiii,0);

  \end{axis}
\end{tikzpicture}
\caption{Case 1b}
\label{fig:case1b}
\end{subfigure}
\hfill
\begin{subfigure}{0.32\textwidth}
\centering
\begin{tikzpicture}[scale=0.7]

  \def\t{3.5}  
  \def\q{1.2}   
  \def\b{0.6}   
  \def\xone{1.2}  
  \def\xtwo{5.2}   
  \def\m{4.2}      

  \pgfmathsetmacro{\yone}{\q*(\t-\xone)}  
  \pgfmathsetmacro{\ytwo}{\b*(\xtwo-\t)}   
  \def\AL{0.12}
  \def\AR{0.12}

  \begin{axis}[
      axis lines=middle,
      xmin=0, xmax=7,
      ymin=0, ymax=6,
      xtick=\empty, ytick=\empty,
      samples=300, clip=false,
      extra x ticks={\xone,\m,\xtwo},
      extra x tick labels={$x_1$,$m$,$x_2$}
    ]

    \addplot[thick, domain=0:\t] {\q*(\t - x)};
    \addplot[thick, domain=\t:7] {\b*(x - \t)};

    \addplot[only marks, mark=*] coordinates {(\t,0)};
    \node[below] at (axis cs:\t,0) {$t$};
    \addplot[only marks, mark=*] coordinates {(\m,0)};

    \addplot[very thick, dashed, domain=0:7, restrict y to domain=0:6]
      {\AL*(x-\xone)^2 - \q*(x-\xone) + \yone};

    \addplot[very thick, dashed, domain=0:7, restrict y to domain=0:6]
      {\AR*(x-\xtwo)^2 + \b*(x-\xtwo) + \ytwo};
    \addplot[only marks, mark=*, mark options={fill=black}] coordinates {(\xone,\yone)};
    \addplot[only marks, mark=*, mark options={fill=black}] coordinates {(\xtwo,\ytwo)};

    \draw[dashed] (axis cs:\xone,\yone) -- (axis cs:\xone,0);
    \draw[dashed] (axis cs:\xtwo,\ytwo) -- (axis cs:\xtwo,0);

  \end{axis}
\end{tikzpicture}
\caption{Case 1c}
\label{fig:case1c}
\end{subfigure}
\caption{Case 1 $m_l \geq t_g$ illustrations}
\end{figure}

\paragraph{Case 2.} Suppose $m_l \leq t_g$.  In this case, $g_2(x) = q(x-t_g)$ and \[
g_1(x)=
\begin{cases}
q(x-t_g), & x \geq t_g\\
b(t_g-x), & m_l < x \leq t_g\\
\end{cases}
\]
\subparagraph{$(2a)$} Assume $f_2$ coincide with $g_2$ on $[0,m_l)$. On this domain, we have
\[
\beta + \alpha(x-m_l) + y_2(x-m_l)^2 = b(t_g - x).
\]
Since $f_2$ is linear in this region, we set $y_2 = 0$. Then the matching coefficients are $\alpha = -b,\ \beta = bt_g - bm_l$. Let $f_1$ be tangent to $g_1$ on $x_1 > t_g$, where $x_1$ denotes the $x$-coordinate of the tangent point. In this case, we obtain the system of equations
\[
\begin{cases}
    bt_g - bm_l - b(x_1 - m_l) + y_1(x_1 - m_l)^2 = q(x_1 - t_g)\\
    2y_1(x_1 - m_l) - b  = q
\end{cases}
\]
Solving for $y_1 = \frac{b+q}{4(t_g - m_l)}$ and $x_1 = 2t_g - m_l$. Similar to case $(1b)$, we know that the primal optimal distribution $\tilde{\xi}$ is a three-point distribution
\begin{align}
    \tilde{\xi} = 
          \begin{cases}
2t_g-m_l & \text{w.p.} \frac{w_1}{(2m_l-2t)^2}\nonumber\\
m_l-\frac{2(t_g-m_l)}{4(t_g-m_l)^2-w_1}(w_1-\sqrt{4w_2(t_g-m_l)^2-w_1(w_1+w_2)}\sqrt{\frac{1-\pi}{\pi}}) & \text{w.p.} \pi(1-\frac{w_1}{(2m_l-2t_g)^2}) \nonumber\\
m_l-\frac{2(t_g-m_l)}{4(t_g-m_l)^2-w_1}(w_1+\sqrt{4w_2(t_g-m_l)^2-w_1(w_1+w_2)}\sqrt{\frac{\pi}{1-\pi}}) & \text{w.p.} (1-\pi)(1-\frac{w_1}{(2m_l-2t_g)^2}), \nonumber
\end{cases}
\end{align}
with the auxiliary variable $\pi \in \left[1-\frac{w_2^2}{w_1(4(m_l-t_g)^2-w_2)},1\right)$. $t_g$ lies in the interval $\bigl[m_l+\frac{1}{2}\sqrt{\frac{w_2(w_1+w_2)}{w_1}}, m_l + \frac{m_lw_1}{2w_2}\bigr)$ ensuring that $x_2 > x_3$ and that the square root expressions are well defined.

\subparagraph{$(2b)$} Assume $f_1$ is tangent to $g_1$ on $x_1 > t_g$ and $f_2$ is tangent to $g_2$ on $x_2 \in [0,m_l)$, where $x_1,\ x_2$ are the $x$-coordinates of the tangent points. In this case, we obtain the system of equations
\[
\begin{cases}
    \beta + \alpha(x_1 - m_l) + y_1(x_1-m_l)^2 = q(x_1 - t_g)\\
    2y_1(x_1-m_l) + \alpha = q\\
    \beta + \alpha(x_2 - m_l) + y_2(x_2-m_l)^2 = b(t_g- x_2)\\
    2y_2(x_2-m_l) + \alpha = -b,
\end{cases}
\]
with the corresponding two-point primal optimal distribution
\[
\tilde{\xi} = 
\begin{cases}
    x_1,\ \text{w.p.}\ p_1\\
    x_2,\ \text{w.p.}\ p_2.
\end{cases}
\]
This is exactly the same as case $(1c)$. Solving for variables $(x_1,\ x_2,\ p_1,\ p_2,\ y_1,\ y_2,\ \alpha,\ \beta)$, we obtain the same result. $t_g$ lies in the interval of $ \bigl[m_l, m_l+\frac{1}{2}\sqrt{\frac{w_2(w_1+w_2)}{w_1}}\bigr)$ ensuring $y_1,\ y_2>0$.

\subparagraph{$(2c)$} Let $f_2$ be concave and intersect with $g_2$ at $x_1 =0$, where $x_1$ is the $x$-coordinate of the intersection point. Let $f_1$ be tangent to $g_1$ at $m_l\leq x_2 <t_g$ and $x_3 \geq t_g$, where $x_2,\ x_3$ are the $x$-coordinates of the two tangent points. In this case, the primal optimal distribution is a three-point distribution 
\[
\tilde{\xi}=
\begin{cases}
    0,\ \text{w.p.}\ p_1 \\
    x_2,\ \text{w.p.}\ p_2 \\
    x_3,\ \text{w.p.}\ p_3, \\
\end{cases}
\]
and the system of equations is
\[
\begin{cases}
\beta - \alpha m_l + y_2m_l^2 = bt_g \\
\beta + \alpha(x_2 - m_l) + y_1(x_2 - m_l)^2 = b(t_g - x_2)\\
2y_1 (x_2 - m_l) + \alpha = -b \\
\beta + \alpha(x_3-m_l) +y_1(x_3 - m_l)^2 = q(x_3 - t_g) \\
2y_1 (x_3 - m_l) + \alpha = q
\end{cases}
\]
Together with the primal feasibility condition\ \eqref{eq:new primal}, we obtain nine equations in nine variables, which return a unique solution for the primal optimal distribution and the corresponding optimal dual solutions. $t_g^*$ lies in the range of $[m_l + \frac{m_l w_1}{2w_2}, \infty)$ ensuring the concavity of $f_1$.

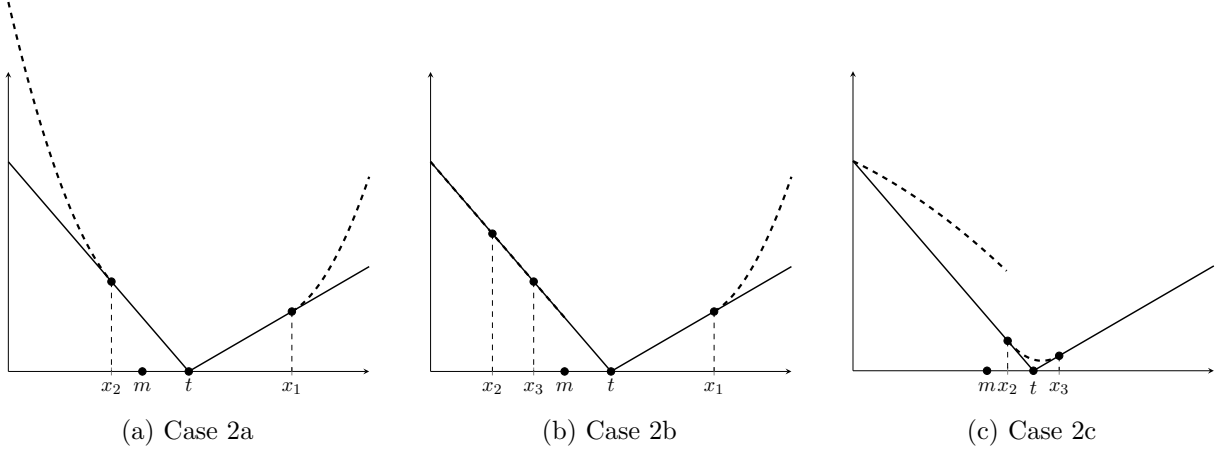
\begin{figure}[ht]
\centering

\begin{subfigure}{0.32\textwidth}
\centering
\begin{tikzpicture}[scale=0.7]
  \def\t{3.5}
  \def\q{1.2}  
  \def\b{0.6}  

  \def\xii{2.0}  
  \def\m{2.6}    
  \def\xi{5.5}   

  \pgfmathsetmacro{\yleft}{\q*(\t-\xii)}   
  \pgfmathsetmacro{\yright}{\b*(\xi-\t)}   

  \def\AL{0.8}   
  \def\AR{0.8}   

  \begin{axis}[
      axis lines=middle,
      xmin=0, xmax=7,
      ymin=0, ymax=6,
      xtick=\empty, ytick=\empty,
      samples=300, clip=false,
      extra x ticks={\xii,\m,\xi},
      extra x tick labels={$x_{2}$,$m$,$x_{1}$}
    ]

    \addplot[thick, domain=0:\t]   {\q*(\t - x)};  
    \addplot[thick, domain=\t:7]   {\b*(x - \t)};  

    \addplot[only marks, mark=*] coordinates {(\t,0)};
    \node[below] at (axis cs:\t,0) {$t$};

    \addplot[only marks, mark=*] coordinates {(\m,0)};

    \addplot[only marks, mark=*] coordinates {(\xii,\yleft)};
    \draw[dashed] (axis cs:\xii,\yleft) -- (axis cs:\xii,0);

    \addplot[only marks, mark=*] coordinates {(\xi,\yright)};
    \draw[dashed] (axis cs:\xi,\yright) -- (axis cs:\xi,0);

    \addplot[very thick, dashed, domain=0:\xii]
      {\AL*(x-\xii)^2 - \q*(x-\xii) + \yleft};

    \addplot[very thick, dashed, domain=\xi:7]
      {\AR*(x-\xi)^2 + \b*(x-\xi) + \yright};

  \end{axis}
\end{tikzpicture}
\caption{Case 2a}
\label{fig:case2a}
\end{subfigure}
\hfill
\begin{subfigure}{0.32\textwidth}
\centering
\begin{tikzpicture}[scale=0.7]
  \def\t{3.5}
  \def\q{1.2}
  \def\b{0.6}

  \def\xii{1.2}   
  \def\xiii{2.0}  
  \def\m{2.6}     
  \def\xi{5.5}    

  \pgfmathsetmacro{\yii}{\q*(\t-\xii)}
  \pgfmathsetmacro{\yiiii}{\q*(\t-\xiii)}
  \pgfmathsetmacro{\yi}{\b*(\xi-\t)}  

  \def\Aint{0.80}

  \begin{axis}[
      axis lines=middle,
      xmin=0, xmax=7,
      ymin=0, ymax=6,
      xtick=\empty, ytick=\empty,
      samples=300, clip=false,
      extra x ticks={\xii,\xiii,\m,\xi},
      extra x tick labels={$x_2$,$x_3$,$m$,$x_1$}
    ]

    \addplot[thick, domain=0:\t] {\q*(\t - x)};
    \addplot[thick, domain=\t:7] {\b*(x - \t)};

    \addplot[very thick, dashed, domain=0:\m] {\q*(\t - x)};

    \addplot[only marks, mark=*] coordinates {(\t,0)};
    \node[below] at (axis cs:\t,0) {$t$};

    \addplot[only marks, mark=*] coordinates {(\m,0)};

    \addplot[only marks, mark=*] coordinates {(\xii,\yii)};
    \draw[dashed] (axis cs:\xii,\yii) -- (axis cs:\xii,0);

    \addplot[only marks, mark=*] coordinates {(\xiii,\yiiii)};
    \draw[dashed] (axis cs:\xiii,\yiiii) -- (axis cs:\xiii,0);

    \addplot[only marks, mark=*] coordinates {(\xi,\yi)};
    \draw[dashed] (axis cs:\xi,\yi) -- (axis cs:\xi,0);

    \addplot[very thick, dashed, domain=\xi:7]
      {\Aint*(x-\xi)^2 + \b*(x-\xi) + \yi};

  \end{axis}
\end{tikzpicture}
\caption{Case 2b}
\label{fig:case2b}
\end{subfigure}
\hfill
\begin{subfigure}{0.32\textwidth}
\centering
\begin{tikzpicture}[scale=0.7]
  \def\t{3.5}
  \def\q{1.2}   
  \def\b{0.6}   

  \def\m{2.6}
  \def\xii{3.0}                 
  \pgfmathsetmacro{\xiii}{2*\t - \xii} 

  \pgfmathsetmacro{\yii}{\q*(\t-\xii)}     
  \pgfmathsetmacro{\yiiii}{\b*(\xiii-\t)}  
  \pgfmathsetmacro{\sii}{-1.0*\q}
  \pgfmathsetmacro{\siii}{\b}

  \pgfmathsetmacro{\Av}{(\siii - \sii)/(2*(\xiii - \xii))}
  \pgfmathsetmacro{\Bv}{\sii - 2*\Av*\xii}
  \pgfmathsetmacro{\Cv}{\yii - \Av*\xii*\xii - \Bv*\xii}

  \def\Ac{-0.08}
  \def\Bc{-0.5}
  \pgfmathsetmacro{\Cc}{\q*\t}

  \def\eps{0.02} 
  \pgfmathsetmacro{\xleftend}{\xii - \eps}

  \begin{axis}[
      axis lines=middle,
      xmin=0, xmax=7,
      ymin=0, ymax=6,
      xtick=\empty, ytick=\empty,
      samples=300, clip=false,
      extra x ticks={\xii,\m,\t,\xiii},
      extra x tick labels={$x_2$,$m$,$t$,$x_3$}
    ]

    \addplot[thick, domain=0:\t] {\q*(\t - x)};
    \addplot[thick, domain=\t:7] {\b*(x - \t)};

    \addplot[only marks, mark=*] coordinates {(\t,0)};
    \addplot[only marks, mark=*] coordinates {(\m,0)};

    \addplot[only marks, mark=*] coordinates {(\xii,\yii)};
    \draw[dashed] (axis cs:\xii,\yii) -- (axis cs:\xii,0);
    \addplot[only marks, mark=*] coordinates {(\xiii,\yiiii)};
    \draw[dashed] (axis cs:\xiii,\yiiii) -- (axis cs:\xiii,0);

    \addplot[very thick, dashed, domain=0:\xleftend]
      {\Ac*x^2 + \Bc*x + \Cc};

    \addplot[very thick, dashed, domain=\xii:\xiii]
      {\Av*x^2 + \Bv*x + \Cv};

  \end{axis}
\end{tikzpicture}
\caption{Case 2c}
\label{fig:case2c}
\end{subfigure}
\caption{Case 2 $m_l \leq t_g$ illustrations}
\end{figure}

\paragraph{Case 1. $m_l \geq t_g$} In this case, $g_1(x)=q(x-t_g)$, $g_2(x)= q(x-t_g)$ for $t_g \leq x \leq m_l$, and $b(t_g-x)$ for $ 0 \leq x \leq t_g$.

\begin{table}[ht!]
    \centering
    \begin{tabular}{|l|}
    \hline
          Case 1(a): $0 \leq t_g \leq \frac{m_l}{2}$.\\
          Primal feasible distribution: \\
          $\tilde{\xi} = 
          \begin{cases}
  0 & \text{w.p.} \frac{w_2}{m^2}\\
  \frac{m_l}{m^2-w_2}(m^2-\sqrt{w_1(m^2-w_2)-w_2^2}\sqrt{\frac{1-\pi}{\pi}}) & \text{w.p.} \pi(1-\frac{w_2}{m^2}) \\
  \frac{m_l}{m^2-w_2}(m^2+\sqrt{w_1(m^2-w_2)-w_2^2}\sqrt{\frac{\pi}{1-\pi}}) & \text{w.p.} (1-\pi)(1-\frac{w_2}{m^2})
\end{cases}$\\
where $\pi \in \left[1-\frac{w_2^2}{w_1(m_l^2-w_2)},1\right)$.\\
Dual feasible solution: $\beta= qm_l-qt_g,\ \alpha=q,\ y_1=0,\  y_2=\frac{(b+q)t_g}{m_l^2}$.\\
Optimal Cost: $qm_l-qt_g + w_2 \frac{(b+q)t_g}{m_l^2}= (\frac{(b+q)(1-s_l)\sigma_l^2}{2m_l^2}-q)t_g + qm_l$.\\
\hline
          Case 1(b): $\frac{m_l}{2} \leq t_g \leq m_l-\frac{1}{2}\sqrt{\frac{w_2(w_1+w_2)}{w_1}}$.\\
          Primal feasible distribution: \\
          $\tilde{\xi} = 
          \begin{cases}
  2t_g-m_l & \text{w.p.} \frac{w_2}{(2m_l-2t_g)^2}\\
  m_l+\frac{2(m_l-t_g)}{4(m_l-t_g)^2-w_2}(w_2-\sqrt{4w_1(m_l-t_g)^2-w_2(w_1+w_2)}\sqrt{\frac{1-\pi}{\pi}}) & \text{w.p.} \pi(1-\frac{w_2}{(2m_l-2t_g)^2}) \\
  m_l+\frac{2(m_l-t_g)}{4(m_l-t_g)^2-w_2}(w_2+\sqrt{4w_1(m_l-t_g)^2-w_2(w_1+w_2)}\sqrt{\frac{\pi}{1-\pi}}) & \text{w.p.} (1-\pi)(1-\frac{w_2}{(2m_l-2t_g)^2})
\end{cases}$\\
where $\pi \in [1-\frac{w_2^2}{w_1[4(m-t_g)^2-w_2]},1)$.\\
Dual feasible solution: $\beta= qm_l-qt_g,\ \alpha=q,\ y_1=0,\  y_2=\frac{(b+q)}{4(m_l-t_g)}$.\\
Optimal Cost: $qm_l-qt_g + w_2 \frac{b+q}{4(m_l-t_g)} = \frac{(b+q)(1-s_l)\sigma_l^2}{8(m_l-t_g)}-qt_g+qm_l$.\\


\hline
          Case 1(c): $m_l-\frac{1}{2}\sqrt{\frac{w_2(w_1+w_2)}{w_1}} \leq t_g \leq m_l$.\\
          Primal feasible distribution: \\
          $\tilde{\xi} = 
          \begin{cases}
  m_l+\sqrt{\frac{w_1(w_1+w_2)}{w_2}} & \text{w.p.} \frac{w_2}{w_1+w_2}\\
m_l-\sqrt{\frac{w_2(w_1+w_2)}{w_1}} & \text{w.p.} \frac{w_1}{w_1+w_2}\\
\end{cases}$\\
Dual feasible solution: \\
$\beta = \frac{(m_l-t_g)(qw_2-bw_1)}{w_1+w_2}+\frac{b+q}{2}\sqrt{\frac{w_1w_2}{w_1+w_2}}$,\\
$\alpha=q-\frac{b+q}{w_1+w_2}[2(t_g-m_l)\sqrt{\frac{w_1w_2}{w_1+w_2}}-w_2]$,\\
$y_1=\frac{w_2(t_g-m_l)(b+q)}{(w_1+w_2)^2}+\frac{b+q}{2}\frac{w_2}{w_1+w_2}\frac{w_2}{w_1(w_1+w_2)}$,\\
$y_2=\frac{w_1(m_l-t_g)(b+q)}{(w_1+w_2)^2}+\frac{b+q}{2}\frac{w_1}{w_1+w_2}\frac{w_1}{w_2(w_1+w_2)}$,\\

Optimal Cost: $(b+q)\sqrt{\frac{w_1w_2}{w_1+w_2}}+(m_l-t_g)\frac{qw_2-bw_1}{w_1+w_2}=\frac{(b+q)\sigma_l}{2}\sqrt{1-s_l^2} + (m_l-t_g)\frac{(q-b)-(q+b)s_l}{2}$.\\
\hline
    \end{tabular}
    \caption{Primal and dual optimal solution when $m_l \geq t_g$}
    \label{case1}
\end{table}

\paragraph{Case 2. $m_l < t_g$} In this case, $g_2(x)=b(t_g-x)$, and $g_1(x) = b(t_g-x)$ for $m_l \leq x \leq t_g$, and $q(x-t_g)$ for $x \geq t_g$.
\begin{table}[h!]
    \centering
    \begin{tabular}{|l|}
\hline
    Case 2(a) $ m_l \leq t_g \leq m_l+\frac{1}{2}\sqrt{\frac{w_2(w_1+w_2)}{w_1}}$.\\
    Primal feasible distribution: \\
          $\tilde{\xi} = 
          \begin{cases}
  m_l+\sqrt{\frac{w_1(w_1+w_2)}{w_2}} & \text{w.p.} \frac{w_2}{w_1+w_2}\\
m_l-\sqrt{\frac{w_2(w_1+w_2)}{w_1}} & \text{w.p.} \frac{w_1}{w_1+w_2}\\
\end{cases}$\\
Dual feasible solution: \\
$\beta = \frac{(m_l-t_g)(qw_2-bw_1)}{w_1+w_2}+\frac{b+q}{2}\sqrt{\frac{w_1w_2}{w_1+w_2}}$,\\
$\alpha=q-\frac{b+q}{w_1+w_2}[2(t_g-m_l)\sqrt{\frac{w_1w_2}{w_1+w_2}}-w_2]$,\\
$y_1=\frac{w_2(t_g-m_l)(b+q)}{(w_1+w_2)^2}+\frac{b+q}{2}\frac{w_2}{w_1+w_2}\frac{w_2}{w_1(w_1+w_2)}$,\\
$y_2=\frac{w_1(m_l-t_g)(b+q)}{(w_1+w_2)^2}+\frac{b+q}{2}\frac{w_1}{w_1+w_2}\frac{w_1}{w_2(w_1+w_2)}$,\\

Optimal Cost: $(b+q)\sqrt{\frac{w_1w_2}{w_1+w_2}}+(t_g-m_l)\frac{bw_1-qw_2}{w_1+w_2}=\frac{(b+q)\sigma_l}{2}\sqrt{1-s_l^2} + (m_l-t_g)\frac{(q-b)-(q+b)s_l}{2}$.\\
\hline
    Case 2(b): $m_l+\frac{1}{2}\sqrt{\frac{w_2(w_1+w_2)}{w_1}} \leq t_g \leq m_l + \frac{m_lw_1}{2w_2}$\\
          Primal feasible distribution: \\
          $\tilde{\xi} = 
          \begin{cases}
  2t_g-m_l & \text{w.p.} \frac{w_1}{(2m_l-2t)^2}\\
  m_l-\frac{2(t_g-m_l)}{4(t_g-m_l)^2-w_1}(w_1-\sqrt{4w_2(t_g-m_l)^2-w_1(w_1+w_2)}\sqrt{\frac{1-\pi}{\pi}}) & \text{w.p.} \pi(1-\frac{w_1}{(2m_l-2t_g)^2}) \\
  m_l-\frac{2(t_g-m_l)}{4(t_g-m_l)^2-w_1}(w_1+\sqrt{4w_2(t_g-m_l)^2-w_1(w_1+w_2)}\sqrt{\frac{\pi}{1-\pi}}) & \text{w.p.} (1-\pi)(1-\frac{w_1}{(2m_l-2t_g)^2})
\end{cases}$\\

where $\pi \in \left[1-\frac{w_2^2}{w_1(4(m_l-t_g)^2-w_2)},1\right)$.\\
Dual feasible solution: $\beta= bt_g-bm_l,\ \alpha=-b,\ y_1=\frac{(b+q)}{4(t_g-m_l)},\  y_2=0$\\
Optimal Cost: $bt_g-bm_l + w_1 \frac{b+q}{4(t_g-m_l)} = \frac{(b+q)(1+s_l)\sigma_l^2}{8(t_g-m_l)}+bt_g-bm_l$.\\
\hline
Case 2(c) $m_l + \frac{m_lw_1}{2w_2} \leq t_g.$\\
Primal feasible distribution: \\
          $\tilde{\xi} = 
          \begin{cases}
  0 & \text{w.p.} \frac{w_2}{m_l^2}\\
  t_g - \sqrt{\frac{(t_g-m_l)^2(m_l^2 - w_2) + m_l^2 w_1 - 2m_l w_2 (t_g-m_l)}{\,m_l^2 - w_2\,}} & \text{w.p.} 1-\frac{w_2}{m_l^2}-\frac{w_2^2}{m_l^2w_1} \\
  t_g + \sqrt{\frac{(t_g-m_l)^2(m_l^2 - w_2) + m_l^2 w_1 - 2m_l w_2 (t_g-m_l)}{\,m_l^2 - w_2\,}} & \text{w.p.} \frac{w_2^2}{m^2w_1}
\end{cases}$\\

Dual feasible solution: \\
$\beta
= \frac{t_g-m_l}{2} \Big[
\bigl(q - 3b\bigr)\sqrt{\tfrac{2A+B}{B}}
- 3(b+q)\sqrt{\tfrac{B}{2A+B}} - \bigl(6b + 4q\bigr)
\Big]$ \\
$\alpha = b + (b+q)\sqrt{\frac{B}{2A+B}},$ \\
$y_1 = \frac{b+q}{2(t_g - m_l)}\sqrt{\frac{B}{2A+B}},$ \\
$y_2
=\frac{1}{2m_l^{2}}\Big[
4\big(t_g(2b+q)-m_l(b+q)\big)
+(m_l-t_g)(q-3b)\sqrt{\tfrac{2A+B}{B}}
+(3t_g-m_l)(b+q)\sqrt{\tfrac{B}{2A+B}}
\Big]$\\
where $A := \frac{w_{1}}{2(t_g-m_l)^{2}} \;-\; \frac{w_{2}}{m_l(t_g-m)},\ B = 1 - \frac{w_2}{m_l^2}$\\
Optimal Cost: $bt_g+qm_l-\frac{b+q}{2}(m_l+(1-\frac{w_2}{m_l^2}) t_g-(t_g-m_l)\sqrt{(1-\frac{w_2}{m_l^2}) ((1-\frac{w_2}{m_l^2})+\frac{m_lw_1-2w_2(t_g-m_l)}{m_l(t_g-m_l)^2})})$.\\
\hline

    \end{tabular}
    \caption{Primal and dual optimal solution when $m_l \leq t$}
    \label{case2}
\end{table}
\end{proof}

\subsection{Proof for Theorem\ \ref{thm:convexity}}
\begin{proof}
First, Problem\ \eqref{eq:heuristic problem} becomes a convex problem since the objective function $\Pi(t_g)=\sum_{l=1}^L p_l\,\Pi_l(t_g)$ is a convex combination of convex functions and is therefore convex in $t_g$. 
By the optimality condition for convex programs with a box constraint $\mathcal{T}=[0,T]$,
$t_g^* \in \arg\min_{t_g \in \mathcal{T}} \Pi(t_g)$ if and only if $0 \in \partial_{t_g} \Pi(t_g^*)+ \mathcal{N}_{\mathcal{T}}(t_g^*)$, where the normal cone of $\mathcal{T}$ is defined as follows
\[
\mathcal{N}_{\mathcal{T}}(t_g)=
\begin{cases}
\{0\}, & t_g \in (0,T),\\
(-\infty,0], & t_g=0,\\
[0,\infty), & t_g=T.
\end{cases}
\]

Because each piece $f_{li}(t_g)$ is differentiable on its interval $\mathcal T_{li}$ for each mode $l=1,\ldots,L$ and $i=1,\ldots,5$, $\Pi(t_g)$ is also differentiable on each of the intervals $\mathcal{T}_j$ for $j=1,\ldots,J$. If there exists $t_g^* \in int(\mathcal{T}_{j})$ that satisfies
$\sum_{l:\,x_{lg}=1} p_l\,\Pi_l'(t_g^*)=0$ (i.e., a stationary point of a convex function piece),
then at this point, $\mathcal{N}_{\mathcal{T}}(t_g^*)=\{0\}$ and $0 \in \partial_{t_g} \Pi(t_g^*)+ \mathcal{N}_{\mathcal{T}}(t_g^*)$. As a result, $t_g^*$ is optimal. If no such interior point exists on any $\mathcal T_j$, the minimizer must be attained at an endpoint of the intervals, i.e., $t_g^*=0$ if the one-sided gradient $\nabla_{t_g}\Pi(0)\ge 0$, $t_g=T$ if the one-sided gradient $\nabla_{t_g}\Pi(T)\le 0$, and $t_g^*$ is achieved at an endpoint other than $0$ and $T$ if $0\in \partial_{t_g} \Pi(t_g^*)$. 


\end{proof}
\subsection{Proof for Proposition\ \ref{prop:parameter analysis pieces}}
\begin{proof}
Recall that MMA-DRO\ \eqref{eq:MMA-DRO} is convex in $t_g$.  Since $\Pi(t_g) = \sum_{l:x_{lg}=1} p_l \Pi_l(t_g)$ is
bounded and convex for group $g \in \mathcal{G}$, the optimizer is determined by \textbf{(i)} the sign of the slope on the linear pieces and \textbf{(ii)} whether the
nonlinear pieces admit interior stationary points. 

\begin{enumerate}
\item \emph{Initial decrease on $\mathcal T_{1}$ (first piece).} 
Define the common first-piece domain
$\mathcal T_{1}\ :=\ \bigcap_{\,l:\,x_{lg}=1}\ \mathcal T_{l1}=\Bigl[0, \min_{\,l:\,x_{lg}=1}\tfrac{m_l}{2}\,\Bigr].$
On $\mathcal T_{1}$, we have 
\[
\Pi(t_g)=\sum_{l:\,x_{lg}=1} p_l\,f_{l1}(t_g)=\Bigl(\sum_{l:\,x_{lg}=1} p_l\,k_l\Bigr)\,t_g\ +q\sum_{l:\,x_{lg}=1} p_lm_l,\] 
which is linear in $t_g$ with $k_l:=\frac{(b+q)(1-s_l)\sigma_l^2}{2 m_l^2}-q$.
Therefore, the gradient of $\Pi(t_g)$ on $\mathcal{T}_1$ is
\[
\sum_{l:x_{lg}=1} p_l k_l = \sum_{l:x_{lg}=1} p_l\left( \frac{(b+q) (1-s_l)\sigma_l^2}{2m_l^2}\,-q\right).
\]

When $\sum_{l:x_{lg}=1} p_l (1-s_l)\left(\frac{\sigma_l}{m_l}\right)^2 \;\ge\; \frac{2q}{\,b+q\,}$, then the objective function $\Pi(t_g)$ is nondecreasing in $t_g$ for $t_g \in \mathcal{T}_1$. Because $\Pi(t_g)$ is convex in $t_g$, it must be nondecreasing for all other pieces. As a result, the global minimizer is attained at the left endpoint of $\mathcal{T}_1$, i.e., $t_g^*=0$.

\item \emph{Last decrease on $\mathcal T_{5}$ (last piece).} 
Define the common last-piece domain
$\mathcal T_{5}\ :=\ \bigcap_{\,l:\,x_{lg}=1}\ \mathcal T_{l5}=\Bigl[\max_{l:x_{lg}=1} m_l + \frac{m_l(1+s_l)}{2(1-s_l)}, T \Bigr].$ 
On $\mathcal{T}_5$, we have
\[
\Pi(t)= \sum_{\,l:\,x_{lg}=1} p_l\,f_{l5}^{(l)}(t_g),
\quad\text{and}\ \partial\Pi(t_g) \;=\; \sum_{\,l:\,x_{lg}=1} p_l\,\partial f_{l5}(t_g).
\]
The subgradients of $\Pi(t_g)$ are monotone nondecreasing in $t_g$ when $t_g \in \mathcal{T}_5$. Since $\partial \Pi(t_g) \in \bigl[\Pi_-'(t_g), \Pi_+'(t_g)\bigr]$, and the largest subgradient is attained at the right endpoint $T$.
\[
\sup_{t_g \in \mathcal{T}_5} \Pi_+'(t_g) = \Pi_+'(T^-)\le 0
\]

Denote $\Delta_{t_g,l}=t_g-m_l>0$. Differentiating $f_{l5}$ with respect to $t_g$ and evaluate at $T$ gives:
\[
f'_{l5,+}(T^-)
= b-\frac{b+q}{2}\!\left[
\beta_l-\sqrt{S_l(T)}
+\frac{\beta_l}{\sqrt{S_l(T)}}
\left(\frac{(1+s_l)\sigma_l^2}{2\,\Delta_{T,l}^{2}}
-\frac{(1-s_l)\sigma_l^2}{2\,m_l\,\Delta_{T,l}}\right)
\right],
\]
where $S_l(T)=\beta_l\!\left(\beta_l+\frac{(1+s_l)\sigma_l^2}{2\,\Delta_{T,l}^{2}}
-\frac{(1-s_l)\sigma_l^2}{m_l\,\Delta_{T,l}}\right)$.
Consequently, $\Pi(t_g)$ is nonincreasing on $\mathcal{T}_5$ iff \(\sum_{l:x_{lg}=1} p_l f'_{l5,+}(T^-)\le 0\), i.e.
\begin{align}\label{eq:T5 inequality}
\sum_{l:x_{lg}=1} p_l
b-\sum_{l:x_{lg}=1} p_l\frac{b+q}{2}\!\left[
\beta_l-\sqrt{S_l(T)}
+\frac{\beta_l}{\sqrt{S_l(T)}}\!\left(
\frac{(1+s_l)\sigma_l^2}{2\,\Delta_{T,l}^{2}}
-\frac{(1-s_l)\sigma_l^2}{2\,m_l\,\Delta_{T,l}}
\right)
\right]\ \le\ 0.
\end{align}
Denote $H_l(T) :=\ \beta_l-\sqrt{S_l(T)}
+\frac{\beta_l}{\sqrt{S_l(T)}}\!\left(
\frac{(1+s_l)\sigma_l^2}{2\,\Delta_{l}^{2}}
-\frac{(1-s_l)\sigma_l^2}{2\,m_l\,\Delta_{l}}
\right)$. Then Eq.\ \eqref{eq:T5 inequality} is equivalent to \begin{align}\label{eq: T5 result}
    \frac{2b}{b+q} \le \sum_{l:x_{lg}=1} p_l H_l(T)
\end{align}
\end{enumerate}
This completes the proof.
\end{proof}

\subsection{Proof for Corollary\ \ref{coro:parameter}}
\begin{proof}
From Proposition\ \ref{prop:parameter analysis pieces}, we know 
\[
\sum_{l:\,x_{lg}=1} p_l(1-s_l)\!\left(\frac{\sigma_l}{m_l}\right)^{\!2} \;\ge\; \frac{2q}{\,b+q\,}
\quad\Longrightarrow\quad t_g^*=0,
\]
and
\[
\sum_{l:\,x_{lg}=1} p_lH_l(T) \;\ge\; \frac{2b}{\,b+q\,}
\quad\Longrightarrow\quad t_g^*=T .
\]
When $q=0$ or $b \rightarrow +\infty$, $\frac{2q}{b+q} \rightarrow0$ and $\sum_{l:\,x_{lg}=1} p_l(1-s_l)\!\left(\frac{\sigma_l}{m_l}\right)^{\!2}$ is always nonnegative under Assumption\ \ref{ass: feasible assumption}. As a result, the model chooses $t_g^*=0$ to be the optimal solution. 

Now we show $H_l(T) \ge 0$. Under the Assumption\ \ref{ass: feasible assumption}, $\beta_l > 0$. Denote $\Delta = T - m_l$ and
\[
A_l(\Delta):=\frac{w_{1l}}{\Delta_{l,T}^2}-\frac{2w_{2l}}{m_l\Delta_{l,T}},
\qquad
B_l(\Delta):=\frac{w_{1l}}{\Delta_{l,T}^2}-\frac{w_{2l}}{m_l\Delta_{l,T}},
\]
and we can write \(S_l(T)=\beta_l\big(\beta_l+A_l(\Delta)\big)\).

Furthermore, on $\mathcal{T}_5$, we have
\(\Delta\ge m_l + \frac{m_l(1+s_l)}{2(1-s_l)} - m_l =
\frac{m_l(1+s_l)}{2(1-s_l)}:= \Delta_l^{\min}\). Since \(A_l(\Delta_l^{\min})=\frac{2\sigma_l^2(1-s_l)^2}{m_l^2 (1+s_l)}- \frac{2\sigma_l^2(1-s_l)^2}{m_l^2 (1+s_l)}= 0\) and \(A_l'(\Delta)<0\) (nonincreasing), obtaining the minimum value $A_l(\Delta^*) = -\frac{(1-s_l)^2 \sigma_l^2}{2(1+s_l)m_l^2}$ at $\Delta^* = \frac{m_l(1+s_1)}{1-s_l} = 2\Delta_l^{\min}$. Therefore, \[\beta_l + A_l(\Delta) \ge \beta_l - \frac{(1-s_l)^2 \sigma_l^2}{2(1+s_l)m_l^2} = 1 - \frac{(1-s_l)\sigma_l^2}{m_l^2 ( 1+s_l)} \ge 0
\]
under the Assumption\ \ref{ass: feasible assumption}. As a result, $S_l(T) \ge 0$.

Based on the analysis, we can get
\[
S_l(T)=\beta_l\big(\beta_l+A_l(\Delta)\big)\le \beta_l^2
\quad\Rightarrow\quad
\sqrt{S_l(T)}\le \beta_l .
\]

Multiplying \(H_l(T)\) by \(\sqrt{S_l(T)}{>0}\) and using \(S_l(T)=\beta_l(\beta_l+A_l(\Delta)\) we have
\[
\begin{aligned}
H_l(T)\ge 0
&\iff \beta_l\sqrt{S_l(T)}-S_l(T)+\beta_l B_l(\Delta)\ \ge\ 0\\
&\iff \sqrt{S_l(T)}-(\beta_l+A_l(\Delta))+B_l(\Delta)\ \ge\ 0\\
&\iff B_l(T)\ \ge\ (\beta_l+A_l(\Delta))-\sqrt{S_l(T)}.
\end{aligned}
\]
Since \(A_l(\Delta)\le 0\) and \(\sqrt{S_l(T)}\le \beta_l\),
\[
(\beta_l+A_l(\Delta))-\sqrt{S_l(T)}
\ \ge\ (\beta_l+A_l(\Delta))-\beta_l \;=\; A_l(\Delta).
\]
As a result, it suffices to show \(B_l(\Delta)\ge A_l(\Delta)\). Notice that
\[
B_l(\Delta)-A_l(\Delta)=\frac{w_{2l}}{m_l\Delta_{l,T}}\ \ge\ 0,
\]
since all $w_{2l},\ m_l$, and $\Delta_{l,T}$ are all nonnegative, which indicates \(H_l(T)\ge 0\) for all $t_g\in\mathcal{T}_5$.

Therefore, Eq.\ \eqref{eq: T5 result} can be satisfied as $q \rightarrow +\infty$ or $b=0$, which both drive $\frac{2b}{b+q}$ to $0$ and in such cases, $t_g^* = T$ is optimal.





\end{proof}

\subsection{Proof for Proposition\ \ref{prop:lb}}\label{app: lower bound pf}
\begin{proof}\label{pf: lower bound pf}
For any $t_g\in [0,T]$, denote $\phi_{t_g}(\tilde{\xi})\triangleq q(\tilde{\xi}-t_g)^+ + b(t_g-\tilde{\xi})^+$, which is convex in $\tilde{\xi}$. By Jensen’s inequality, we have
\[
\Pi_l(t_g)=\max_{f_l\in U_l}\;\mathbb E_{f_l}\bigl[\phi_{t_g}(\tilde\xi)\bigr]
\;\ge\; \max_{f_l\in U_l}\; \phi_{t_g}\!\bigl(\mathbb E_{f_l}[\tilde\xi]\bigr)
= q(m_l-t_g)^+ + b(t_g-m_l)^+ \eqqcolon \Phi_{t_g}(m_l).
\]

Hence
\[
\sum_{l:\,x_{lg}=1} p_l \Phi_{t_g}(m_l)
\;\ge\;
b\sum_{l:\,x_{lg}=1} p_l (t_g-m_l)^+ \;+\;
q\sum_{l:\,x_{lg}=1} p_l (m_l-t_g)^+ .
\]

Since all terms are nonnegative, keeping only the extreme means yields
\[
\sum_{l:\,x_{lg}=1} p_l \Phi_{t_g}(m_l)
\ \ge\ b \sum_{\substack{l:\,x_{lg}=1\\ m_l=m_{\min}(g)}}\!\! p_l\,\bigl(t_g-m_{\min}(g)\bigr)^+
\;+\;
q \sum_{\substack{l:\,x_{lg}=1\\ m_l=m_{\max}(g)}}\!\! p_l\,\bigl(m_{\max}(g)-t_g\bigr)^+ .
\]
For any $t_g$, maximizing over $\boldsymbol p$ on both sides gives two valid lower bounds,
\[
\max_{\boldsymbol p}\sum_{l:\,x_{lg}=1} p_l \Phi_{t_g}(m_l)\ \ge\ b\,\overline p_{\min}(g)\,(t_g-m_{\min}(g))^+,
\
\max_{\boldsymbol p}\sum_{l:\,x_{lg}=1} p_l \Phi_{t_g}(m_l)\ \ge\ q\,\overline p_{\max}(g)\,(m_{\max}(g)-t_g)^+ .
\]
Therefore, for every $t_g$,
\begin{equation}\label{eq:intermediate}
\max_{\boldsymbol p}\sum_{l:\,x_{lg}=1} p_l \Phi_{t_g}(m_l)
\ \ge \max\!\Big\{\,b\,\overline p_{\min}(g)\,(t_g-m_{\min}(g))^+,\ q\,\overline p_{\max}(g)\,(m_{\max}(g)-t_g)^+ \Big\}.
\end{equation}

Denote the right-hand side of\ \eqref{eq:intermediate} as $F(t_g)$, and minimize $F(t_g)$ over $t_g \in \mathcal{T}$. If $t_g<m_{\min}(g) < m_{\max}(g)$ then
$F(t_g)=q\,\overline p_{\max}(g)\,(m_{\max}(g)-t_g)\ge q\,\overline p_{\max}(g)\,(m_{\max}-m_{\min})=F(m_{\min}(g))$, which shows that $F(t_g)$ is decreasing on $t_g \in (0,m_{\min}(g))$.
Similarly, if $t_g>m_{\max}(g)>m_{\min}(g)$ then
$F(t_g)=b\,\overline p_{\min}(g)\,(t_g-m_{\min}(g))\ge b\,\overline p_{\min}(g)\,(m_{\max}-m_{\min})=F(m_{\max}(g))$,
which shows that $F(t_g)$ is increasing on $t_g \in (m_{\max}(g),T]$.
Thus, the minimizer lies in the range $[m_{\min}(g),m_{\max}(g)]$.

On this interval, $F(t_g)$ is the maximum of two affine functions:
\[
F(t_g)=\max\big\{\,b\,\overline p_{\min}(g)\,(t_g-m_{\min}),\ q\,\overline p_{\max}(g)\,(m_{\max}-t_g)\big\}.
\]
The minimum of the maximum of these two affine functions on a closed interval is attained at the intersection point of these two lines if such a point lies within the interval. Solving
\[
b\,\overline p_{\min}(g)\,(t_g-m_{\min})= q\,\overline p_{\max}(g)\,(m_{\max}-t_g)
\]
gives us \[t_g^* =  \frac{b\,\overline p_{\min}(g) m_{\min}(g) + q\,\overline p_{\max}(g) m_{\max}(g)}{b\,\overline p_{\min}(g) + q\,\overline p_{\max}(g)}\in [m_{\min}(g),m_{\max}(g)].\]

Therefore, at this point, 
\[
\min_{t_g \in [m_{\min}(g),m_{\max}(g)]} F(t_g) = F(t_g^*) =  \frac{b\,\overline p_{\min}(g) q\,\overline p_{\max}(g)}{b\,\overline p_{\min}(g)+q\,\overline p_{\max}(g)} (m_{\max}(g)- m_{\min}(g))
\]

As a result, the problem\ \eqref{eq:heuristic problem} is bounded below by $\min_{t_g}F(t_g)
=\frac{b\,\overline p_{\min}(g) q\,\overline p_{\max}(g)}{b\,\overline p_{\min}(g)+q\,\overline p_{\max}(g)} (m_{\max}(g)- m_{\min}(g))$.
\end{proof}

\subsection{Proof for Proposition\ \ref{prop:upper bound}}\label{app:ub proof}
\begin{proof}\label{pf:ub proof}
Suppose we are given $t_g\in\mathcal{T}$. For any $\tilde\xi$,
$\,q(\tilde\xi-t_g)^+ + b(t_g-\tilde\xi)^+\le \max\{b,q\}\,\lvert \tilde\xi - t_g\rvert$.
Therefore, for each $l$,
    \begin{align}
        \Pi_l(t_g) & \le \max\{b,q\} \max_{f_l \in U_l} \mathbb{E}[\bigl| \tilde{\xi} - t_g\bigr|] \nonumber \\
        & \le \max\{b,q\}\left(\max_{f_l \in U_l}\mathbb{E}[|\tilde{\xi} - m_l|] + |m_l - t_g|\right) \qquad\text{(triangle inequality in $L^1$)} \nonumber \\
    & \le \max\{b,q\} \Bigl(\max_{f_l \in U_l}\sqrt{\mathbb{E}[\tilde{\xi} - m_l]^2} + |m_l - t_g| \Bigr) \qquad\text{(Cauchy--Schwarz)}\nonumber \\
   & \le \max\{b,q\} \Bigl(\sigma_l + |m_l - t_g|\Bigr) \nonumber.
  \end{align}

Therefore, for each group $g$,
\begin{align}\label{eq:UB_step}
\sum_{l:\,x_{lg}=1} p_l\,\Pi_l(t_g)
&\le \max\{b,q\}\left(\sum_{l:\,x_{lg}=1} p_l\sigma_l + \sum_{l:\,x_{lg}=1} p_l\lvert m_l - t_g\rvert\right)\nonumber\\
&\le \max\{b,q\}\,y_g\left(\sigma_{\max}(g) + \max_{\,l:\,x_{lg}=1}\lvert m_l - t_g\rvert\right).
\end{align}

We then take the minimum over $t_g$. Suppose $t_g \le m_{\min}(g)$, we have $\max_{l:x_{lg}=1} |m_l-t_g| = m_{\max}(g) - t_g \ge m_{\max}(g) - m_{\min}(g)$. Similarly, suppose $t_g \ge m_{\max}(g)$, we have $\max_{l:x_{lg}=1} |m_l-t_g| =  t_g - m_{\min}(g) \ge m_{\max}(g) - m_{\min}(g)$. Therefore, the right-hand side of\ \eqref{eq:UB_step} achieves the minimum when $t_g \in [m_{\min}(g),m_{\max}(g)]$. Within this interval, $\min_{t_g \in [m_{\min}(g),m_{\max}(g)]}\max_{l:x_{lg}=1} |m_l - t_g| = \max\{t_g - m_{\min}(g), m_{\max}(g) - t_g\}$, and the minimizer is obtained when two terms are equal, i.e., $t_g = \frac{m_{\min}(g) + m_{\max}(g)}{2}$. Plugging this into the right-hand side of\ \eqref{eq:UB_step} yields
\[
\min_{t_g}\ \max_{p\in\Delta_g}\ \sum_{l:\,x_{lg}=1} p_l\,\Pi_l(t_g)
\ \le\
y_g\,\max\{b,q\Bigl(\sigma_{\max}(g) + \tfrac12\big(m_{\max}(g)-m_{\min}(g)\big)\Bigr).
\]
This completes the proof.
\end{proof}

\section{Comparison of Different Input Features for Heuristics}\label{app:rho =0.1 feature selection}
We compare the results of different input features for K-Means and K-Medoids heuristics when $\rho=0.1$ in Tables\ \ref{tab: compare input feature kmeans eq14} and\ \ref{tab: compare input feature kmedoids eq14}.
\begin{table}[ht!]
\centering
\caption{Comparison of different input features for K-Means when $\rho = 0.1$}
\resizebox{\textwidth}{!}{
\begin{tabular}{lccccc}
\toprule
{Feature} & {IS} & {OOS} & {CV Time (sec.)} & {Solve Time (sec.)} & {Total Time (sec.)} \\
\midrule
$(m_l,\sigma_l,s_l)$ & $4501.05 $ & $4444.69$ & $1.20$ & $0.05$ & $1.25$ \\
$(m_l)$  & $ 4527.38$ & $4452.32$ & $1.20$ & $0.06$ & $1.26$ \\
$(\sigma_l)$  & $4597.18$& $4499.95$ & $1.19$ & $0.05 $& $1.24$ \\
$(s_l)$   & $4616.57$ &$ 4477.08$ &$ 1.20$ & $0.06$ & $1.26$ \\
$(m_l,\sigma_l)$  & $4524.63 $ & $4459.03$ & $1.20$ & $0.05$ & $1.25$ \\
$(m_l,s_l)$ & $4513.48$ & $\bf{ 4433.79}$ & $1.22$ & $0.05$ & $1.27$ \\
$(\sigma_l,s_l)$  & $4513.57$ & $4466.95$ & $1.20$ & $0.05$ & $1.25$ \\
\bottomrule
\end{tabular}}
\label{tab: compare input feature kmeans eq14}
\end{table}

\begin{table}[ht!]
\centering
\caption{Comparison of different input features for K-Medoids when $\rho = 0.1$}
\resizebox{\textwidth}{!}{
\begin{tabular}{lccccc}
\toprule
{Feature} & {IS} & {OOS} & {CV Time (sec.)} & {Solve Time (sec.)} & {Total Time (sec.)} \\
\midrule
$(m_l,\sigma_l,s_l)$ & $ 4692.99$ & $\bf{4485.05}$ & $0.83$ & $0.04$ & $0.87$ \\
$(m_l)$  & $4625.24$ & $4493.93$ & $0.83$ & $0.04$ & $0.87$ \\
$(\sigma_l)$  & $4655.13$& $4520.65$ & $0.83$ & $0.03 $& $0.86$ \\
$(s_l)$   & $4712.51$ &$ 4513.85$ &$ 0.84$ & $0.04$ & $0.88$ \\
$(m_l,\sigma_l)$  & $4668.22$ & $4495.23$ & $0.85$ & $0.04$ & $0.89$ \\
$(m_l,s_l)$ & $4639.48$ & $4498.57$ & $0.83$ & $0.04$ & $0.87$ \\
$(\sigma_l,s_l)$  & $ 4656.40$ & $ 4517.69$ & $0.83$ & $0.05$ & $0.88$ \\
\bottomrule
\end{tabular}}
\label{tab: compare input feature kmedoids eq14}
\end{table}

\section{Results from different $\rho$ for one run}\label{app:one run rho}

We investigate how the optimal solutions and optimal cost change for different $\rho$ in a single run in Table\ \ref{tab:seed0 rho comparison}.

\begin{table}[ht!]
\centering
\caption{Results across different $\rho$ values with a fixed seed.}
\begin{tabular}{lcccc}
\toprule
$\rho$ & Optimal duration (min.) & IS & OOS & Runtime (sec.)\\
\midrule
$0$   & $164(1),\ 516(2,3,4,5)$            & $\bf{5112.99}$ & $6125.17$      & $0.96$ \\
$0.1$ & $164(1),\ 512(2,3,4,5)$            & $5155.85$ & $6101.67$      & $1.62$ \\
$0.5$ & $164(1),\ 466(2,3,4),\ 541(4)$     & $5202.27$ & $5966.09$   & $1.01$ \\
$1.0$ & $164(1),\ 455(2,3,4),\ 541(4)$     & $5214.00$ & $\bf{5916.96}$      & $0.97$ \\
\bottomrule
\end{tabular}

\label{tab:seed0 rho comparison}
\end{table}

Table\ \ref{tab:seed0 rho comparison} shows that although sometimes the assignments $x_{lg},\ y_g$ remain unchanged when $\rho$ increases - for example, both $\rho=0.5$ and $\rho=1.0$ keep the same grouping assignments, the corresponding optimal treatment time $t_g^*$ and IS cost still shift. This is because the grouping decision variables $x_{lg},\ y_g$ are discrete choices that govern the trade-off between the fixed assignment cost of opening another group and the variability among the modes in that group. Therefore, a small increase in $\rho$ could not be enough to change the grouping. However, the increase in the size of the ambiguity set allows re-weighting the modes, assigning more probability mass to those with higher costs to obtain the worst-case cost so that we can minimize. As a result, different results can be achieved even if the groupings are the same. The run time comparison is inconclusive since adding more mode probabilities as $\rho$ increases changes the problem constraints and doesn't increase monotonically.

\section{MMA-$\rm{DRO^{mis}}$ Problem}\label{app:mmadro_mis}
We relax the ambiguity set $U_l$ in\ \eqref{eq:ul ambiguity moment constrs} to allow a wider range of distributions and thus a more robust model as follows
\begin{align}\label{eq:Ul perturbed ambiguity}
\Pi_l(t_g) = \quad & \max_{\pi_{lk}} \sum_{k=1}^K \pi_{lk} \left(q(\xi_{lk}-t_g)^+ + b(t_g-\xi_{lk})^+\right) \nonumber\\
    \textrm{s.t.} \quad & \sum_{k=1}^K \pi_{lk}\,=1,\nonumber\\
    \quad & \sum_{k=1}^K \pi_{lk} \xi_{lk} \, \le m_l (1+\delta),\nonumber \\
    \quad & \sum_{k=1}^K \pi_{lk} \xi_{lk} \, \ge m_l (1-\delta),\nonumber \\
    \quad & \sum_{k=1}^K \pi_{lk} (\xi_{lk}-m_l)^2 \le \sigma_l^2 (1+\delta),\nonumber\\
    \quad & \sum_{k=1}^K \pi_{lk} (\xi_{lk}-m_l)^2 \ge \sigma_l^2 (1-\delta),\nonumber\\
    \quad & \frac{1}{\sigma_l^2} \left(\sum_{k=1}^K \pi_{lk} (\xi_{lk}-m_l)^{+2} -\sum_{k=1}^K \pi_{lk} (m_l - \xi_{lk})^{+2} \right) \le s_l (1+\delta) , \nonumber \\
    \quad & \frac{1}{\sigma_l^2} \left(\sum_{k=1}^K \pi_{lk} (\xi_{lk}-m_l)^{+2} -\sum_{k=1}^K \pi_{lk} (m_l - \xi_{lk})^{+2} \right) \ge s_l (1-\delta) , \nonumber \\
    \quad & \pi_{lk} \geq 0,\ \forall k \in \mathcal{K}.
\end{align}
Because an exact continuous formulation is unavailable, we discretize the distribution on a finite support set $\Xi \in \mathbb{R}^{\mathcal{L} \times \mathcal{K}}$. Similar to Theorem\ \ref{thm: dual of Pi}, we derive the dual formulation and construct the MMA-$\rm{DRO^{mis}}$ model.

Assign dual variables $\gamma_{lg},\ \nu_{lg}^1,\ \nu_{lg}^2,\ \kappa_{lg}^1,\ \kappa_{lg}^2,\ \zeta_{lg}^1,\ \zeta_{lg}^2$ to the constraints in\ \eqref{eq:Ul perturbed ambiguity}, we formulate the MMA-$\rm{DRO^{mis}}$ problem as follows
\begin{subequations}\label{eq:continuous overall min misspecified}
\begin{align}
\min \quad & \Biggl\{\sum_{g \in \mathcal{G}} c_gy_g + \frac{1}{\sum_{g \in \mathcal{G}} y_g}\sum_{g \in \mathcal{G}} \Bigl(\mu_g y_g + \sum_{l=1}^L \lambda_{lg} x_{lg} + \sum_{l=1}^L (\alpha_{lg} - \beta_{lg}) \hat{p}_l x_{lg}+\tau_g (\rho \sum_{l=1}^L x_{lg}\hat{p}_l) \Bigr) \Biggr\} \nonumber\\
    \textrm{s.t.}  \quad & \mu_g + \lambda_{lg} + (\alpha_{lg}-\beta_{lg})(\sum_{l=1}^L x_{lg}\hat{p}_l) \geq \gamma_{lg} + \Bigl((1-\delta)\nu_{lg}^1 - (1-\delta)\nu_{lg}^2\Bigr) m_l +\Bigl((1-\delta)\kappa_{lg}^1 - (1-\delta)\kappa_{lg}^2\Bigr) \sigma_l^2 \nonumber \\
    \quad & \hspace{5.2cm}  +\Bigl((1-\delta)\zeta_{lg}^1 - (1-\delta)\zeta_{lg}^2\Bigr) s_l \sigma_l^2,\ \forall l \in \mathcal{L},\ g\in \mathcal{G}\nonumber\\
    \quad & \gamma_{lg} + (\nu_{lg}^1 - \nu_{lg}^2)\xi_{lk} + (\kappa_{lg}^1 - \kappa_{lg}^2)(\xi_{lk}-m_l)^2 + \frac{1}{\sigma^2}(\zeta_{lg}^1 - \zeta_{lg}^2) \left((\xi_{lk}-m_l)^{+2} - (m_l-\xi_{lk})^{+2}\right) \nonumber \\
    \quad & \hspace{5.2cm} \geq q(\tilde{\xi}-t_g)^+ + b(t_g-\tilde{\xi})^+,\ \forall l \in \mathcal{L},\ \forall g \in \mathcal{G},\ \forall k \in \mathcal{K} \nonumber\\
    \quad & -\alpha_{lg} - \beta_{lg} + \tau_g \geq 0,\ \forall l \in \mathcal{L},\ g\in \mathcal{G} \nonumber\\
    \quad & \text{Constraints}\ \eqref{eq:sum to one},\ \eqref{eq:x less than y},\ \eqref{eq:feasibility},\ \eqref{eq:t range} \nonumber \\
    \quad & \nu_{lg}^1,\ \nu_{lg}^2,\ \kappa_{lg}^1,\ \kappa_{lg}^2,\ \zeta_{lg}^1,\ \zeta_{lg}^2,\ \lambda_{lg},\  \alpha_{lg},\ \beta_{lg},\ \tau_g,\ t_g \geq 0\nonumber\\
    \quad & \gamma_l,\ \mu_g \in \mathbb{R} \nonumber\\
    \quad & x_{lg},\ y_g \in \{0,1\},\ \forall l \in \mathcal{L},\ g \in \mathcal{G} \nonumber
\end{align}
\end{subequations}

\end{document}